# A PARTIALLY LINEAR FRAMEWORK FOR MASSIVE HETEROGENEOUS DATA

By Tianqi Zhao[§], Guang Cheng[††] and Han Liu[**]

[§**]*Princeton University and* [††]*Purdue University*

We consider a partially linear framework for modelling massive heterogeneous data. The major goal is to extract common features across all sub-populations while exploring heterogeneity of each sub-population. In particular, we propose an aggregation type estimator for the commonality parameter that possesses the (non-asymptotic) minimax optimal bound and asymptotic distribution as if there were no heterogeneity. This oracle result holds when the number of sub-populations does not grow too fast. A plug-in estimator for the heterogeneity parameter is further constructed, and shown to possess the asymptotic distribution as if the commonality information were available. We also test the heterogeneity among a large number of sub-populations. All the above results require to regularize each sub-estimation as though it had the entire sample. Our general theory applies to the divide-and-conquer approach that is often used to deal with massive homogeneous data. A technical by-product of this paper is statistical inferences for general kernel ridge regression. Thorough numerical results are also provided to back up our theory.

**1. Introduction.** In this paper, we propose a partially linear regression framework for modelling massive heterogeneous data. Let $\big\{(Y_i, \boldsymbol{X}_i, Z_i)\big\}_{i=1}^N$ be samples from an underlying distribution that may change with $N$. We assume that there exist $s$ independent sub-populations, and the data from the $j$th sub-population follow a partially linear model:

$$Y = \boldsymbol{X}^T \boldsymbol{\beta}_0^{(j)} + f_0(Z) + \varepsilon, \tag{1.1}$$

[§]PhD student.

[††]Associate Professor. Research Sponsored by NSF CAREER Award DMS-1151692, DMS-1418042, Simons Fellowship in Mathematics, Office of Naval Research (ONR N00014-15-1-2331) and a grant from Indiana Clinical and Translational Sciences Institute. Guang Cheng was on sabbatical at Princeton while part of this work was carried out; he would like to thank the Princeton ORFE department for its hospitality.

[**]Assistant Professor. Research Sponsored by NSF IIS1408910, NSF IIS1332109, NIH R01MH102339, NIH R01GM083084, and NIH R01HG06841

where $\varepsilon$ has zero mean and known variance $\sigma^2$. In the above model, $Y$ depends on $\boldsymbol{X}$ through a linear function that may vary across all sub-populations, and depends on $Z$ through a nonlinear function that is common to all sub-populations. The possibly different values of $\boldsymbol{\beta}_0^{(j)}$ are viewed as the source of heterogeneity. In reality, the number of sub-populations grows and some sub-populations may have extremely high sample sizes. Note that (1.1) is a typical "semi-nonparametric" model (Cheng and Shang, 2015) since we infer both commonality and heterogeneity components throughout the paper.

The model (1.1) is motivated by the following scenario: different labs conduct the same experiment on the relationship between a response variable $Y$ (e.g., heart disease) and a set of predictors $Z, X_1, X_2, \ldots, X_p$. It is known from biological knowledge that the dependence structure between $Y$ and $Z$ (e.g., blood pressure) should be homogeneous for all human. However, for the other covariates (e.g., certain genes), we allow their (linear) relations with $Y$ to potentially vary in different labs. For example, the genetic functionality of different races might be heterogenous. The linear relation is assumed here for simplicity, and particularly suitable when the covariates are discrete.

Statistical modelling for massive data has attracted a flurry of recent research. For homogeneous data, the statistical studies of the divide-and-conquer method currently focus on either parametric inferences, e.g., Li et al. (2013), Bag of Little Bootstraps (Kleiner et al., 2012), and parallel MCMC computing (Wang and Dunson, 2013), or nonparametric minimaxity (Zhang et al., 2013). The other relevant work includes high dimensional linear models with variable selection (Chen and Xie, 2012) and structured perceptron (McDonald et al., 2010). Heterogenous data are often handled by fitting mixture models (Aitkin and Rubin, 1985; McLachlan and Peel, 2004; Figueiredo and Jain, 2002), time varying coefficient models (Hastie and Tibshirani, 1993; Fan and Zhang, 1999) or multitask regression (Huang and Zhang, 2010; Nardi and Rinaldo, 2008; Obozinski et al., 2008). The recent high dimensional work includes Städler et al. (2010); Meinshausen and Bühlmann (2014). However, as far as we are aware, *semi-nonparametric inference* for massive homogeneous/heterogeneous data still remains untouched.

In this paper, our primary goal is to extract common features across all sub-populations while exploring heterogeneity of each sub-population. Specifically, we employ a simple aggregation procedure, which averages commonality estimators across all sub-populations, and then construct a plug-in estimator for each heterogeneity parameter based on the combined estimator for commonality. A similar two-stage estimation method was proposed in Li and Liang (2008), but for the purpose of variable selection in $\boldsymbol{\beta}$ based on a single data set. The secondary goal is to apply the divide-and-conquer method to the sub-population with a huge sample size that is unable to be processed in one single computer. The above purposes are achieved by estimating our statistical model (1.1) with the kernel ridge regression (KRR) method. In the partially linear literature, there also exist other estimation and inference methods (based on a single dataset) such as profile least squares method, partial residual method and backfitting method; see Härdle et al. (2000); Ruppert et al. (2003); Yatchew (2003).

The KRR framework is known to be very flexible and well supported by the general reproducing kernel Hilbert space (RKHS) theory (Mendelson, 2002; Steinwart et al., 2009; Zhang, 2005). In particular, partial smoothing spline models (Cheng et al., 2015) can be viewed as a special case of our general framework. An important technical contribution of this paper is statistical inferences for general kernel ridge regression by extending smoothing spline results developed in Cheng and Shang (2015). This theoretical innovation makes our work go beyond the existing statistical study on the KRR for large datastes, which mainly focus on their nonparametric minimaxity, e.g., Zhang et al. (2013); Bach (2012); Raskutti et al. (2014).

Our theoretical studies are mostly concerned with the so-called "oracle rule" for massive data. Specifically, we define the "oracle estimate" for commonality (heterogeneity) as the one computed when all the heterogeneity information are given (the commonality information is given in each-subpopulation), i.e., $\boldsymbol{\beta}_0^{(j)}$'s are known ($f_0$ is known). We claim that a commonality estimator satisfies the oracle rule if it possesses the same minimax optimality and asymptotic distribution as the "oracle estimate" defined above. A major contribution of this paper is to derive the largest possible diverging rate of $s$ under which our combined estimator for commonality satisfies the oracle rule. In other words, our aggregation procedure is shown to "filter out" the heterogeneity in data when $s$ does not grow too fast with $N$. On the other hand, we have to set a lower bound on $s$ for our heterogeneity estimate to possess the asymptotic distribution as if the commonality information were available, i.e., oracle rule. Our second contribution is to test the heterogeneity among a large number of sub-populations by employing a recent Gaussian approximation theory (Chernozhukov et al., 2013).

In the standard implementation of KRR, we must invert a kernel matrix, which requires costs $O(N^3)$ in time and $O(N^2)$ in memory, respectively; see Saunders et al. (1998). This is computationally prohibitive for a large $N$. Hence, when some sub-population has a huge sample size, we may apply the divide-and-conquer approach whose statistical analysis directly follows from the above results. In this case, the "oracle estimate" is defined as those computed based on the entire (homogeneous) data in those sub-populations.

A rather different goal here is to explore the most computationally efficient way to split the whole sample while performing the best possible statistical inference. Specifically, we derive the largest possible number of splits under which the averaged estimators for both components enjoy the same statistical properties as the oracle estimators.

In both homogeneous and heterogeneous settings above, we note that the upper bounds established for $s$ increase with the smoothness of $f_0$. Hence, our aggregation procedure favors smoother regression functions in the sense that more sub-populations/splits are allowed in the massive data. On the other hand, we have to admit that our upper and lower bound results for $s$ are only sufficient conditions although empirical results show that our bounds are quite sharp. Another interesting finding is that even the semi-nonparametric estimation is applied to only one fraction of the entire data, it is nonetheless essential to regularize each sub-estimation as if it had the entire sample.

In the end, we highlight two key technical challenges: (i) delicate interplay between the parametric and nonparametric components in the *semi-nonparametric estimation*. In particular, we observe a "bias propagation" phenomenon: the bias introduced by the penalization of the nonparametric component propagates to the parametric component, and the resulting parametric bias in turn propagates back to the nonparametric component. To analyze this complicated propagation mechanism, we extend the existing RKHS theory to an enlarged partially linear function space by defining a novel inner product under which the expectation of the Hessian of the objective function becomes identity; see Proposition 2.2. (ii) double asymptotics framework in terms of diverging $s$ and $N$. In this challenging regime, we develop more refined concentration inequalities in characterizing the concentration property of an averaged empirical process. These very refined theoretical analyses show that an average of $s$ asymptotic linear expansions is still a valid one as $s \wedge N \to \infty$.

The rest of the paper is organized as follows: Section 2 briefly introduces the general RKHS theory and discusses its extension to an enlarged partially linear function space. Section 3 describes our aggregation procedure, and studies the "oracle" property of this procedure from both asymptotic and non-asymptotic perspectives. The efficiency boosting of heterogeneity estimators and heterogenous testing results are also presented in this section. Section 4 applies our general theory to various examples with different smoothness. Section 5 is devoted to the analysis of divide-and-conquer algorithms for homogeneous data. Section 6 presents some numerical experiments. All the technical details are deferred to Section 7 or Online Supplementary (Zhao et al., 2015).

**Notation:** Denote $\|\cdot\|_2$ and $\|\cdot\|_\infty$ as the Euclidean $L_2$ and infinity norm in $\mathbb{R}^p$, respectively. For any function $f : \mathcal{S} \mapsto \mathbb{R}$, let $\|f\|_{\sup} = \sup_{x\in\mathcal{S}} |f(x)|$. We use $\|\cdot\|$ to denote the spectral norm of matrices. For positive sequences $a_n$ and $b_n$, we write $a_n \lesssim b_n$ ($a_n \gtrsim b_n$) if there exists some universal constant constant $c > 0$ ($c' > 0$) independent of $n$ such that $a_n \le cb_n$ ($a_n \ge c'b_n$) for all $n \in \mathbb{N}$. We denote $a_n \asymp b_n$ if both $a_n \lesssim b_n$ and $a_n \gtrsim b_n$.

## 2. Preliminaries.

In this section, we briefly introduce the general RKHS theory, and then extend it to a partially linear function space. Below is a generic definition of RKHS (Berlinet and Thomas-Agnan, 2004):

**Definition 2.1.** Denote by $\mathcal{F}(\mathcal{S}, \mathbb{R})$ a vector space of functions from a general set $\mathcal{S}$ to $\mathbb{R}$. We say that $\mathcal{H}$ is a reproducing kernel Hilbert space (RKHS) on $\mathcal{S}$, provided that:

(i) $\mathcal{H}$ is a vector subspace of $\mathcal{F}(\mathcal{S}, \mathbb{R})$;
(ii) $\mathcal{H}$ is endowed with an inner product, denoted as $\langle\cdot,\cdot\rangle_{\mathcal{H}}$, under which it becomes a Hilbert space;
(iii) for every $y \in \mathcal{S}$, the linear evaluation functional defined by $E_y(f) = f(y)$ is bounded.

If $\mathcal{H}$ is a RKHS, by Riesz representation, we have that for every $y \in \mathcal{S}$, there exists a unique vector, $K_y \in \mathcal{H}$, such that for every $f \in \mathcal{H}$, $f(y) = \langle f, K_y\rangle_{\mathcal{H}}$. The reproducing kernel for $\mathcal{H}$ is defined as $K(x, y) = K_y(x)$.

Denote $U := (\boldsymbol{X}, Z) \in \mathcal{X} \times \mathcal{Z} \subset \mathbb{R}^p \times \mathbb{R}$, and $\mathbb{P}_U$ as the distribution of $U$ ($\mathbb{P}_X$ and $\mathbb{P}_Z$ are defined similarly). According to Definition 2.1, if $\mathcal{S} = \mathcal{Z}$ and $\mathcal{F}(\mathcal{Z}, \mathbb{R}) = L_2(\mathbb{P}_Z)$, then we can define a RKHS $\mathcal{H} \subset L_2(\mathbb{P}_Z)$ (endowed with a proper inner product $\langle\cdot,\cdot\rangle_{\mathcal{H}}$), in which the true function $f_0$ is believed to lie. The corresponding kernel for $\mathcal{H}$ is denoted by $K$ such that for any $z \in \mathcal{Z}$: $f(z) = \langle f, K_z\rangle_{\mathcal{H}}$. By Mercer theorem, this kernel function has the following unique eigen-decomposition:

$$K(z_1, z_2) = \sum_{\ell=1}^{\infty} \mu_\ell \phi_\ell(z_1)\phi_\ell(z_2),$$

where $\mu_1 \ge \mu_2 \ge \ldots > 0$ are eigenvalues and $\{\phi_\ell\}_{\ell=1}^\infty$ are an orthonormal basis in $L_2(\mathbb{P}_Z)$. Mercer theorem together with the reproducing property implies that $\langle\phi_i, \phi_j\rangle_{\mathcal{H}} = \delta_{ij}/\mu_i$, where $\delta_{ij}$ is the Kronecker's delta. The smoothness of the functions in RKHS can be characterized by the decaying rate of $\{\mu_\ell\}_{\ell=1}^\infty$. Below, we present three different decaying rates together with the corresponding kernel functions.

**Finite rank kernel:** the kernel has finite rank $r$ if $\mu_\ell = 0$ for all $\ell > r$. For example, the linear kernel $K(\boldsymbol{z}_1, \boldsymbol{z}_2) = \langle \boldsymbol{z}_1, \boldsymbol{z}_2\rangle_{\mathbb{R}^d}$ has rank $d$, and generates a $d$-dimensional linear function space. The eigenfunctions are given by $\phi_\ell(\boldsymbol{z}) = z_\ell$ for $\ell = 1, \ldots, d$. The polynomial kernel $K(z_1, z_2) = (1 + z_1 z_2)^d$ has rank $d+1$, and generates a space of polynomial functions with degree at most $d$. The eigenfunctions are given by $\phi_\ell = z^{\ell-1}$ for $\ell = 1, \ldots, d+1$.

**Exponentially decaying kernel:** the kernel has eigenvalues that satisfy $\mu_\ell \asymp c_1 \exp(-c_2 \ell^p)$ for some $c_1, c_2 > 0$. An example is the Gaussian kernel $K(z_1, z_2) = \exp(-|z_1 - z_2|^2)$. The eigenfunctions are given by Sollich and Williams (2005)

$$\phi_\ell(x) = (\sqrt{5}/4)^{1/4} (2^{\ell-2}(\ell-1)!)^{-1/2} e^{-(\sqrt{5}-1)x^2/4} H_{\ell-1}\big((\sqrt{5}/2)^{1/2} x\big), \tag{2.1}$$

for $\ell = 1, 2, \cdots$, where $H_\ell(\cdot)$ is the $\ell$-th Hermite polynomial.

**Polynomially decaying kernel:** the kernel has eigenvalues that satisfy $\mu_\ell \asymp \ell^{-2\nu}$ for some $\nu \geq 1/2$. Examples include those underlying for Sobolev space and Besov space (Birman and Solomjak, 1967). In particular, the eigenfunctions of a $\nu$-th order periodic Sobolev space are trigonometric functions as specified in Section 4.3. The corresponding Sobolev kernels are given in Gu (2013).

In this paper, we consider the following penalized estimation:

$$(\widehat{\boldsymbol{\beta}}^\dagger, \widehat{f}^\dagger) = \operatorname*{argmin}_{(\boldsymbol{\beta}, f) \in \mathcal{A}} \left\{ \frac{1}{n} \sum_{i=1}^{n} \big(Y_i - \boldsymbol{X}_i^T \boldsymbol{\beta} - f(Z_i)\big)^2 + \lambda \|f\|_{\mathcal{H}}^2 \right\}, \tag{2.2}$$

where $\lambda > 0$ is a regularization parameter and $\mathcal{A}$ is defined as the parameter space $\mathbb{R}^p \times \mathcal{H}$. For simplicity, we do not distinguish $m = (\boldsymbol{\beta}, f) \in \mathcal{A}$ from its associated function

$$m \in \mathcal{M} := \big\{ m \mid m(u) = \boldsymbol{\beta}^T \boldsymbol{x} + f(z), \text{ for } u = (\boldsymbol{x}, z) \in \mathcal{X} \times \mathcal{Z}, \ (\boldsymbol{\beta}, f) \in \mathcal{A} \big\},$$

throughout the paper. We call $(\widehat{\boldsymbol{\beta}}^\dagger, \widehat{f}^\dagger)$ as partially linear kernel ridge regression (KRR) estimate in comparison with the nonparametric KRR estimate in Shawe-Taylor and Cristianini (2004). In particular, when $\mathcal{H}$ is a $\nu$-th order Sobolev space endowed with $\langle f, \widetilde{f} \rangle_{\mathcal{H}} := \int_{\mathcal{Z}} f^{(\nu)}(z) \widetilde{f}^{(\nu)}(z) dz$, $(\widehat{\boldsymbol{\beta}}^\dagger, \widehat{f}^\dagger)$ becomes the commonly used partial smoothing spline estimate.

We next illustrate that $\mathcal{A}$ can be viewed as a partially linear extension of $\mathcal{H}$ in the sense that it shares some nice reproducing properties as this RKHS $\mathcal{H}$ under the following inner product: for any $m = (\boldsymbol{\beta}, f) \in \mathcal{A}$ and $\widetilde{m} = (\widetilde{\boldsymbol{\beta}}, \widetilde{f}) \in \mathcal{A}$, define

$$\langle m, \widetilde{m} \rangle_{\mathcal{A}} := \langle m, \widetilde{m} \rangle_{L_2(\mathbb{P}_{\boldsymbol{X}, Z})} + \lambda \langle f, \widetilde{f} \rangle_{\mathcal{H}}, \tag{2.3}$$

where $\langle m, \widetilde{m} \rangle_{L_2(\mathbb{P}_{\boldsymbol{X}, Z})} = \mathbb{E}_{\boldsymbol{X}, Z} \big[ (\boldsymbol{X}^T \boldsymbol{\beta} + f(Z))(\boldsymbol{X}^T \widetilde{\boldsymbol{\beta}} + \widetilde{f}(Z)) \big]$. Note that $m$ and $\widetilde{m}$ in $\langle m, \widetilde{m} \rangle_{L_2(\mathbb{P}_{\boldsymbol{X}, Z})}$ are both functions in the set $\mathcal{M}$. Similar as the kernel function $K_z$, we can construct a linear operator $R_u(\cdot) \in \mathcal{A}$ such that $\langle R_u, m\rangle_{\mathcal{A}} = m(u)$ for any $u \in \mathcal{X} \times \mathcal{Z}$. Also, construct another linear operator $P_\lambda : \mathcal{A} \mapsto \mathcal{A}$ such that $\langle P_\lambda m, \widetilde{m}\rangle_{\mathcal{A}} = \lambda\langle f, \widetilde{f}\rangle_{\mathcal{H}}$ for any $m$ and $\widetilde{m}$. See Proposition 2.3 for the construction of $R_u$ and $P_\lambda$.

We next present a proposition illustrating the rationale behind the definition of $\langle\cdot,\cdot\rangle_{\mathcal{A}}$. Denote $\otimes$ as the outer product on $\mathcal{A}$. Hence, $\mathbb{E}_U[R_U \otimes R_U] + P_\lambda$ is an operator from $\mathcal{A}$ to $\mathcal{A}$.

**Proposition 2.2.** $\mathbb{E}_U[R_U \otimes R_U] + P_\lambda = id$, where $id$ is an identity operator on $\mathcal{A}$.

Proof. For any $m = (\boldsymbol{\beta}, f) \in \mathcal{A}$ and $\widetilde{m} = (\widetilde{\boldsymbol{\beta}}, \widetilde{f}) \in \mathcal{A}$, we have

$$\begin{aligned}\langle(\mathbb{E}_U[R_U \otimes R_U] + P_\lambda)m, \widetilde{m}\rangle_{\mathcal{A}} &= \langle\mathbb{E}_U[R_U \otimes R_U]m, \widetilde{m}\rangle_{\mathcal{A}} + \langle P_\lambda m, \widetilde{m}\rangle_{\mathcal{A}} \\ &= \mathbb{E}_U[m(U)\widetilde{m}(U)] + \lambda\langle f, \widetilde{f}\rangle_{\mathcal{H}} \\ &= \langle m, \widetilde{m}\rangle_{\mathcal{A}}.\end{aligned}$$

Since the choice of $(m, \widetilde{m})$ is arbitrary, we conclude our proof. □

As will be seen in the subsequent analysis, e.g., in Theorem 3.4, the operator $\mathbb{E}[R_U \otimes R_U] + P_\lambda$ is essentially the expectation of the Hessian of the objective function (w.r.t. Fréchet derivative) minimized in (2.2). Proposition 2.2 shows that the inversion of this Hessian matrix is trivial when the inner product is designed as in (2.3). Due to that, the theoretical analysis of $\widehat{m}^\dagger = (\widehat{\boldsymbol{\beta}}^\dagger, \widehat{f}^\dagger)$ based on the first order optimality condition becomes much more transparent.

To facilitate the construction of $R_u$ and $P_\lambda$, we need to endow a new inner product with $\mathcal{H}$:

$$\langle f, \widetilde{f}\rangle_{\mathcal{C}} = \langle f, \widetilde{f}\rangle_{L_2(\mathbb{P}_Z)} + \lambda\langle f, \widetilde{f}\rangle_{\mathcal{H}}, \tag{2.4}$$

for any $f, \widetilde{f} \in \mathcal{H}$. Under (2.4), $\mathcal{H}$ is still a RKHS as the evaluation functional is bounded by Lemma A.1. We denote the new kernel function as $\widetilde{K}(\cdot,\cdot)$, and define a positive definite self-adjoint operator $W_\lambda : \mathcal{H} \mapsto \mathcal{H}$:

$$\langle W_\lambda f, \widetilde{f}\rangle_{\mathcal{C}} = \lambda\langle f, \widetilde{f}\rangle_{\mathcal{H}} \text{ for any } f, \widetilde{f} \in \mathcal{H}', \tag{2.5}$$

whose existence is proven in Lemma A.2. We next define two crucial quantities needed in the construction: $B_k := \mathbb{E}[X_k \mid Z]$ and its Riesz representer $A_k \in \mathcal{H}$ satisfying $\langle A_k, f\rangle_{\mathcal{C}} = \langle B_k, f\rangle_{L_2(\mathbb{P}_Z)}$ for all $f \in \mathcal{H}$. Here, we implicity assume $B_k$ is square integrable. The existence of $A_k$ follows from the boundedness of the linear functional $\mathcal{B}_k f := \langle B_k, f\rangle_{L_2(\mathbb{P}_Z)}$ (by Riesz's representer theorem) as follows:

$$|\mathcal{B}_k f| = |\langle B_k, f\rangle_{L_2(\mathbb{P}_Z)}| \le \|B_k\|_{L_2(\mathbb{P}_Z)}\|f\|_{L_2(\mathbb{P}_Z)} \le \|B_k\|_{L_2(\mathbb{P}_Z)}\|f\|_{\mathcal{C}}.$$

We are now ready to construct $R_u$ and $P_\lambda$ based on $\widetilde{K}_z$, $W_\lambda$, $\boldsymbol{B}$ and $\boldsymbol{A}$

introduced above, where $\boldsymbol{B} = (B_1, \ldots, B_p)^T$ and $\boldsymbol{A} = (A_1, \ldots, A_p)^T$. Define $\boldsymbol{\Omega} = \mathbb{E}\big[(\boldsymbol{X} - \boldsymbol{B})(\boldsymbol{X} - \boldsymbol{B})^T\big]$ and $\boldsymbol{\Sigma}_\lambda = \mathbb{E}[\boldsymbol{B}(Z)(\boldsymbol{B}(Z) - \boldsymbol{A}(Z))^T]$.

**Proposition 2.3.** For any $u = (\boldsymbol{x}, z)$, $R_u$ can be expressed as $R_u : u \mapsto (L_u, N_u) \in \mathcal{A}$, where

$$L_u = (\boldsymbol{\Omega} + \boldsymbol{\Sigma}_\lambda)^{-1}(\boldsymbol{x} - \boldsymbol{A}(z)) \text{ and } N_u = \widetilde{K}_z - \boldsymbol{A}^T L_u,$$

Moreover, for any $m = (\boldsymbol{\beta}, f) \in \mathcal{A}$, $P_\lambda m$ can be expressed as $P_\lambda m = (L_\lambda f, N_\lambda f) \in \mathcal{A}$, where

$$L_\lambda f = -(\boldsymbol{\Omega} + \boldsymbol{\Sigma}_\lambda)^{-1}\langle \boldsymbol{B}, W_\lambda f\rangle_{L_2(\mathbb{P}_Z)} \text{ and } N_\lambda f = W_\lambda f - \boldsymbol{A}^T L_\lambda f.$$

The quantities $R_u$ and $P_\lambda$ correspond to the variance, i.e., $n^{-1}\sum_{i=1}^n R_{U_i}\varepsilon_i$, and bias, i.e., $P_\lambda m_0$, in the stochastic expansion of $\widehat{m}^\dagger - m_0$, where $\widehat{m}^\dagger = (\widehat{\boldsymbol{\beta}}^\dagger, \widehat{f}^\dagger)$, $m_0 = (\boldsymbol{\beta}_0, f_0)$; see Eq. (7.2) in Section 7.1. We remark that the penalized loss function in (2.2) can be written as $(1/n)\sum_{i=1}^n \big(Y_i - \langle R_{U_i}, m\rangle_{\mathcal{A}}\big)^2 + \langle P_\lambda m, m\rangle_{\mathcal{A}}$. This explains why $R_u$ and $P_\lambda$ show up in the stochastic expansion, which is derived from the KKT condition of the above loss function and Proposition 2.2. Moreover, $R_u = (L_u, N_u)$ and $P_\lambda m = (L_\lambda f, N_\lambda f)$. Hence, $L_u, L_\lambda f_0$ and $N_u, N_\lambda f_0$ appear in the stochastic expansions of $\widehat{\boldsymbol{\beta}}^\dagger - \boldsymbol{\beta}_0$ and $\widehat{f}^\dagger - f_0$; see Lemma 3.1.

**3. Heterogeneous Data: Aggregation of Commonality.** In this section, we start from describing our aggregation procedure and model assumptions in Section 3.1. The main theoretical results are presented in Sections 3.2 – 3.4 showing that our combined estimate for commonality enjoys the "oracle property". To be more specific, we show that it possesses the same (non-asymptotic) minimax optimal bound (in terms of mean-squared error) and asymptotic distribution as the "oracle estimate" $\widehat{f}_{or}$ computed when all the heterogeneity information are available:

$$\widehat{f}_{or} = \underset{f\in\mathcal{H}}{\operatorname{argmin}} \left\{ \frac{1}{N}\sum_{j=1}^{s}\sum_{i\in S_j}\big(Y_i - (\boldsymbol{\beta}_0^{(j)})^T\boldsymbol{X}_i - f(Z_i)\big)^2 + \lambda\|f\|_{\mathcal{H}}^2 \right\}, \tag{3.1}$$

where $S_j$ denotes the index set of all samples from the subpopulation $j$. The above nice properties hold when the number of sub-populations does not grow too fast and the smoothing parameter is chosen according to the entire sample size $N$. Based on this combined estimator, we further construct a plug-in estimator for each heterogeneity parameter $\boldsymbol{\beta}_0^{(j)}$, which possesses the asymptotic distribution as if the commonality were known, in Section 3.5. Interestingly, this oracular result holds when the number of sub-population is not too small. In the end, Section 3.6 tests the possible heterogeneity among a large number of sub-populations.

3.1. *Method and Assumptions.* The heterogeneous data setup and averaging procedure are described below:

1. Observe data $(\boldsymbol{X}_i, Z_i, Y_i)$ with the known labels indicating the sub-population it belongs to, for $i = 1, \ldots, N$. The size of samples from each sub-population is assumed to be the same, denoted by $n$, for simplicity. Hence, $N = n \times s$.
2. On the $j$-th sub-population, obtain the following penalized estimator:

$$(\widehat{\boldsymbol{\beta}}_{n,\lambda}^{(j)}, \widehat{f}_{n,\lambda}^{(j)}) = \underset{(\boldsymbol{\beta}, f)\in\mathcal{A}}{\operatorname{argmin}} \left\{ \frac{1}{n} \sum_{i \in S_j} \big(Y_i - \boldsymbol{X}_i^T \boldsymbol{\beta} - f(Z_i)\big)^2 + \lambda \|f\|_{\mathcal{H}}^2 \right\}. \tag{3.2}$$

3. Obtain the final nonparametric estimate[1] for commonality by averaging:

$$\bar{f}_{N,\lambda} = \frac{1}{s} \sum_{j=1}^{s} \widehat{f}_{n,\lambda}^{(j)}. \tag{3.3}$$

We point out that $\widehat{\boldsymbol{\beta}}_{n,\lambda}^{(j)}$ is not our final estimate for heterogeneity. In fact, it can be further improved based on $\bar{f}_{N,\lambda}$; see Section 3.5.

For simplicity, we will drop the subscripts $(n, \lambda)$ and $(N, \lambda)$ in those notation defined in (3.2) and (3.3) throughout the rest of this paper. Moreover, we make the technical assumption that $s \lesssim N^{\psi}$, although $\psi$ could be very close to 1. The main assumptions of this section are stated below.

**Assumption 3.1** (Regularity Condition)**.** (i) $\varepsilon_i$'s are i.i.d. sub-Gaussian random variables independent of the designs; (ii) $B_k \in L_2(\mathbb{P}_Z)$ for all $k$, and $\boldsymbol{\Omega} := \mathbb{E}[(\boldsymbol{X} - \boldsymbol{B}(Z))(\boldsymbol{X} - \boldsymbol{B}(Z))^T]$ is positive definite; (iii) $\boldsymbol{X}_i$'s are uniformly bounded by a constant $c_x$.

Conditions in Assumption 3.1 are fairly standard in the literature. For example, the positive definiteness of $\boldsymbol{\Omega}$ is needed for obtaining semiparametric efficient estimation; see Mammen and van de Geer (1997). Note that we do not require the independence between $\boldsymbol{X}$ and $Z$ throughout the paper.

**Assumption 3.2** (Kernel Condition)**.** We assume that there exist $0 < c_\phi < \infty$ and $0 < c_K < \infty$ such that $\sup_\ell \|\phi_\ell\|_{\sup} \leq c_\phi$ and $\sup_z K(z, z) \leq c_K$.

Assumption 3.2 is commonly assumed in kernel ridge regression literature (Zhang et al., 2013; Lafferty and Lebanon, 2005; Guo, 2002). In the case of finite rank kernel, e.g., linear and polynomial kernels, the eigenfunctions are uniformly bounded as long as $\mathcal{Z}$ has finite support. As for the exponentially decaying kernels such as Gaussian kernel, we prove in Section 4.2 that the eigenfunctions given in (2.1) are uniformly bounded by 1.336. Lastly, for the polynomially decaying kernels, Proposition 2.2 in Shang and Cheng (2013) showed that the eigenfunctions induced from a $\nu$-th order Sobolev space (under a proper inner product $\langle\cdot,\cdot\rangle_{\mathcal{H}}$) are uniformly bounded under mild smoothness conditions for the density of $Z$.

[1] The commonality estimator $\bar{f}_{N,\lambda}$ can be adjusted as a weighted sum $\sum_{j=1}^{s}(n_j/N)\widehat{f}_{n,\lambda}^{(j)}$ if sub-sample sizes are different. In particular, the divide-and-conquer method can be applied to those sub-populations with huge sample sizes; see Section 5.

**Assumption 3.3.** For each $k = 1, \dots, p$, $B_k(\cdot) \in \mathcal{H}$. This is equivalent to

$$\sum_{\ell=1}^{\infty} \mu_\ell^{-1} \langle B_k, \phi_\ell \rangle^2_{L_2(\mathbb{P}_Z)} < \infty.$$

Assumption 3.3 requires the conditional expectation of $X_k$ given $Z = z$ is as smooth as $f_0(z)$. As can be seen in Section 3.4, this condition is imposed to control the bias of the parametric component, which is caused by penalization on the nonparametric component. We call this interaction as the "bias propagation phenomenon", and study it in Section 3.4.

Before laying out our main theoretical results, we define a key quantity used throughout the paper:

$$d(\lambda) := \sum_{\ell=1}^{\infty} \frac{1}{1 + \lambda/\mu_\ell}. \tag{3.4}$$

The quantity $d(\lambda)$ is essentially the "effective dimension", which was introduced in (Zhang, 2005). For a finite dimensional space, $d(\lambda)$ corresponds to the true dimension, e.g., $d(\lambda) \asymp r$ for the finite rank kernel (with rank $r$). For an infinite-dimensional space, $d(\lambda)$ is jointly determined by the size of that space and the smoothing parameter $\lambda$. For example, $d(\lambda) \asymp (-\log \lambda)^{1/p}$ for exponentially decaying kernel, and $d(\lambda) \asymp \lambda^{-1/(2m)}$ for polynomially decaying kernels. More details are provided in Section 4, where the three RKHS examples are carefully discussed.

In the end, we state a technical lemma that is crucially important in the subsequent theoretical derivations. For any function space $\mathcal{F}$, define an entropy integral

$$J(\mathcal{F}, \delta) =: \int_0^\delta \sqrt{\log \mathcal{N}(\mathcal{F}, \|\cdot\|_{\sup}, \epsilon)}\, d\epsilon,$$

where $\mathcal{N}(\mathcal{F}, \|\cdot\|_{\sup}, \epsilon)$ is an $\epsilon$-covering number of $\mathcal{F}$ w.r.t. supreme norm. Define the following sets of functions: $\mathcal{F}_1 = \{f \mid f(\boldsymbol{x}) = \boldsymbol{x}^T \boldsymbol{\beta}$ for $\boldsymbol{x} \in \mathcal{X}, \boldsymbol{\beta} \in \mathbb{R}^p, \|f\|_{\sup} \leq 1\}$, $\mathcal{F}_2 = \{f \in \mathcal{H} : \|f\|_{\sup} \leq 1, \|f\|_{\mathcal{H}} \leq d(\lambda)^{-1/2} \lambda^{-1/2}\}$, $\mathcal{F} := \{f = f_1 + f_2 : f_1 \in \mathcal{F}_1, f_2 \in \mathcal{F}_2, \|f\|_{\sup} \leq 1/2\}$.

**Lemma 3.1.** For any fixed $j = 1, \dots, s$, we have

$$\widehat{\boldsymbol{\beta}}^{(j)} - \boldsymbol{\beta}_0^{(j)} = \frac{1}{n}\sum_{i\in S_j} L_{U_i}\varepsilon_i - L_\lambda f_0 - Rem_\beta^{(j)}, \tag{3.5}$$

and

$$\widehat{f}^{(j)} - f_0 = \frac{1}{n}\sum_{i\in S_j} N_{U_i}\varepsilon_i - N_\lambda f_0 - Rem_f^{(j)}. \tag{3.6}$$

where $(Rem_\beta^{(j)}, Rem_f^{(j)}) \in \mathcal{A}$. Moreover, suppose Assumptions 3.1 and 3.2 hold, and $d(\lambda)n^{-1/2}(J(\mathcal{F},1) + \log n) = o(1)$, then we have

$$\mathbb{E}\big[\|Rem_\beta^{(j)}\|_2^2\big] \le a(n,\lambda,J), \tag{3.7}$$

and

$$\mathbb{P}\Big(\|Rem_\beta^{(j)}\|_2 \ge b(n,\lambda,J)\Big) \lesssim n\exp(-c\log^2 n), \tag{3.8}$$

where $a(n,\lambda,J) = Cd(\lambda)^2 n^{-1} r_{n,\lambda}^2\big(J(\mathcal{F},1)^2+1\big) + Cd(\lambda)^2\lambda^{-1}n\exp(-c\log^2 n)$, $b(n,\lambda,J) = Cd(\lambda)n^{-1/2}r_{n,\lambda}(J(\mathcal{F},1)+\log n)$ and $r_{n,\lambda} = (\log n)^2(d(\lambda)/n)^{1/2} + \lambda^{1/2}$. The same inequalities also hold for $\|Rem_f^{(j)}\|_{\mathcal{C}}$ under the same set of conditions.

Eq. (3.5) and (3.6) in the above lemma are fundamentally important in deriving the subsequent asymptotic and non-asymptotic results. In particular, by (3.6) and the definition of $\bar{f}$, we obtain the stochastic expansion

$$\bar{f} - f_0 = \frac{1}{N}\sum_{i=1}^{N} N_{U_i}\varepsilon_i - N_\lambda f_0 - \frac{1}{s}\sum_{j=1}^{s} Rem_f^{(j)}, \tag{3.9}$$

which is the starting point for deriving the results in Theorems 3.2 and 3.4. These two equations are also of independent interest. For example, they trivially apply to the classical setup where there is only one dataset, i.e., $s = 1$. In addition, they can be used in the other model settings where sub-populations do not share the same sample size or are not independent. As far as we know, the (non-asymptotic) moment and probability inequalities (3.7) and (3.8) on the remainder term are new. They are useful in determining a proper growth rate of $s$ such that the aggregated remainder term $(1/s)\sum_{j=1}^{s} Rem_f^{(j)}$ still vanishes in probability.

3.2. *Non-Asymptotic Bound for Mean-Squared Error.* The primary goal of this section is to evaluate the estimation quality of the combined estimate from a *non-asymptotic* point of view. Specifically, we derive a finite sample upper bound for the mean-squared error $\mathrm{MSE}(\bar{f}) := \mathbb{E}\big[\|\bar{f} - f_0\|_{L_2(\mathbb{P}_Z)}^2\big]$. When $s$ does not grow too fast, we show that $\mathrm{MSE}(\bar{f})$ is of the order $O\big(d(\lambda)/N+\lambda\big)$, from which the aggregation effect on $f$ can be clearly seen. If $\lambda$ is chosen in the order of $N$, the mean-squared error attains the (un-improvable) optimal minimax rate. As a by-product, we establish a *non-asymptotic* upper bound for the mean-squared error of $\widehat{\boldsymbol\beta}^{(j)}$, i.e., $\mathrm{MSE}(\widehat{\boldsymbol\beta}^{(j)}) := \mathbb{E}\big[\|\widehat{\boldsymbol\beta}^{(j)} - \boldsymbol\beta_0^{(j)}\|_2^2\big]$. The results in this section together with Theorem 3.6 in Section 3.4 determine an upper bound of $s$ under which $\bar{f}$ *enjoys the same statistical properties (minimax optimality and asymptotic distribution) as the oracle estimate* $\widehat{f}_{or}$.

Define $\tau_{\min}(\boldsymbol\Omega)$ as the minimum eigenvalue of $\boldsymbol\Omega$ and $\mathrm{Tr}(K) := \sum_{\ell=1}^{\infty}\mu_\ell$ as the trace of $K$. Moreover, let $C_1' = 2\tau_{\min}^{-2}(\boldsymbol\Omega)\Big(c_x^2 p + c_\phi^2\,\mathrm{Tr}(K)\sum_{k=1}^p\|B_k\|_{\mathcal{H}}^2\Big)$, $C_1 = 2c_\phi^2\Big(1 + C_1'\sum_{k=1}^p\|B_k\|_{L_2(\mathbb{P}_Z)}^2\Big)$, $C_2' = \tau_{\min}^{-2}(\boldsymbol\Omega)\|f_0\|_{\mathcal{H}}^2\sum_{k=1}^p\|B_k\|_{\mathcal{H}}^2$ and $C_2 = 2C_2'\sum_{k=1}^p\|B_k\|_{L_2(\mathbb{P}_Z)}^2$.

**Theorem 3.2.** Under Assumptions 3.1 – 3.3, if $s = o\Big(Nd(\lambda)^{-2}\big(J(\mathcal{F},1) + \log N\big)^{-2}\Big)$, then we have

$$\mathrm{MSE}(\bar{f}) \le C_1\sigma^2 d(\lambda)/N + 2\|f_0\|_{\mathcal{H}}^2\lambda + C_2\lambda^2 + s^{-1}a(n,\lambda,J), \tag{3.10}$$

where $a(n,\lambda,J)$ is defined in Lemma 3.1.

Typically, we require an upper bound for $s$ so that the fourth term in the R.H.S. of (3.10) can be dominated by the first two terms, which correspond to variance and bias, respectively. To attain the optimal *bias-variance trade-off*, we choose $\lambda \asymp d(\lambda)/N$. Solving this equation yields the choice of regularization parameter $\lambda$, which varies in different RKHS. The resulting rate of convergence for $\mathrm{MSE}(\bar{f})$ coincides with the minimax optimal rate of the oracle estimate in different RKHS; see Section 4. This can be viewed as a non-asymptotic version of the "oracle property" of $\bar{f}$. In comparison with the nonparametric KRR result in Zhang et al. (2013), we realize that adding one parametric component does not affect the finite sample upper bound (3.10).

As a by-product, we obtain a *non-asymptotic* upper bound for $\mathrm{MSE}(\widehat{\boldsymbol\beta}^{(j)})$. This result is new, and also of independent interest.

**Theorem 3.3.** Suppose that Assumptions 3.1 – 3.3 hold. Then we have

$$\mathrm{MSE}(\widehat{\boldsymbol\beta}^{(j)}) \le C_1'\sigma^2 n^{-1} + C_2'\lambda^2 + a(n,\lambda,J), \tag{3.11}$$

where $a(n,\lambda,J)$ is defined in Lemma 3.1, and $C_1'$ and $C_2'$ are defined before Theorem 3.2.

Again, the first term and second term in the R.H.S. of (3.11) correspond to the variance and bias, respectively. In particular, the second term comes from the bias propagation effect to be discussed in Section 3.4. By choosing $\lambda = o(n^{-1/2})$, we can obtain the optimal rate of $\mathrm{MSE}(\widehat{\boldsymbol{\beta}}^{(j)})$, i.e., $O(n^{-1/2})$, but may lose the minimax optimality of $\mathrm{MSE}(\bar{f})$ in most cases.

3.3. *Joint Asymptotic Distribution.* In this section, we derive a preliminary result on the joint limit distribution of $(\widehat{\boldsymbol{\beta}}^{(j)}, \bar{f}(z_0))$ at any $z_0 \in \mathcal{Z}$. A key issue with this result is that their centering is not at the true value. However, we still choose to present it here since we will observe an interesting phenomenon when removing the bias in Section 3.4.

**Theorem 3.4** (Joint Asymptotics I)**.** Suppose that Assumptions 3.1 and 3.2 hold, and that as $N \to \infty$, $\|\widetilde{K}_{z_0}\|_{L_2(\mathbb{P}_Z)}/d(\lambda)^{1/2} \to \sigma_{z_0}$, $(W_\lambda \boldsymbol{A})(z_0)/d(\lambda)^{1/2} \to \boldsymbol{\alpha}_{z_0} \in \mathbb{R}^p$, and $\boldsymbol{A}(z_0)/d(\lambda)^{1/2} \to -\boldsymbol{\gamma}_{z_0} \in \mathbb{R}^p$. Suppose the following conditions are satisfied

$$s = o\Big(N d(\lambda)^{-2}\big(J(\mathcal{F},1) + \log N\big)^{-2}\Big), \tag{3.12}$$

$$s d(\lambda)/N \log^4 N + \lambda = o\Big(d(\lambda)^{-2}\big(J(\mathcal{F},1) + \log N\big)^{-2} \log^{-2} N\Big). \tag{3.13}$$

Denote $(\boldsymbol{\beta}_0^{(j)*}, f_0^*)$ as $(id - P_\lambda) m_0^{(j)}$, where $m_0^{(j)} = (\boldsymbol{\beta}_0^{(j)}, f_0)$. We have for any $z_0 \in \mathcal{Z}$ and $j = 1, \ldots, s$,

(i) if $s \to \infty$ then

$$\begin{pmatrix} \sqrt{n}(\widehat{\boldsymbol{\beta}}^{(j)} - \boldsymbol{\beta}_0^{(j)*}) \\ \sqrt{N/d(\lambda)}\big(\bar{f}(z_0) - f_0^*(z_0)\big) \end{pmatrix} \rightsquigarrow N\left(\mathbf{0}, \sigma^2 \begin{pmatrix} \boldsymbol{\Omega}^{-1} & \mathbf{0} \\ \mathbf{0} & \Sigma_{22} \end{pmatrix}\right), \tag{3.14}$$

where $\Sigma_{22} = \sigma_{z_0}^2 + 2\boldsymbol{\gamma}_{z_0}^T \boldsymbol{\Omega}^{-1} \boldsymbol{\alpha}_{z_0} + \boldsymbol{\gamma}_{z_0}^T \boldsymbol{\Omega}^{-1} \boldsymbol{\gamma}_{z_0}$;

(ii) if $s$ is fixed, then

$$\begin{pmatrix} \sqrt{n}(\widehat{\boldsymbol{\beta}}^{(j)} - \boldsymbol{\beta}_0^{(j)*}) \\ \sqrt{N/d(\lambda)}\big(\bar{f}(z_0) - f_0^*(z_0)\big) \end{pmatrix} \rightsquigarrow N\left(\mathbf{0}, \sigma^2 \begin{pmatrix} \boldsymbol{\Omega}^{-1} & s^{-1/2}\boldsymbol{\Sigma}_{21} \\ s^{-1/2}\boldsymbol{\Sigma}_{12} & \Sigma_{22} \end{pmatrix}\right), \tag{3.15}$$

where $\boldsymbol{\Sigma}_{12} = \boldsymbol{\Sigma}_{12}^T = \boldsymbol{\Omega}^{-1}(\boldsymbol{\alpha}_{z_0} + \boldsymbol{\gamma}_{z_0})$.

Part (i) of Theorem 3.4 says that $\sqrt{n}\widehat{\boldsymbol{\beta}}^{(j)}$ and $\sqrt{N/d(\lambda)}\bar{f}(z_0)$ are asymptotically independent as $s \to \infty$. This is not surprising since only samples in one sub-population (with size $n$) contribute to the estimation of the heterogeneity component while the entire sample (with size $N$) to commonality. As $n/N = s^{-1} \to 0$, the former data becomes asymptotically independent of (or asymptotically ignorable to) the latter data. So are these two estimators. The estimation bias $P_\lambda m_0^{(j)}$ can be removed by placing a smoothness condition on $B_k$, i.e., Assumption 3.3. Interestingly, given this additional condition, even when $s$ is fixed, these two estimators can still achieve the asymptotic independence if $d(\lambda) \to \infty$. Please see more details in next section.

The norming $d(\lambda)^{1/2}$ is needed in these conditions: $\|\widetilde{K}_{z_0}\|_{L_2(\mathbb{P}_Z)}/d(\lambda)^{1/2} \to \sigma_{z_0}$, $(W_\lambda \boldsymbol{A})(z_0)/d(\lambda)^{1/2} \to \boldsymbol{\alpha}_{z_0} \in \mathbb{R}^p$, and $\boldsymbol{A}(z_0)/d(\lambda)^{1/2} \to -\boldsymbol{\gamma}_{z_0} \in \mathbb{R}^p$ due to the following variance calculation:

$$\begin{aligned}
&\mathrm{Var}\left(\sqrt{N/d(\lambda)}(\bar{f} - f_0^*)\right)\\
&\approx \Big\{\big[\|\widetilde{K}_{z_0}\|_{L_2(\mathbb{P}_Z)}/d(\lambda)^{1/2}\big]^2 + 2\big[\boldsymbol{A}(z_0)/d(\lambda)^{1/2}\big]^T \boldsymbol{\Omega}^{-1}\big[W_\lambda \boldsymbol{A}(z_0)/d(\lambda)^{1/2}\big]\\
&\qquad + \big[\boldsymbol{A}(z_0)/d(\lambda)^{1/2}\big]^T \boldsymbol{\Omega}^{-1}\big[\boldsymbol{A}(z_0)/d(\lambda)^{1/2}\big]\Big\},
\end{aligned}$$

where the first term is dominating and $\|\widetilde{K}_{z_0}\|_{L_2(\mathbb{P}_Z)} = O(d(\lambda)^{1/2})$. So, by the norming $d(\lambda)^{1/2}$, we obtain the order of $\sqrt{N/d(\lambda)}(\bar{f} - f_0^*)$ as $O_{\mathbb{P}}(1)$.

Our last result in this section is the joint asymptotic distribution of $\{\widehat{\boldsymbol{\beta}}^{(j)}\}_{j=1}^s$ (expressed in a linear contrast form). Denote

$$\underset{\sim}{\widehat{\boldsymbol{\beta}}} = \big(\widehat{\boldsymbol{\beta}}^{(1)T}, \ldots, \widehat{\boldsymbol{\beta}}^{(s)T}\big)^T \in \mathbb{R}^{ps} \ \text{ and } \ \underset{\sim}{\boldsymbol{\beta}}_0 = \big(\boldsymbol{\beta}_0^{(1)T}, \ldots, \boldsymbol{\beta}_0^{(s)T}\big)^T \in \mathbb{R}^{ps}.$$

**Theorem 3.5.** Suppose Assumptions 3.1 - 3.3 hold. If $\lambda = o(N^{-1/2})$, and $s$ satisfies (3.12) and

$$s^2 d(\lambda)/N \log^4 N + s\lambda = o\Big(d(\lambda)^{-2}\big(J(\mathcal{F}, 1) + \log N\big)^{-2} \log^{-2} N\Big), \tag{3.16}$$

then for any $\underset{\sim}{\boldsymbol{u}} = (\boldsymbol{u}_1^T, \ldots, \boldsymbol{u}_s^T) \in \mathbb{R}^{ps}$ with $\|\underset{\sim}{\boldsymbol{u}}\|_2 = 1$, it holds

$$\sqrt{n} V_s^{-1} \underset{\sim}{\boldsymbol{u}}^T (\underset{\sim}{\widehat{\boldsymbol{\beta}}} - \underset{\sim}{\boldsymbol{\beta}}_0) \rightsquigarrow N(0, \sigma^2),$$

where $V_s^2 = \sum_{j=1}^s \boldsymbol{u}_j^T \boldsymbol{\Omega}^{-1} \boldsymbol{u}_j$, as $N \to \infty$.

Note that the upper bound condition on $s$ is slightly different from that in Theorem 3.4.

3.4. *Bias Propagation.* In this section, we first analyze the source of estimation bias observed in the joint asymptotics Theorem 3.4. In fact, these analysis leads to a bias propagation phenomenon, which intuitively explains how Assumption 3.3 removes the estimation bias. More importantly, we show that $\bar{f}$ shares exactly the same asymptotic distribution as $\widehat{f}_{or}$, i.e., oracle rule, when $s$ does not grow too fast and $\lambda$ is chosen in the order of $N$.

Our study on propagation mechanism is motivated by the following simple observation. Denote $\mathbb{X} \in \mathbb{R}^{n \times p}$ and $\mathbb{Z} \in \mathbb{R}^n$ as the designs based on the samples from the $j$th sub-population and let $\boldsymbol{\varepsilon}^{(j)} = [\varepsilon_i]_{i \in S_j} \in \mathbb{R}^n$. The first order optimality condition (w.r.t. $\boldsymbol{\beta}$) gives

$$\widehat{\boldsymbol{\beta}}^{(j)} - \boldsymbol{\beta}_0^{(j)} = (\mathbb{X}^T\mathbb{X})^{-1}\mathbb{X}^T \boldsymbol{\varepsilon}^{(j)} - (\mathbb{X}^T\mathbb{X})^{-1}\mathbb{X}^T(\widehat{f}^{(j)}(\mathbb{Z}) - f_0(\mathbb{Z})), \tag{3.17}$$

where $f_0(\mathbb{Z})$ is a $n$-dimensional vector with entries $f_0(Z_i)$ for $i \in S_j$ and $\widehat{f}^{(j)}(\mathbb{Z})$ is defined similarly. Hence, the estimation bias of $\widehat{\boldsymbol{\beta}}^{(j)}$ inherits from that of $\widehat{f}^{(j)}$. A more complete picture on the propagation mechanism can be seen by decomposing the total bias $P_\lambda m_0^{(j)}$ into two parts:

$$\text{parametric bias: } L_\lambda f_0 = -(\boldsymbol{\Omega} + \boldsymbol{\Sigma}_\lambda)^{-1} \langle \boldsymbol{B}, W_\lambda f_0 \rangle_{L_2(\mathbb{P}_Z)}, \tag{3.18}$$

$$\text{nonparametric bias: } N_\lambda f_0 = W_\lambda f_0 - \boldsymbol{A}^T L_\lambda f_0 \tag{3.19}$$

according to Proposition 2.3. The first term in (3.19) explains the bias introduced by penalization; see (2.5). This bias propagates to the parametric component through $\boldsymbol{B}$, as illustrated in (3.18). The parametric bias $L_\lambda f_0$ propagates back to the nonparametric component through the second term of (3.19). Therefore, by strengthening $B_k \in L_2(\mathbb{P}_Z)$ to $B_k \in \mathcal{H}$, i.e., Assumption 3.3, it can be shown that the order of $L_\lambda f_0$ in (3.18) reduces to that of $\lambda$. And then we can remove $L_\lambda f_0$ asymptotically by choosing a sufficiently small $\lambda$. In this case, the nonparametric bias becomes $W_\lambda f_0$.

We summarize the above discussions in the following theorem:

**Theorem 3.6** (Joint Asymptotics II)**.** Suppose Assumption 3.3 and the conditions in Theorem 3.4 hold. If we choose $\lambda = o\big(\sqrt{d(\lambda)/N} \wedge n^{-1/2}\big)$, then (i) if $s \to \infty$ then

$$\begin{pmatrix} \sqrt{n}(\widehat{\boldsymbol{\beta}}^{(j)} - \boldsymbol{\beta}_0^{(j)}) \\ \sqrt{N/d(\lambda)}\big(\bar{f}(z_0) - f_0(z_0) - W_\lambda f_0(z_0)\big) \end{pmatrix} \rightsquigarrow N\left(\boldsymbol{0}, \sigma^2 \begin{pmatrix} \boldsymbol{\Omega}^{-1} & \boldsymbol{0} \\ \boldsymbol{0} & \Sigma_{22}^* \end{pmatrix}\right), \tag{3.20}$$

where $\Sigma_{22}^* = \sigma_{z_0}^2 + \boldsymbol{\gamma}_{z_0}^T \boldsymbol{\Omega}^{-1} \boldsymbol{\gamma}_{z_0}$;
(ii) if $s$ is fixed, then

$$\begin{pmatrix} \sqrt{n}(\widehat{\boldsymbol{\beta}}^{(j)} - \boldsymbol{\beta}_0^{(j)}) \\ \sqrt{N/d(\lambda)}\big(\bar{f}(z_0) - f_0(z_0) - W_\lambda f_0(z_0)\big) \end{pmatrix} \rightsquigarrow N\left(\boldsymbol{0}, \sigma^2 \begin{pmatrix} \boldsymbol{\Omega}^{-1} & s^{-1/2}\boldsymbol{\Sigma}_{21}^* \\ s^{-1/2}\boldsymbol{\Sigma}_{12}^* & \Sigma_{22}^* \end{pmatrix}\right), \tag{3.21}$$

where $\boldsymbol{\Sigma}_{12}^* = \boldsymbol{\Sigma}_{12}^{*T} = \boldsymbol{\Omega}^{-1}\boldsymbol{\gamma}_{z_0}$ and $\Sigma_{22}^*$ is the same as in (i).
Moreover, if $d(\lambda) \to \infty$, then $\boldsymbol{\Sigma}_{12}^* = \boldsymbol{\Sigma}_{21}^* = \boldsymbol{0}$ and $\Sigma_{22}^* = \sigma_{z_0}^2$ in (i) and (ii).

The nonparametric estimation bias $W_\lambda f_0(z_0)$ can be further removed by performing undersmoothing, a standard procedure in nonparametric inference, e.g., Shang and Cheng (2013). We will illustrate this point in Section 4.

By examining the proof for case (ii) of Theorem 3.6 (and taking $s = 1$), we know that the oracle estimate $\widehat{f}_{or}$ defined in (3.1) attains the same asymptotic distribution as that of $\bar{f}$ in (3.20) when $s$ grows at a proper rate. Therefore, we claim that our combined estimate $\bar{f}$ satisfies the desirable oracle property.

In Section 4, we apply Theorem 3.6 to several examples, and find that even though the minimization (3.2) is based only on one fraction of the entire sample, it is nonetheless essential to regularize each sub-estimation as if it had the entire sample. In other words, $\lambda$ should be chosen in the order of $N$. Similar phenomenon also arises in analyzing minimax optimality of each sub-estimation; see Section 3.2.

When $s = 1$, our model reduces to the standard partially linear model. The joint distribution of parametric and nonparametric estimators is shown in Part (ii) of Theorem 3.6, where their asymptotic covariance is derived as $\boldsymbol{\Omega}^{-1}\boldsymbol{\gamma}_{z_0}$ with $\boldsymbol{\gamma}_{z_0} = \lim_{N\to\infty} -\boldsymbol{A}(z_0)/d(\lambda)^{1/2}$. In Section 7.5, we show that $\boldsymbol{A}(z_0)$ is bounded for any $z_0 \in \mathcal{Z}$ (uniformly over $\lambda$). Hence, the asymptotic correlation disappears, i.e., $\boldsymbol{\gamma}_{z_0} = \mathbf{0}$, as $d(\lambda) \to \infty$. This corresponds to the exponentially decaying kernel and polynomially decaying kernel. In fact, this finding generalizes the *joint asymptotics phenomenon* recently discovered for partial smoothing spline models; see Cheng and Shang (2015). However, when $d(\lambda)$ is finite, e.g., finite rank kernel, the asymptotic correlation remains. This is not surprising since the semiparametric estimation in this case essentially reduces to a parametric one.

**Remark 3.1.** Theorem 3.6 implies that $\sqrt{n}(\widehat{\boldsymbol{\beta}}^{(j)} - \boldsymbol{\beta}_0^{(j)}) \rightsquigarrow N(\mathbf{0}, \sigma^2\boldsymbol{\Omega}^{-1})$ when $\lambda = o(n^{-1/2})$. When the error $\varepsilon$ follows a Gaussian distribution, it is well known that $\widehat{\boldsymbol{\beta}}^{(j)}$ achieves the semiparametric efficiency bound (Kosorok, 2007). Hence, the semiparametric efficient estimate can be obtained by applying the kernel ridge method. However, we can further improve its estimation efficiency to a parametric level by taking advantage of $\bar{f}$ (built on the whole samples). This represents an important feature of massive data: strength-borrowing.

**Remark 3.2.** We can also construct a simultaneous confidence band for $f_0$ based on the stochastic expansion of $\bar{f}$ and strong approximation techniques (Bickel and Rosenblatt, 1973). Specifically, we start from (3.9) that implies

$$\Big\|\bar{f} - f_0^* - \frac{1}{N}\sum_{i=1}^{N}\varepsilon_i N_{U_i}\Big\|_{\sup} = \Big\|\frac{1}{s}\sum_{j=1}^{s} Rem_f^{(j)}\Big\|_{\sup}. \tag{3.22}$$

Similar to the pointwise case, we can show that the remainder term on the R.H.S. of (3.22) is $o_P(1)$ once $s$ does not grow too fast. Hence, the distribution of $\sup_z |\bar{f}(z) - f_0^*(z)|$ can be approximated by that of $\sup_z |N^{-1}\sum_{i=1}^N \varepsilon_i N_{U_i}(z)|$. We next apply strong approximation techniques to prove that $N^{-1}\sum_{i=1}^N \varepsilon_i N_{U_i}$ can be further approximated by a proper Gaussian process. This would yield a simultaneous confidence band. More rigorous arguments can be adapted from the proof of Theorem 5.1 in Shang and Cheng (2013).

3.5. *Efficiency Boosting: from semiparametric level to parametric level.* The previous sections show that the combined estimate $\bar{f}$ achieves the "oracle property" in both asymptotic and non-asymptotic senses when $s$ does not grow too fast and $\lambda$ is chosen according to the entire sample size. In this section, we further employ $\bar{f}$ to boost the estimation efficiency of $\widehat{\boldsymbol{\beta}}^{(j)}$ from semiparametric level to parametric level. This leads to our final estimate for heterogeneity, i.e., $\check{\boldsymbol{\beta}}^{(j)}$ defined in (3.23). More importantly, $\check{\boldsymbol{\beta}}^{(j)}$ possesses the limit distribution as if the commonality in each sub-population were known, and hence satisfies the "oracle rule". This interesting efficiency boosting phenomenon will be empirically verified in Section 6. A similar two-stage estimation method was proposed in Li and Liang (2008), but for the purpose of variable selection in $\boldsymbol{\beta}$ based on a single data set.

Specifically, we define the following improved estimator for $\boldsymbol{\beta}_0$:

$$\check{\boldsymbol{\beta}}^{(j)} = \underset{\boldsymbol{\beta}\in\mathbb{R}^p}{\operatorname{argmin}} \frac{1}{n} \sum_{i\in S_j} \big(Y_i - \boldsymbol{X}_i^T\boldsymbol{\beta} - \bar{f}(Z_i)\big)^2. \tag{3.23}$$

Theorem 3.7 below shows that $\check{\boldsymbol{\beta}}^{(j)}$ achieves the parametric efficiency bound as if the nonparametric component $f$ were known. This is not surprising given that the nonparametric estimate $\bar{f}$ now possesses a faster convergence rate after aggregation. What is truly interesting is that we need to set a lower bound for $s$, i.e., (3.24), which slows down the convergence rate of $\check{\boldsymbol{\beta}}^{(j)}$, i.e., $\sqrt{n}$, such that $\bar{f}$ can be treated as if it were known. Note that the homogeneous data setting is trivially excluded in this case.

**Theorem 3.7.** Suppose Assumptions 3.1 and 3.2 hold. If $s$ satisfies conditions (3.12), (3.13) and

$$s^{-1} = o\big(d(\lambda)^{-2}\log^{-4} N\big), \tag{3.24}$$

and we choose $\lambda = o(d(\lambda)/N)$, then we have

$$\sqrt{n}(\check{\boldsymbol{\beta}}^{(j)} - \boldsymbol{\beta}_0^{(j)}) \rightsquigarrow N(0, \sigma^2\boldsymbol{\Sigma}^{-1}),$$

where $\boldsymbol{\Sigma} = \mathbb{E}[\boldsymbol{X}\boldsymbol{X}^T]$.

Note that $\boldsymbol{X}$ and $Z$ are not assumed to be independent. Hence, the parametric efficiency bound $\boldsymbol{\Sigma}^{-1}$ is not larger than the semiparametric efficiency bound $\boldsymbol{\Omega}^{-1}$. The intuition for this lower bound of $s$ is that the total sample size $N$ should grow much faster than the subsample size $n$, so that the nonparametric estimator $\bar{f}$ converges faster than the parametric estimator $\check{\boldsymbol{\beta}}^{(j)}$. In this case, the error of estimating $f_0$ is negligible so that $\check{\boldsymbol{\beta}}^{(j)}$ behaves asymptotically as if $f$ were known, resulting in *parametric* efficiency.

3.6. *Testing for Heterogeneity.* The heterogeneity across different sub-populations is a crucial feature of massive data. However, there is still some chance that some sub-populations may share the same underlying distribution. In this section, we consider testing for the heterogeneity among sub-populations. We start from a simple pairwise testing, and then extend it to a more challenging simultaneous testing that can be applied to a large number of sub-populations.

Consider a general class of pairwise heterogeneity testing:

$$H_0 : Q(\boldsymbol{\beta}_0^{(j)} - \boldsymbol{\beta}_0^{(k)}) = \mathbf{0} \ \text{ for } j \neq k, \tag{3.25}$$

where $Q = (Q_1^T, \dots, Q_q^T)^T$ is a $q \times p$ matrix with $q \le p$. The general formulation (3.25) can test either the whole vector or one fraction of $\boldsymbol{\beta}_0^{(j)}$ is equal to that of $\boldsymbol{\beta}_0^{(k)}$. A test statistic can be constructed based on either $\widehat{\boldsymbol{\beta}}$ or its improved version $\check{\boldsymbol{\beta}}$. Let $C_\alpha \subset \mathbb{R}^q$ be a confidence region satisfying $\mathbb{P}(\boldsymbol{b} \in C_\alpha) = 1 - \alpha$ for any $\boldsymbol{b} \sim N(0, I_q)$. Specifically, we have the following $\alpha$-level Wald tests:

$$\begin{aligned} \Psi_1 &= I\{Q\big(\widehat{\boldsymbol{\beta}}^{(j)} - \widehat{\boldsymbol{\beta}}^{(k)}\big) \notin \sqrt{2/n}\sigma(Q\boldsymbol{\Omega}^{-1}Q^T)^{1/2}C_\alpha\}, \\ \Psi_2 &= I\{Q\big(\check{\boldsymbol{\beta}}^{(j)} - \check{\boldsymbol{\beta}}^{(k)}\big) \notin \sqrt{2/n}\sigma(Q\boldsymbol{\Sigma}^{-1}Q^T)^{1/2}C_\alpha\}. \end{aligned}$$

The consistency of the above tests are guaranteed by Theorem 3.8 below. In addition, we note that the power of the latter test is larger than the former; see the analysis below Theorem 3.8. The price we need to pay for this larger power is to require a lower bound on $s$.

**Theorem 3.8.** Suppose that the conditions in Theorem 3.6 are satisfied. Under the null hypothesis specified in (3.25), we have

$$\sqrt{n}Q\big(\widehat{\boldsymbol{\beta}}^{(j)} - \widehat{\boldsymbol{\beta}}^{(k)}\big) \rightsquigarrow N(\mathbf{0}, 2\sigma^2 Q\boldsymbol{\Omega}^{-1}Q^T).$$

Moreover, under the conditions in Theorem 3.7, we have

$$\sqrt{n}Q\big(\check{\boldsymbol{\beta}}^{(j)} - \check{\boldsymbol{\beta}}^{(k)}\big) \rightsquigarrow N(\mathbf{0}, 2\sigma^2 Q\boldsymbol{\Sigma}^{-1}Q^T),$$

where $\boldsymbol{\Sigma} = \mathbb{E}[\boldsymbol{X}\boldsymbol{X}^T]$.

The larger power of $\Psi_2$ is due to the smaller asymptotic variance of $\check{\boldsymbol{\beta}}^{(j)}$, and can be deduced from the following power function. For simplicity, we consider $H_0 : \beta_{01}^{(j)} - \beta_{01}^{(k)} = 0$, i.e., $Q = (1, 0, 0 \dots, 0)$. In this case, we have $\Psi_1 = I\big\{|\widehat{\beta}_1^{(j)} - \widehat{\beta}_1^{(k)}| > \sqrt{2}\sigma[\boldsymbol{\Omega}^{-1}]_{11}^{1/2} z_{\alpha/2}/\sqrt{n}\big\}$, and $\Psi_2 = I\big\{|\check{\beta}_1^{(j)} - \check{\beta}_1^{(k)}| > \sqrt{2}\sigma[\boldsymbol{\Sigma}^{-1}]_{11}^{1/2} z_{\alpha/2}/\sqrt{n}\big\}$. The (asymptotic) power function under the alternative that $\beta_{01}^{(j)} - \beta_{01}^{(k)} = \beta^*$ for some non-zero $\beta^*$ is

$$\text{Power}(\beta^*) = 1 - \mathbb{P}\Big(W \in \Big[-\frac{\beta^*\sqrt{n}}{\sigma^*} \pm z_{\alpha/2}\Big]\Big),$$

where $W \sim N(0,1)$ and $\sigma^*$ is $\sqrt{2}\sigma[\boldsymbol{\Omega}^{-1}]_{11}^{1/2}$ for $\Psi_1$ and $\sqrt{2}\sigma[\boldsymbol{\Sigma}^{-1}]_{11}^{1/2}$ for $\Psi_2$. Hence, a smaller $\sigma^*$ gives rise to a larger power, and $\Psi_2$ is more powerful than $\Psi_1$. Please see Section 6 for empirical support for this power comparison.

We next consider the problem of heterogeneous testing for a large number of sub-populations:

$$H_0 : \boldsymbol{\beta}^{(j)} = \widetilde{\boldsymbol{\beta}}^{(j)} \text{ for all } j \in \mathcal{G}, \tag{3.26}$$

where $\mathcal{G} \subset \{1, 2, \ldots, s\}$, versus the alternative:

$$H_1 : \boldsymbol{\beta}^{(j)} \neq \widetilde{\boldsymbol{\beta}}^{(j)} \text{ for some } j \in \mathcal{G}. \tag{3.27}$$

The above $\widetilde{\boldsymbol{\beta}}^{(j)}$'s are pre-specified for each $j \in \mathcal{G}$. If all $\widetilde{\boldsymbol{\beta}}^{(j)}$'s are the same, then it becomes a type of heterogeneity test for the group of sub-populations indexed by $\mathcal{G}$. Here we allow $|\mathcal{G}|$ to be as large as $s$, and thus it can increase with $n$. Let $\widehat{\boldsymbol{\Sigma}}^{(j)}$ be the sample covariance matrix of $\boldsymbol{X}$ for the $j$-th sub-population, i.e., $n^{-1} \sum_{i \in S_j} \boldsymbol{X}_i \boldsymbol{X}_i^T$. Define the test statistic

$$T_{\mathcal{G}} := \max_{j \in \mathcal{G}, 1 \le k \le p} \sqrt{n}(\check{\boldsymbol{\beta}}_k^{(j)} - \widetilde{\boldsymbol{\beta}}_k^{(j)}).$$

We approximate the distribution of the above test statistic using multiplier bootstrap. Define the following quantity:

$$W_{\mathcal{G}} := \max_{j \in \mathcal{G}, 1 \le k \le p} \frac{1}{\sqrt{n}} \sum_{i \in S_j} \big(\widehat{\boldsymbol{\Sigma}}^{(j)}\big)_k^{-1} \boldsymbol{X}_i e_i,$$

where $e_i$'s are i.i.d. $N(0, \sigma^2)$ independent of the data and $\big(\widehat{\boldsymbol{\Sigma}}^{(j)}\big)_k^{-1}$ is the $k$-th row of $\big(\widehat{\boldsymbol{\Sigma}}^{(j)}\big)^{-1}$. Let $c_{\mathcal{G}}(\alpha) = \inf\{t \in \mathbb{R} : \mathbb{P}(W_{\mathcal{G}} \le t \,|\, \mathbb{X}) \ge 1 - \alpha\}$. We employ the recent Gaussian approximation and multiplier bootstrap theory (Chernozhukov et al., 2013) to obtain the following theorem.

**Theorem 3.9.** Suppose Assumptions 3.1 and 3.2 hold. In addition, suppose (3.12) and (3.13) in Theorem 3.4 hold. For any $\mathcal{G} \subset \{1, 2, \ldots, s\}$ with $|\mathcal{G}| = d$, if (i) $s \gtrsim d(\lambda)^2 \log(pd) \log^4 N$, (ii) $(\log(pdn))^7/n \le C_1 n^{-c_1}$ for some constants $c_1, C_1 > 0$, and (iii) $p^2 \log(pd)/\sqrt{n} = o(1)$, then under $H_0$ and choosing $\lambda = o(d(\lambda)/N)$, we have

$$\sup_{\alpha \in (0,1)} \left| \mathbb{P}\Big(T_{\mathcal{G}} > c_{\mathcal{G}}(\alpha)\Big) - \alpha \right| = o(1).$$

**Remark 3.3.** We can perform heterogeneity testing even without specifying $\widetilde{\boldsymbol{\beta}}^{(j)}$'s. This can be done by simply reformulating the null hypothesis as follows (for simplicity we set $\mathcal{G} = [s]$): $H_0 : \alpha^{(j)} = 0$ for $j \in [s-1]$, where $\alpha^{(j)} = \boldsymbol{\beta}^{(j)} - \boldsymbol{\beta}^{(j+1)}$ for $j = 1, \ldots, s-1$. The test statistic is $T'_{\mathcal{G}} = \max_{1 \le j \le s-1} \max_{1 \le k \le p} \alpha_k^{(j)}$. The bootstrap quantity is defined as

$$W'_{\mathcal{G}} := \max_{1 \le j \le s-1, 1 \le k \le p} \frac{1}{\sqrt{n}} \sum_{i \in S_j} \big(\widehat{\boldsymbol{\Sigma}}^{(j)}\big)_k^{-1} \boldsymbol{X}_i e_i - \frac{1}{\sqrt{n}} \sum_{i \in S_{j+1}} \big(\widehat{\boldsymbol{\Sigma}}^{(j+1)}\big)_k^{-1} \boldsymbol{X}_i e_i.$$

The proof is similar to that of Theorem 3.9 and is omitted.

## 4. Examples.

In this section, we consider three specific classes of RKHS with different smoothness, characterized by the decaying rate of the eigenvalues: finite rank, exponential decay and polynomial decay. In particular, we give explicit upper bounds for $s$ under which the combined estimate enjoys the oracle property, and also explicit lower bounds for obtaining efficiency boosting studied in Section 3.5. Interestingly, we find that the upper bound for $s$ increases for RKHS with faster decaying eigenvalues. Hence, our aggregation procedure favors smoother regression functions in the sense that more sub-populations are allowed to be included in the observations. The choice of $\lambda$ is also explicitly characterized in terms of the entire sample size and the decaying rate of eigenvalues. In all three examples, the undersmoothing is implicitly assumed for removing the nonparametric estimation bias. Our bounds on $s$ and $\lambda$ here are not the most general ones, but are those that can easily deliver theoretical insights.

### 4.1. *Example I: Finite Rank Kernel.*

The RKHS with finite rank kernels includes linear functions, polynomial functions, and, more generally, functional classes with finite dictionaries. In this case, the effective dimension is simply proportional to the rank $r$. Hence, $d(\lambda) \asymp r$. Combining this fact with Theorem 3.6, we get the following corollary for finite rank kernels:

**Corollary 4.1.** Suppose Assumption 3.1 – 3.3 hold and $s \to \infty$. For any $z_0 \in \mathcal{Z}$, if $\lambda = o\big(N^{-1/2}\big)$, $\log(\lambda^{-1}) = o\big(N^2 \log^{-12} N\big)$ and $s = o\Big(\frac{N}{\sqrt{\log \lambda^{-1}} \log^6 N}\Big)$, then

$$\begin{pmatrix} \sqrt{n}(\widehat{\boldsymbol{\beta}}^{(j)} - \boldsymbol{\beta}_0^{(j)}) \\ \sqrt{N}(\bar{f}(z_0) - f_0(z_0)) \end{pmatrix} \rightsquigarrow N\left(\mathbf{0}, \sigma^2 \begin{pmatrix} \boldsymbol{\Omega}^{-1} & \mathbf{0} \\ \mathbf{0} & \Sigma_{22}^* \end{pmatrix}\right),$$

where $\Sigma_{22}^* = \sum_{\ell=1}^r \phi_\ell(z_0)^2 + \boldsymbol{\gamma}_{z_0}^T \boldsymbol{\Omega}^{-1} \boldsymbol{\gamma}_{z_0}$ and $\boldsymbol{\gamma}_{z_0} = \sum_{\ell=1}^r \langle \boldsymbol{B}, \phi_\ell \rangle_{L_2(\mathbb{P}_Z)} \phi_\ell(z_0)$.

From the above Corollary, we can easily tell that the upper bound for $s$ can be as large as $o\big(N \log^{-7} N\big)$ by choosing a sufficiently large $\lambda$. Hence, $s$ can be chosen nearly as large as $N$. As for the lower bound of $s$ for boosting the efficiency, we have $s \gtrsim r^2 \log^4 N$ by plugging $d(\lambda) \asymp r$ into (3.24). This lower bound is clearly smaller than the upper bound. Hence, the efficiency boosting is feasible.

Corollary 4.2 below specifies conditions on $s$ and $\lambda$ under which $\bar{f}$ achieves the nonparametric minimaxity.

**Corollary 4.2.** Suppose that Assumptions 3.1 - 3.3 hold. When $\lambda = r/N$

and $s = o\big(N\log^{-5} N\big)$, we have

$$\mathbb{E}\big[\|\bar{f} - f_0\|^2_{L_2(\mathbb{P}_Z)}\big] \le Cr/N,$$

for some constant $C$.

4.2. *Example II: Exponential Decay Kernel.* We next consider the RKHS for which the kernel has exponentially decaying eigenvalues, i.e., $\mu_\ell = \exp(-\alpha\ell^p)$ for some $\alpha > 0$. In this case, we have $d(\lambda) \asymp (\log\lambda^{-1})^{1/p}$ by explicit calculations.

**Corollary 4.3.** Suppose Assumptions 3.1 – 3.3 hold, and for any $z_0 \in \mathcal{Z}$, $f_0 \in \mathcal{H}$ satisfies $\sum_{\ell=1}^{\infty} |\phi_\ell(z_0)\langle f_0, \phi_\ell\rangle_{\mathcal{H}}| < \infty$. If $\lambda = o\big(N^{-1/2}\log^{1/(2p)} N \wedge n^{-1/2}\big)$, $\log(\lambda^{-1}) = o\big(N^{p/(p+4)}\log^{-6p/(p+4)} N\big)$ and $s = o\Big(\frac{N}{\log^6 N \log^{(p+4)/p}(\lambda^{-1})}\Big)$, then

$$\begin{pmatrix} \sqrt{n}(\widehat{\boldsymbol\beta}^{(j)} - \boldsymbol\beta_0^{(j)}) \\ \sqrt{N/d(\lambda)}(\bar{f}(z_0) - f_0(z_0)) \end{pmatrix} \rightsquigarrow N\left(\mathbf{0}, \sigma^2\begin{pmatrix} \boldsymbol\Omega^{-1} & \mathbf{0} \\ \mathbf{0} & \sigma^2_{z_0} \end{pmatrix}\right),$$

where $\sigma^2_{z_0} = \lim_{\lambda\to 0} d(\lambda)^{-1}\sum_{\ell=1}^{\infty}\frac{\phi_\ell^2(z_0)}{(1+\lambda/\mu_\ell)^2}$.

Corollary 4.3 implies the shrinking rate of the confidence interval for $f_0(z_0)$ as $\sqrt{d(\lambda)/N}$. This motivates us to choose $\lambda$ (equivalently $d(\lambda)$) as large as possible (as small as possible). Plugging such a $\lambda$ into the upper bound of $s$ yields $s = o(N\log^{-(7p+4)/p} N)$. For example, when $p = 1 (p = 2)$, the upper bound is $s = o(N\log^{-11} N) (s = o(N\log^{-9} N))$. Note that this upper bound for $s$ only differs from that for the finite rank kernel up to some logrithmic term. This is mainly because RKHS with exponentially decaying eigenvalues has an effective dimension $(\log N)^{1/p}$ (for the above $\lambda$). Again, by (3.24) we get the lower bound of $s \gtrsim (\log\lambda^{-1})^{2/p}\log^2 N$. When $\lambda \asymp N^{-1/2}\log^{1/(2p)} N \wedge n^{-1/2}$, it is approximately $s \gtrsim \log^{(4p+2)/p} N$.

As a concrete example, we consider the Gaussian kernel $K(z_1, z_2) = \exp\big(-|z_1 - z_2|^2/2\big)$. The eigenfunctions are given in (2.1), and the eigenvalues are exponentially decaying, as $\mu_\ell = \eta^{2\ell+1}$, where $\eta = (\sqrt{5} - 1)/2$. According to Krasikov (2004), we can get that

$$c_\phi = \sup_{\ell\in\mathbb{N}} \|\phi_\ell\|_{\sup} \le \frac{2e^{15/8}(\sqrt{5}/4)^{1/4}}{3\sqrt{2\pi}2^{1/6}} \le 1.336.$$

Thus, Assumption 3.2 is satisfied. We next give an upper bound of $\sigma^2_{z_0}$ in Corollary 4.3 as follows:

$$\sigma^2_{z_0} \le \lim_{N\to\infty} \sigma^2 c_\phi^2 h \sum_{\ell=0}^{\infty}(1+\lambda\eta\exp(-2(\log\eta)\ell))^{-2} = c_\phi^2{\cdot}2\sigma^2\log(1/\eta) \le 1.7178\sigma^2,$$

where equality follows from Lemma C.1 in Appendix C for the case $t = 2$. Hence, a (conservative) $100(1-\alpha)\%$ confidence interval for $f_0(z_0)$ is given by $\bar{f}(z_0) \pm 1.3106\sigma z_{\alpha/2}\sqrt{d(\lambda)/N}$.

**Corollary 4.4.** Suppose that Assumptions 3.1 – 3.3 hold. By choosing $\lambda = (\log N)^{1/p}/N$ and $s = o\big(N \log^{-(5p+3)/p} N\big)$, we have

$$\mathbb{E}\big[\|\bar{f} - f_0\|^2_{L_2(\mathbb{P}_Z)}\big] \le C(\log N)^{1/p}/N.$$

We know that the above rate is minimax optimal according to Zhang et al. (2013). Note that the upper bound for $s$ required here is similar as that for obtaining the joint limiting distribution in Corollary 4.3.

4.3. *Example III: Polynomial Decay Kernel.* We now consider the RKHS for which the kernel has polynomially decaying eigenvalues, i.e., $\mu_\ell = c\ell^{-2\nu}$ for some $\nu > 1/2$. Hence, we can explicitly calculate that $d(\lambda) = \lambda^{-1/(2\nu)}$. The resulting penalized estimate is called as "partial smoothing spline" in the literature; see Gu (2013); Wang (2011).

**Corollary 4.5.** Suppose Assumptions 3.1 – 3.3 hold, and $\sum_{\ell=1}^\infty |\phi_\ell(z_0)\langle f_0, \phi_\ell\rangle_{\mathcal{H}}| < \infty$ for any $z_0 \in \mathcal{Z}$ and $f_0 \in \mathcal{H}$. For any $\nu > 1 + \sqrt{3}/2 \approx 1.866$, if $\lambda \asymp N^{-d}$ for some $\frac{2\nu}{4\nu+1} < d < \frac{4\nu^2}{10\nu-1}$, $\lambda = o(n^{-1/2})$ and $s = o\Big(\lambda^{\frac{10\nu-1}{4\nu^2}} N \log^{-6} N\Big)$, then

$$\begin{pmatrix} \sqrt{n}(\widehat{\boldsymbol{\beta}}^{(j)} - \boldsymbol{\beta}_0^{(j)}) \\ \sqrt{N/d(\lambda)}(\bar{f}(z_0) - f_0(z_0)) \end{pmatrix} \rightsquigarrow N\left(\mathbf{0}, \sigma^2\begin{pmatrix} \boldsymbol{\Omega}^{-1} & \mathbf{0} \\ \mathbf{0} & \sigma^2_{z_0} \end{pmatrix}\right).$$

where $\sigma^2_{z_0} = \lim_{\lambda\to 0} d(\lambda)^{-1}\sum_{\ell=1}^\infty \frac{\phi^2_\ell(z_0)}{(1+\lambda/\mu_\ell)^2}$.

Similarly, we choose $\lambda \asymp N^{-\frac{2\nu}{4\nu+1}} \wedge n^{-1/2}$ to get the fastest shrinking rate of the confidence interval. Plugging the above $\lambda$ into the upper bound for $s$, we get

$$s = o\Big(N^{\frac{8\nu^2-8\nu+1}{2\nu(4\nu+1)}} \log^{-6} N \wedge N(\log N)^{-\frac{48\nu^2}{8\nu^2+10\nu+1}}\Big).$$

When $N$ is large, the above bound reduces to $s = o\Big(N^{\frac{8\nu^2-8\nu+1}{2\nu(4\nu+1)}} \log^{-6} N\Big)$. We notice that the upper bound for $s$ increases as $\nu$ increases, indicating that the aggregation procedure favors smoother functions. As an example, for the case that $\nu = 2$, we have the upper bound for $s = o(N^{17/36}\log^{-6} N) \approx o(N^{0.47}\log^{-6} N)$. Again, we obtain the lower bound $s \gtrsim \lambda^{-1/\nu}\log^4 N$ by plugging $d(\lambda) \asymp \lambda^{-\frac{1}{2\nu}}$ into (3.24). When $\lambda \asymp N^{-\frac{2\nu}{4\nu+2}}$, we get $s \gtrsim N^{\frac{1}{4\nu+1}}\log^2 N$. For $\nu = 2$, this is approximately $s \gtrsim N^{0.22}\log^4 N$.

As a concrete example, we consider the periodic Sobolev space $H_0^\nu[0,1]$

with the following eigenfunctions:

$$\phi_\ell(x) = \begin{cases} 1, & \ell = 0, \\ \sqrt{2}\cos(\ell\pi x), & \ell = 2k \text{ for } k = 1, 2, \ldots, \\ \sqrt{2}\sin((\ell+1)\pi x), & \ell = 2k-1 \text{ for } k = 1, 2, \ldots, \end{cases} \tag{4.1}$$

and eigenvalues

$$\mu_\ell = \begin{cases} \infty, & \ell = 0, \\ (\ell\pi)^{-2\nu}, & \ell = 2k \text{ for } k = 1, 2, \ldots, \\ ((\ell+1)\pi)^{-2\nu}, & \ell = 2k-1 \text{ for } k = 1, 2, \ldots, \end{cases} \tag{4.2}$$

Hence, Assumption 3.2 trivially holds. Under the above eigensystem, the following lemma gives an explicit expression of $\sigma_{z_0}^2$.

**Lemma 4.1.** Under the eigen-system defined by (4.1) and (4.2), we can explicitly calculate:

$$\sigma_{z_0}^2 = \lim_{\lambda\to 0} d(\lambda)^{-1}\sum_{\ell=1}^{\infty}\frac{\phi_\ell^2(z_0)}{(1+\lambda/\mu_\ell)^2} = \int_0^\infty \frac{1}{(1+x^{2\nu})^2}dx = \frac{\pi}{2\nu\sin(\pi/(2\nu))}.$$

Therefore, by Corollary 4.5, we have that when $\lambda \asymp N^{-\frac{2\nu}{4\nu+1}}$ and $s = o\Big(N^{\frac{8\nu^2-8\nu+1}{2\nu(4\nu+1)}}\log^{-6} N\Big)$,

$$\begin{pmatrix} \sqrt{n}(\widehat{\boldsymbol\beta}^{(j)} - \boldsymbol\beta_0^{(j)}) \\ \sqrt{N/d(\lambda)}(\bar f(z_0) - f_0(z_0)) \end{pmatrix} \rightsquigarrow N\left(\mathbf{0}, \sigma^2\begin{pmatrix} \boldsymbol\Omega^{-1} & \mathbf{0} \\ \mathbf{0} & \sigma_{z_0}^2 \end{pmatrix}\right). \tag{4.3}$$

where $\sigma_{z_0}^2$ is given in Lemma 4.1. When $\nu = 2$, $\lambda \asymp N^{-4/9}$ and the upper bound for $s = o(N^{17/36}\log^{-6} N)$.

**Corollary 4.6.** Suppose that Assumption 3.1 - 3.3 hold. If we choose $\lambda = N^{-\frac{2\nu}{2\nu+1}}$, and $s = o\Big(N^{\frac{4\nu^2-4\nu+1}{4\nu^2+2\nu}}\log^{-4} N\Big)$, the combined estimator achieves optimal rate of convergence, i.e.,

$$\mathbb{E}\big[\|\bar f - f_0\|_{L_2(\mathbb{P}_Z)}^2\big] \le CN^{-\frac{2\nu}{2\nu+1}}. \tag{4.4}$$

The above rate is known to be minimax optimal for the class of functions in consideration (Stone, 1985).

## 5. Application to Homogeneous Data: Divide-and-Conquer Approach.

In this section, we apply the divide-and-conquer approach, which is commonly used to deal with massive homogeneous data, to some sub-populations that have huge sample sizes. A general goal of this section is to explore the most computationally efficient way to split the sample in those sub-populations while preserving the best possible statistical inference. Specifically, we want to derive the largest possible number of splits under which the averaged estimators for both components enjoy the same statistical performances as the "oracle" estimator that is computed based on the entire sample. Without loss of generality, we assume the entire sample to be homogeneous by setting all $\boldsymbol{\beta}_0^{(j)}$'s to be equal throughout this section. It is worth mentioning that Li et al. (2013) have done an earlier and interesting work on parametric or nonparametric models.

The divide-and-conquer method *randomly* splits the massive data into $s$ mutually exclusive subsamples. For simplicity, we assume all the subsamples share the same sample size, denoted as $n$. Hence, $N = n \times s$. With a bit abuse of notation, we define the divide-and-conquer estimators as $\widehat{\boldsymbol{\beta}}^{(j)}$ and $\widehat{f}^{(j)}$ when they are based on the $j$-th subsample. Thus, the averaged estimator is defined as

$$\bar{\boldsymbol{\beta}} = (1/s)\sum_{j=1}^{s}\widehat{\boldsymbol{\beta}}^{(j)} \text{ and } \bar{f}(\cdot) = (1/s)\sum_{j=1}^{s}\widehat{f}^{(j)}(\cdot).$$

Comparing to the oracle estimator, the aggregation procedure reduces the computational complexity in terms of the entire sample size $N$ to the subsample size $N/s$. In the case of kernel ridge regression, the complexity is $O(N^3)$, while our aggregation procedure (run in one single machine) reduces it to $O(N^3/s^2)$. Propositions 5.1 below state conditions under which the divide-and-conquer estimators maintain the same statistical properties as oracle estimate, i.e., so-called oracle property.

**Proposition 5.1.** Suppose that the conditions in Theorem 3.6 hold. If we choose $\lambda = o(N^{-1/2})$, then

$$\begin{pmatrix} \sqrt{N}(\bar{\boldsymbol{\beta}} - \boldsymbol{\beta}_0) \\ \sqrt{N/d(\lambda)}(\bar{f}(z_0) - f_0(z_0) - W_\lambda f_0(z_0)) \end{pmatrix} \rightsquigarrow N\left(\mathbf{0}, \begin{pmatrix} \sigma^2\boldsymbol{\Omega}^{-1} & \boldsymbol{\Sigma}_{12}^* \\ \boldsymbol{\Sigma}_{21}^* & \Sigma_{22}^* \end{pmatrix}\right),$$

where $\boldsymbol{\Sigma}_{12}^* = \boldsymbol{\Sigma}_{21}^{*T} = \sigma^2\boldsymbol{\Omega}^{-1}\boldsymbol{\gamma}_{z_0}$ and $\Sigma_{22}^* = \sigma^2\left(\sigma_{z_0}^2 + \boldsymbol{\gamma}_{z_0}^T\boldsymbol{\Omega}^{-1}\boldsymbol{\gamma}_{z_0}\right)$. Moreover, if $d(\lambda) \to \infty$, then $\boldsymbol{\gamma}_{z_0} = \mathbf{0}$. In this case, $\boldsymbol{\Sigma}_{12}^* = \boldsymbol{\Sigma}_{21}^{*T} = \mathbf{0}$ and $\Sigma_{22}^* = \sigma^2\sigma_{z_0}^2$.

The conclusion of Proposition 5.1 holds no matter $s$ is fixed or diverges (once the condition for $s$ in Theorem 3.6 are satisfied). In view of Proposition 5.1, we note that the above joint asymptotic distribution is exactly the same as that for the oracle estimate, i.e., $s = 1$.

**Remark 5.1.** We can also derive the minimax rate of $\mathrm{MSE}(\bar{f})$, which is exactly the same as that in Theorem 3.2, based on similar proof techniques.

**6. Numercial Experiment.** In this section, we empirically examine the impact of the number of sub-populations on the statistical inference built on $(\widehat{\boldsymbol{\beta}}^{(j)}, \bar{f})$. As will be seen, the simulation results strongly support our general theory.

Specifically, we consider the partial smoothing spline models in Section 4.3. In the simulation setup, we let $\varepsilon \sim N(0,1)$, $p = 1$ and $\nu = 2$ (cubic spline). Moreover, $Z \sim \text{Uniform}(-1,1)$ and $X = (W+Z)/2$, where $W \sim \text{Uniform}(-1,1)$, such that $X$ and $Z$ are dependent. It is easy to show that $\mathbf{\Omega} = E\big[(X - E[X \mid Z])^2)\big] = 1/12$ and $\mathbf{\Sigma} = E[X^2] = 1/6$. To design the heterogenous data setting, we let $\boldsymbol{\beta}_0^{(j)} = j$ for $j = 1, 2, \ldots, s$ on the $j$-th subpopulation. The nonparametric function $f_0(z)$, which is common across all subpopulations, is assumed to be $0.6b_{30,17}(z) + 0.4b_{3,11}$, where $b_{a_1,a_2}$ is the density function for $Beta(a_1, a_2)$.

We start from the 95% predictive interval (at $(x_0, z_0)$) implied by the joint asymptotic distribution (4.3):

$$\left[\widehat{Y}^{(j)} \pm 1.96\sigma\sqrt{x_0^T\mathbf{\Omega}^{-1}x_0/n + \sigma_{z_0}^2/(N\lambda^{1/(2\nu)}) + 1}\right],$$

where $\widehat{Y}^{(j)} = x_0^T\widehat{\boldsymbol{\beta}}^{(j)} + \bar{f}(z_0)$ is the predicted response. The unknown error variance $\sigma$ is estimated by $(\widehat{\sigma}^{(j)})^2 = n^{-1}\sum_{i\in S_j}(Y_i - X_i^T\widehat{\boldsymbol{\beta}}^{(j)} - \widehat{f}^{(j)}(Z_i))^2/(n - Tr(A(\lambda)))$, where $A(\lambda)$ denotes the smoothing matrix, followed by an aggregation $\bar{\sigma}^2 = 1/s\sum_{j=1}^{s}(\widehat{\sigma}^{(j)})^2$. In the simulations, we fix $x_0 = 0.5$ and choose $z_0 = 0.25, 0.5, 0.75$ and $0.95$. The coverage probability is calculated based on 200 repetitions. As for $N$ and $s$, we set $N = 256, 528, 1024, 2048, 4096$, and choose $s = 2^0, 2^1, \ldots, 2^{t-3}$ when $N = 2^t$. The simulation results are summarized in Figure 1. We notice an interesting phase transition from Figure 1: when $s \leq s^*$ where $s^* \approx N^{0.45}$, the coverage probability is approximately 95%; when $s \geq s^*$, the coverage probability drastically decreases. This empirical observation is strongly supported by our theory developed in Section 4.3 where $s^* \approx N^{0.42}\log^{-6} N$ for $\nu = 2$.

We next compute the mean-squared errors of $\bar{f}$ under different choices of $N$ and $s$ in Figure 2. It is demonstrated that the increasing trends of MSE as $s$ increases are very similar for different $N$. More importantly, all the MSE curves suddenly blow up when $s \approx N^{0.4}$. This is also close to our theoretical result that the transition point is around $N^{0.45}\log^{-6} N$.

We next empirically verify the efficiency boosting theory developed in Section 3.5. Based on $\widehat{\boldsymbol{\beta}}^{(j)}$ and $\check{\boldsymbol{\beta}}^{(j)}$, we construct the following two types of 95% confidence intervals for $\boldsymbol{\beta}_0^{(j)}$:

$$\begin{aligned}\mathrm{CI}_1 &= \left[\widehat{\boldsymbol{\beta}}^{(j)} \pm 1.96\mathbf{\Omega}^{-1/2}n^{-1/2}\bar{\sigma}\right],\\ \mathrm{CI}_2 &= \left[\check{\boldsymbol{\beta}}^{(j)} \pm 1.96\mathbf{\Sigma}^{-1/2}n^{-1/2}\bar{\sigma}\right].\end{aligned}$$

Obviously, $\mathrm{CI}_2$ is shorter than $\mathrm{CI}_1$. However, Theorem 3.7 shows that $\mathrm{CI}_2$ is valid only when $s$ satisfies both a upper bound and a lower bound. This theoretical condition is empirically verified in Figure 3 which exhibits the validity range of $\mathrm{CI}_2$ in terms of $s$. In Figure 4, we further compare $\mathrm{CI}_2$ and $\mathrm{CI}_1$ in terms of their coverage probabilities and lengths. This figure shows that when $s$ is in a proper range, the coverage probabilities of $\mathrm{CI}_1$ and $\mathrm{CI}_2$ are similar, while $\mathrm{CI}_2$ is significantly shorter.

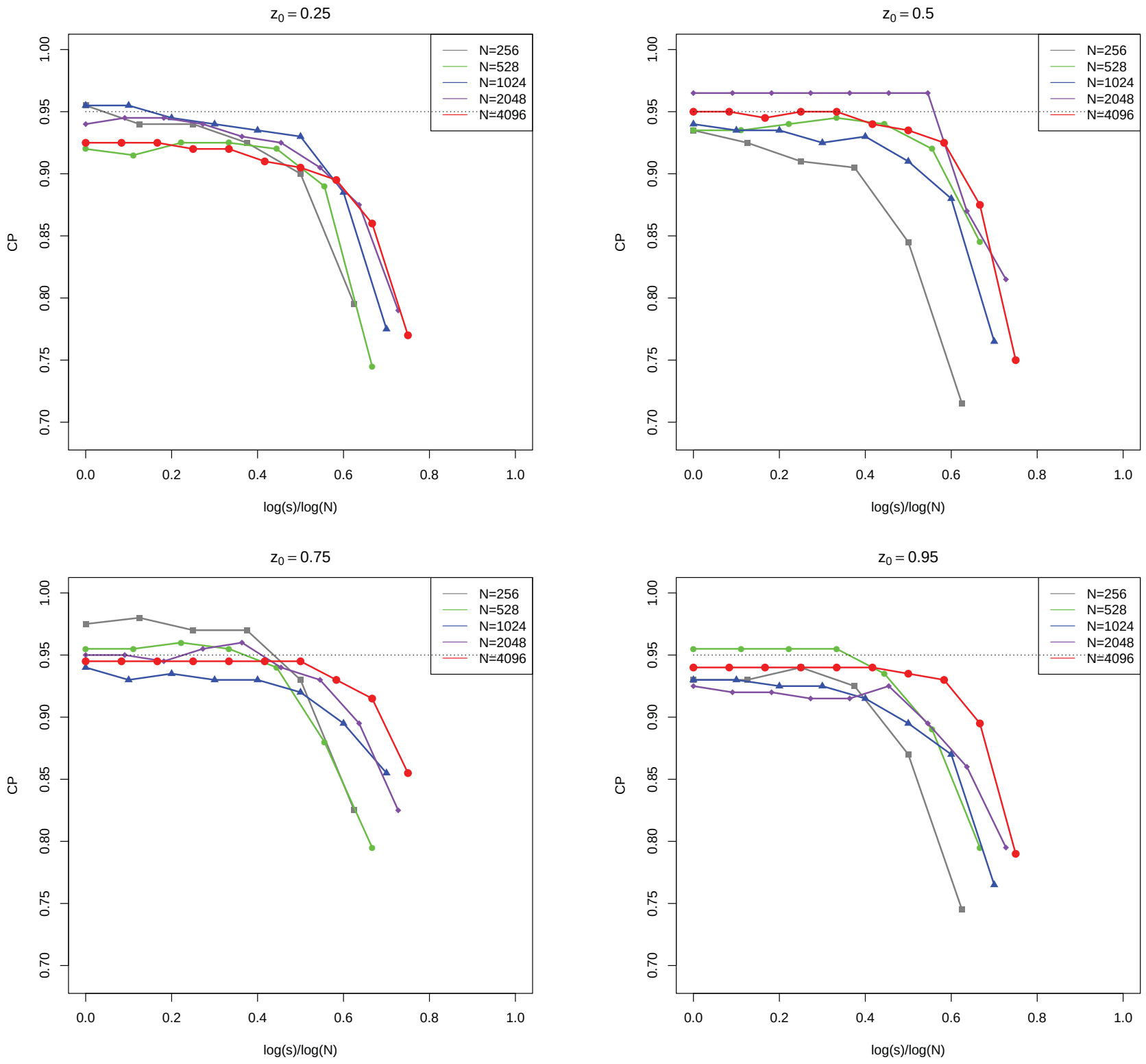

Fig 1. *Coverage probability of 95% predictive interval with different choices of $s$ and $N$*

Lastly, we consider the heterogeneity testing. In Figure 5, we compare tests $\Psi_1$ and $\Psi_2$ under different choices of $N$ and $s \geq 2$. Specifically, Figure 5 (i) compares the nominal levels, while Figure 5 (ii) - (iv) compare the powers under various alternative hypotheses $H_1 : \boldsymbol{\beta}_0^{(j)} - \boldsymbol{\beta}_0^{(k)} = \Delta$, where $\Delta = 0.5, 1, 1.5$. It is clearly seen that both tests are consistent, and their powers increase as $\Delta$ or $N$ increases. In addition, we observe that $\Psi_2$ has uniformly larger powers than $\Psi_1$.

## 7. Proof of Main Results.

In this section, we present main proofs of Lemma 3.1 and Theorems 3.2, 3.4 and 3.6 in the main text.

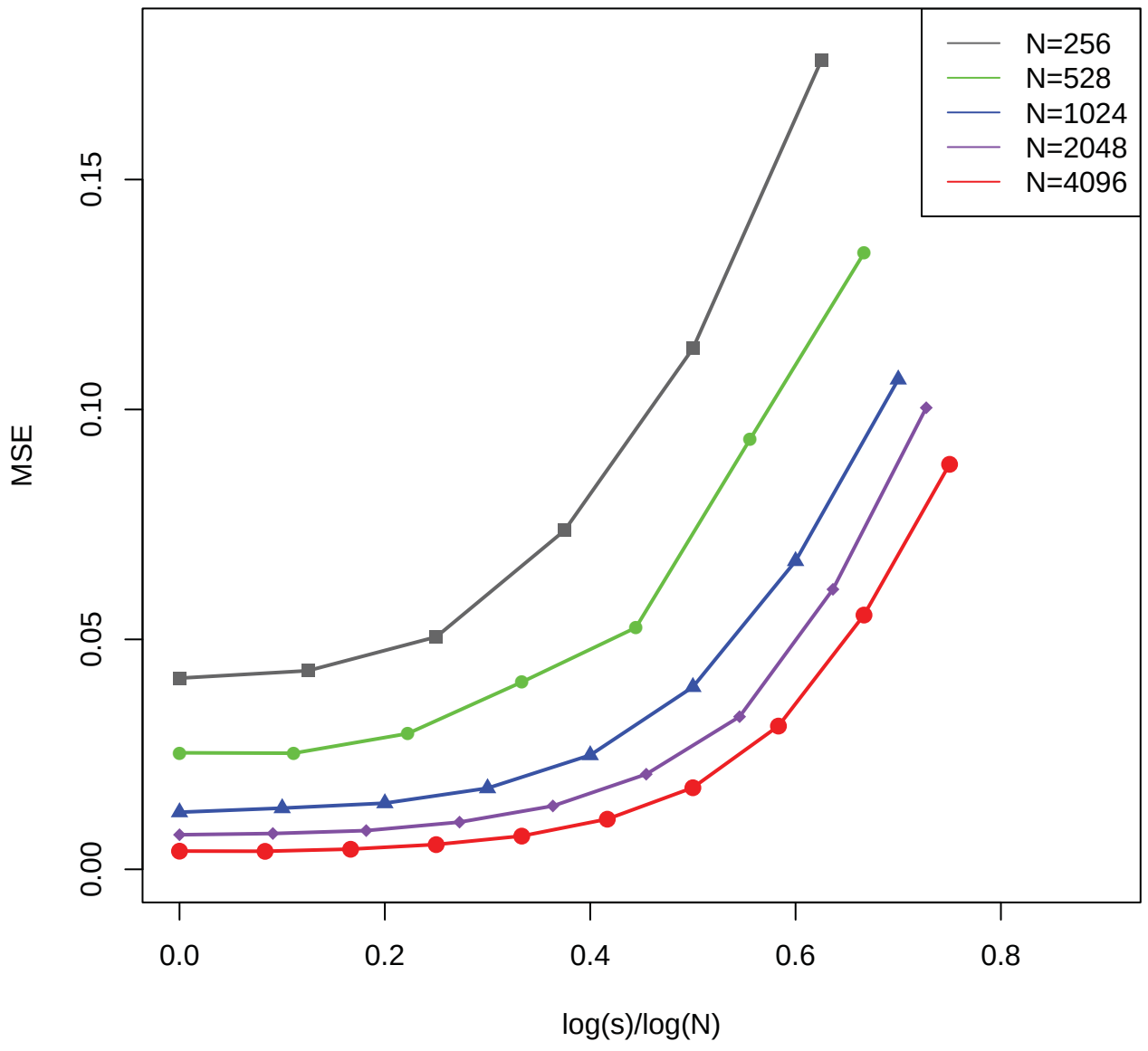

FIG 2. *Mean-square errors of $\bar{f}$ under different choices of $N$ and $s$*

7.1. *Proof of Lemma 3.1.*

PROOF. We start from analyzing the minimization problem (3.2) on each sub-population. Recall $m = (\boldsymbol{\beta}, f)$ and $U = (\boldsymbol{X}, Z)$. The objective function can be rewritten as

$$\frac{1}{n}\sum_{i\in S_j}(Y_i - \boldsymbol{X}_i^T\boldsymbol{\beta} - f(Z_i))^2 + \lambda\|f\|_{\mathcal{H}}^2 = \frac{1}{n}\sum_{i\in S_j}(Y_i - m(U_i))^2 + \langle P_\lambda m, m\rangle_{\mathcal{A}}$$
$$= \frac{1}{n}\sum_{i\in S_j}(Y_i - \langle R_{U_i}, m\rangle_{\mathcal{A}})^2 + \langle P_\lambda m, m\rangle_{\mathcal{A}}$$

The first order optimality condition (w.r.t. Fréchet derivative) gives

$$\frac{1}{n}\sum_{i\in S_j} R_{U_i}(\widehat{m}^{(j)}(U_i) - Y_i) + P_\lambda \widehat{m}^{(j)} = 0,$$

where $\widehat{m}^{(j)} = (\widehat{\boldsymbol{\beta}}^{(j)}, \widehat{f}^{(j)})$. This implies that

$$-\frac{1}{n}\sum_{i\in S_j} R_{U_i}\varepsilon_i + \frac{1}{n}\sum_{i\in S_j} R_{U_i}(\widehat{m}^{(j)}(U_i) - m_0^{(j)}(U_i)) + P_\lambda \widehat{m}^{(j)} = 0,$$

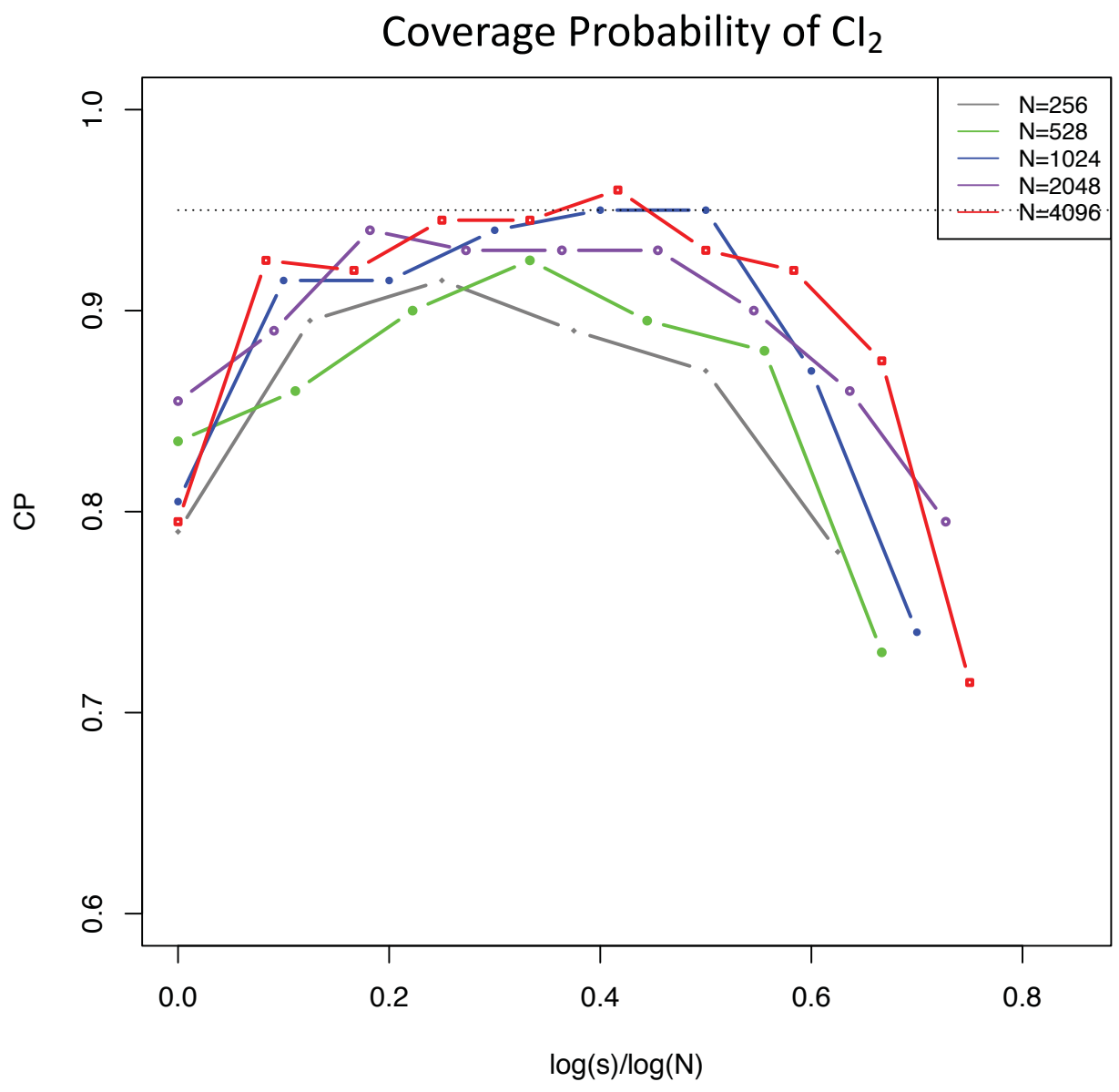

FIG 3. *Coverage probability of* 95% *confidence interval based on* $\check{\boldsymbol{\beta}}^{(j)}$

where $m_0^{(j)} = (\boldsymbol{\beta}_0^{(j)}, f_0)$. Define $\Delta m^{(j)} := \widehat{m}^{(j)} - m_0^{(j)}$. Adding $\mathbb{E}_U[R_U \Delta m^{(j)}(U)]$ on both sides of the above equation, we have

(7.1)

$$\mathbb{E}_U[R_U \Delta m^{(j)}(U)] + P_\lambda \Delta m^{(j)} \\ = \frac{1}{n} \sum_{i \in S_j} R_{U_i} \varepsilon_i - P_\lambda m_0^{(j)} - \frac{1}{n} \sum_{i \in S_j} \big( R_{U_i} \Delta m^{(j)}(U_i) - \mathbb{E}_U[R_U \Delta m^{(j)}(U)] \big).$$

The L.H.S. of (7.1) can be rewritten as

$$\begin{aligned} \mathbb{E}_U[R_U \Delta m^{(j)}(U)] + P_\lambda \Delta m^{(j)} &= \mathbb{E}_U[R_U \langle R_U, \Delta m^{(j)} \rangle_{\mathcal{A}}] + P_\lambda \Delta m^{(j)} \\ &= \big( \mathbb{E}_U[R_U \otimes R_U] + P_\lambda \big) \Delta m^{(j)} \\ &= \Delta m^{(j)}, \end{aligned}$$

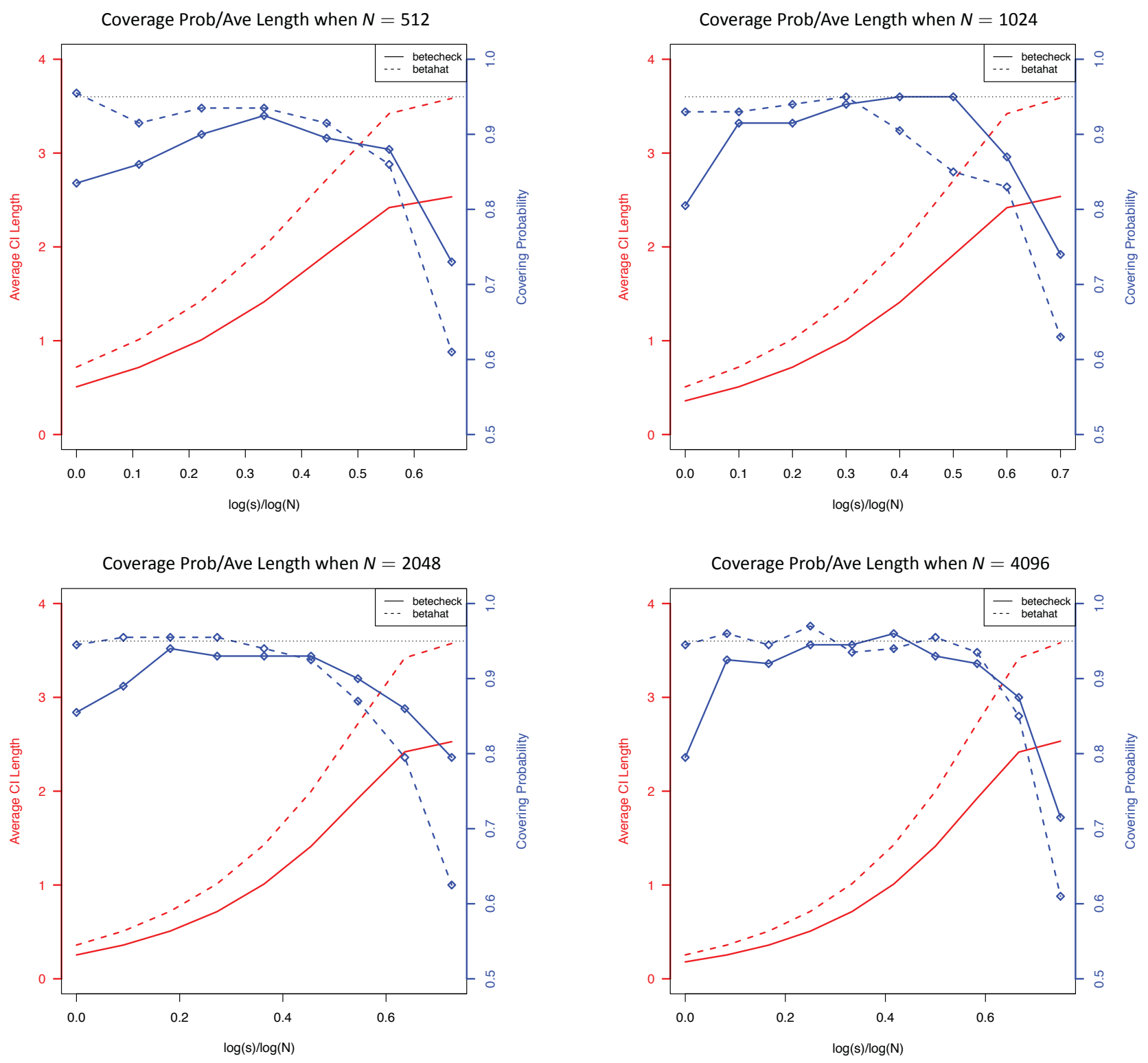

FIG 4. *Coverage probabilities and average lengths of* 95% *for two types of confidence intervals. In the above figures, dashed lines represent* $CI_1$*, which is constructed based on* $\widehat{\boldsymbol{\beta}}^{(j)}$*, and solid lines represent* $CI_2$*, which is constructed based on* $\check{\boldsymbol{\beta}}^{(j)}$.

where the last equality follows from proposition 2.2. Then (7.1) becomes

$$\widehat{m}^{(j)} - m_0^{(j)} = \frac{1}{n}\sum_{i\in S_j} R_{U_i}\varepsilon_i - P_\lambda m_0^{(j)} - \frac{1}{n}\sum_{i\in S_j}\Big(R_{U_i}\Delta m^{(j)}(U_i) - \mathbb{E}_U[R_U\Delta m^{(j)}(U)]\Big). \tag{7.2}$$

We denote the last term in the R.H.S. of (7.2) as $Rem^{(j)} := \frac{1}{n}\sum_{i\in S_j}\big(R_{U_i}\Delta m^{(j)}(U_i) - \mathbb{E}_U[R_U\Delta m^{(j)}(U)]\big)$. Recall that $R_u = (L_u, N_u)$ and $P_\lambda m_0^{(j)} = (L_\lambda f_0, N_\lambda f_0)$.

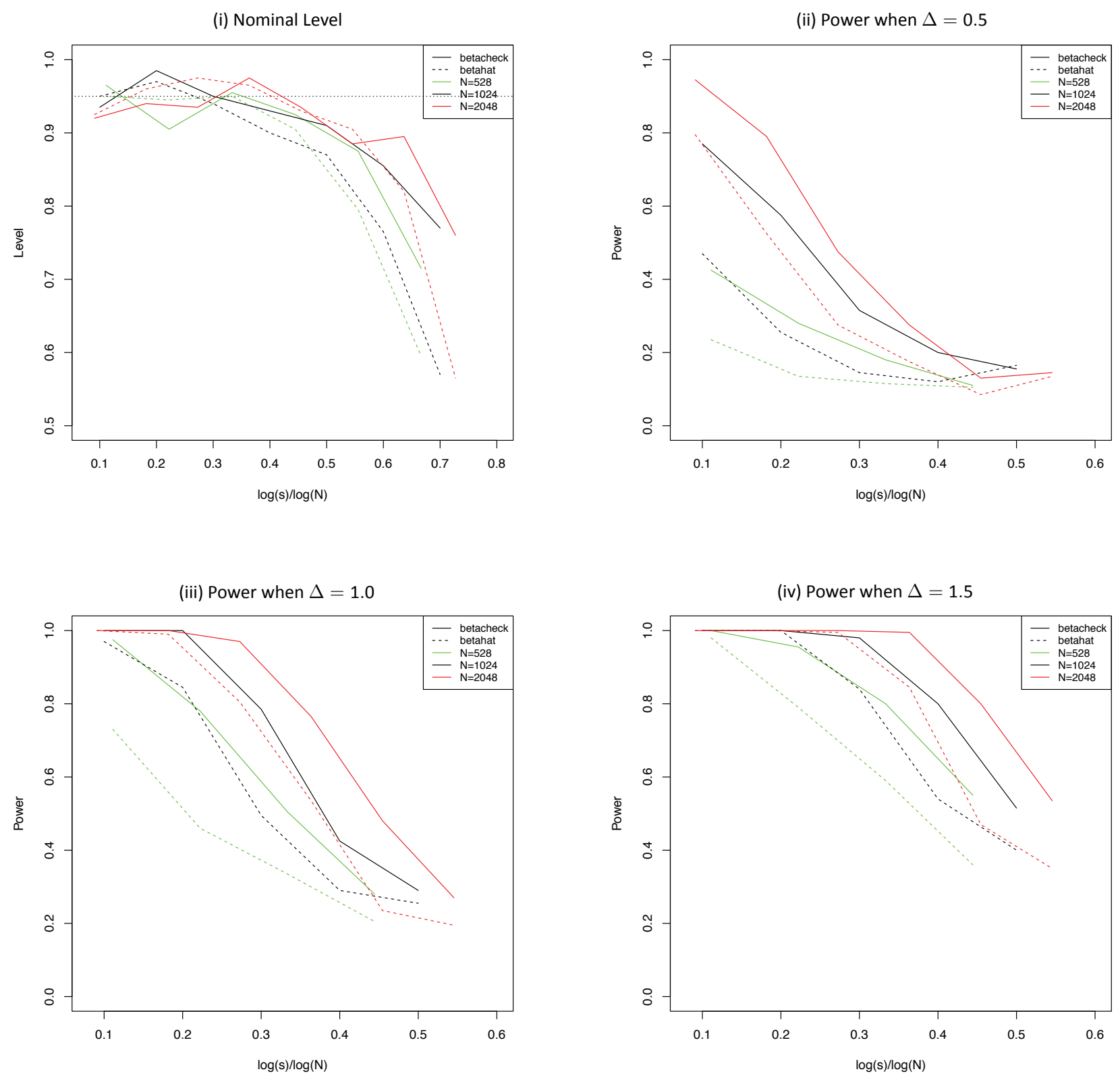

FIG 5. *(i) Nominal level of heterogeneity tests* $\Psi_1$ *and* $\Psi_2$*; (ii) - (iv) Power of heterogeneity tests* $\Psi_1$ *and* $\Psi_2$ *when* $\Delta = 0.5, 1.0, 1.5$. *In the above figures, dashed lines represent* $\Psi_1$*, which is constructed based on* $\check{\boldsymbol{\beta}}$*, and solid lines represent* $\Psi_2$*, which is constructed based on* $\widehat{\boldsymbol{\beta}}$.

Thus the above remainder term decomposes into two components:

$$Rem_\beta^{(j)} := \frac{1}{n}\sum_{i\in S_j}\Big(L_{U_i}\Delta m^{(j)}(U_i) - \mathbb{E}_U[L_U \Delta m^{(j)}(U)]\Big)$$

$$Rem_f^{(j)} := \frac{1}{n}\sum_{i\in S_j}\Big(N_{U_i}\Delta m^{(j)}(U_i) - \mathbb{E}_U[N_U \Delta m^{(j)}(U)]\Big).$$

Therefore, (7.2) can be rewritten into Equations (3.5) and (3.6) for all $j = 1, \ldots, s$. This completes the proof of the first part of Lemma 3.1. Taking average of (3.6) for all $j$ over $s$, and by definition of $\bar{f}$, we have

$$\bar{f} - f_0 = \frac{1}{N}\sum_{i=1}^{N} N_{U_i}\varepsilon_i - N_\lambda f_0 - \frac{1}{s}\sum_{j=1}^{s} Rem_f^{(j)}, \tag{7.3}$$

where we used $1/N\sum_{i=1}^{N} N_{U_i}\varepsilon_i = 1/s\sum_{j=1}^{s} 1/n\sum_{i\in S_j} N_{U_i}\varepsilon_i$. Eq. (3.5) and (7.3) are the basic equalities to derive the finite sample rate of convergence and joint limit distribution of $\widehat{\boldsymbol{\beta}}^{(j)}$ and $\bar{f}$. To this end, we need to control the remainder terms in the above two equalities, which is the second part of Lemma 3.1. We delegate the proofs to the following two lemmas:

**Lemma 7.1.** Suppose the conditions in Lemma 3.1 hold. We have for all $j = 1, \ldots, s$

$$\mathbb{E}\big[\|Rem_\beta^{(j)}\|_2^2\big] \le a(n, \lambda, J),$$

for sufficiently large $n$, where $a(n, \lambda, J)$ is as defined in Lemma 3.1. Moreover, the inequality also holds for $\mathbb{E}\big[\|Rem_f^{(j)}\|_{\mathcal{C}}^2\big]$.

**Lemma 7.2.** Suppose the conditions in Lemma 3.1 hold. We have the following two sets of results that control the remainder terms:

(i) For all $j = 1, \ldots, s$, it holds that

$$\mathbb{P}\Big(\|Rem_\beta^{(j)}\|_2 \ge b(n, \lambda, J)\Big) \lesssim n\exp(-c\log^2 n), \tag{7.4}$$

where $b(n, \lambda, J)$ is as defined in Lemma 3.1.

(ii) In addition, we have

$$\Big\|\frac{1}{s}\sum_{j=1}^{s} Rem_\beta^{(j)}\Big\|_2 = o_P(s^{-1/2} b(n, \lambda, J)\log N). \tag{7.5}$$

Furthermore, (7.4) and (7.5) also hold if $\|Rem_\beta^{(j)}\|_2$ and $\big\|s^{-1}\sum_{j=1}^{s} Rem_\beta^{(j)}\big\|_2$ are replaced by $\|Rem_f^{(j)}\|_{\mathcal{C}}$ and $\big\|s^{-1}\sum_{j=1}^{s} Rem_f^{(j)}\big\|_{\mathcal{C}}$.

By the above two lemmas, we complete the second part of Lemma 3.1. □

7.2. *Proof of Theorem 3.2.*

Proof. By (7.3), it follows that

$$\mathbb{E}[\|\bar{f} - f_0\|_{\mathcal{C}}^2] \le 3\mathbb{E}\Big[\Big\|\frac{1}{N}\sum_{i=1}^{N} N_{U_i}\varepsilon_i\Big\|_{\mathcal{C}}^2\Big] + 3\|N_\lambda f_0\|_{\mathcal{C}}^2 + 3\mathbb{E}\Big[\Big\|\frac{1}{s}\sum_{j=1}^{s} Rem_f^{(j)}\Big\|_{\mathcal{C}}^2\Big]. \tag{7.6}$$

By Lemma A.5 and the fact that each $N_{U_i}\varepsilon_i$ is i.i.d., it follows that

$$\mathbb{E}\left[\left\|\frac{1}{N}\sum_{i=1}^{N} N_{U_i}\varepsilon_i\right\|_{\mathcal{C}}^2\right] = \frac{1}{N}\mathbb{E}\big[\|N_U\varepsilon\|_{\mathcal{C}}^2\big] \le C_1\sigma^2\frac{d(\lambda)}{N}, \tag{7.7}$$

and

$$\|N_\lambda f_0\|_{\mathcal{C}}^2 \le 2\|f_0\|_{\mathcal{H}}^2\lambda + C_2\lambda^2, \tag{7.8}$$

where $C_1$ and $C_2$ are constants specified in Lemma A.5. As for the third term in (7.6), we have by independence across sub-populations that

$$\mathbb{E}\left[\left\|\frac{1}{s}\sum_{j=1}^{s} Rem_f^{(j)}\right\|_{\mathcal{C}}^2\right] = \frac{1}{s^2}\sum_{j=1}^{s}\mathbb{E}\big[\|Rem_f^{(j)}\|_{\mathcal{C}}^2\big]. \tag{7.9}$$

Combining (7.6) - (7.9) and Lemma 7.1, and by the fact that $\|\bar{f}-f_0\|_{L_2(\mathbb{P}_Z)}^2 \le \|\bar{f}-f_0\|_{\mathcal{C}}^2$, we complete the proof of Theorem 3.2. □

7.3. *Proof of Theorem 3.4.*

PROOF. Recall that $m_0^{(j)*} = (\boldsymbol{\beta}_0^{(j)*}, f_0^*) = (id - P_\lambda)m_0^{(j)}$ where $m_0^{(j)} = (\boldsymbol{\beta}_0^{(j)}, f_0)$. This implies that $\boldsymbol{\beta}_0^{(j)*} = \boldsymbol{\beta}_0^{(j)} - L_\lambda f_0$ and $f_0^* = f_0 - N_\lambda f_0$. By (3.5) and (7.3), for arbitrary $\boldsymbol{x}$ and $z_0$,

$$\begin{aligned}
&(\boldsymbol{x}^T, 1)\begin{pmatrix} \sqrt{n}(\widehat{\boldsymbol{\beta}}^{(j)} - \boldsymbol{\beta}_0^{(j)*}) \\ \sqrt{N/d(\lambda)}\big(\bar{f}(z_0) - f_0^*(z_0)\big)\end{pmatrix} \\
&= \sqrt{n}\boldsymbol{x}^T(\widehat{\boldsymbol{\beta}}^{(j)} - \boldsymbol{\beta}_0^{(j)*}) + \sqrt{N/d(\lambda)}(\bar{f}_{N,\lambda}(z_0) - f_0^*(z_0)) \\
&= \underbrace{\frac{1}{\sqrt{n}}\sum_{i\in S_j}\boldsymbol{x}^T L_{U_i}\varepsilon_i + \frac{1}{\sqrt{N}}\sum_{i=1}^{N} d(\lambda)^{-1/2}N_{U_i}(z_0)\varepsilon_i}_{(I)} \\
&\qquad \underbrace{+\sqrt{n}\boldsymbol{x}^T Rem_\beta^{(j)} + \sqrt{N/d(\lambda)}s^{-1}\sum_{j=1}^{s} Rem_f^{(j)}(z_0)}_{(II)}.
\end{aligned}$$

In what follows, we will show that the main term (I) is asymptotically normal and the remainder term (II) is of order $o_P(1)$. Given that $\boldsymbol{x}$ is arbitrary, we apply Wold device to conclude the proof of joint asymptotic normality.

**Asymptotic normality of (I)**: We present the result for showing asymptotic normality of (I) in the following lemma and defer its proof to supplemental material (Zhao et al., 2015).

**Lemma 7.3.** Suppose Assumptions 3.1, 3.2 hold and that $\|\widetilde{K}_{z_0}\|_{L_2(\mathbb{P}_Z)}/d(\lambda)^{1/2} \to \sigma_{z_0}$, $(W_\lambda \boldsymbol{A})(z_0)/d(\lambda)^{1/2} \to \boldsymbol{\alpha}_{z_0} \in \mathbb{R}^p$, and $\boldsymbol{A}(z_0)/d(\lambda)^{1/2} \to -\boldsymbol{\gamma}_{z_0} \in \mathbb{R}^p$ as $N \to \infty$. We have

(i) if $s \to \infty$, then

$$(I) \rightsquigarrow N(0, \sigma^2(\boldsymbol{x}^T \boldsymbol{\Omega}^{-1}\boldsymbol{x} + \Sigma_{22})). \tag{7.10}$$

(ii) if $s$ is fixed, then

$$(I) \rightsquigarrow N(0, \sigma^2(\boldsymbol{x}^T \boldsymbol{\Omega}^{-1}\boldsymbol{x} + \Sigma_{22} + 2s^{-1/2}\boldsymbol{x}^T \boldsymbol{\Sigma}_{12})). \tag{7.11}$$

**Control of the remainder term (II)**: We now turn to bound the remainder term (II). We can show that if (3.12) holds, then $d(\lambda)n^{-1/2}(J(\mathcal{F},1)+\log n) = o(1)$. Hence by Lemma 7.2, we have

$$\begin{aligned} \sqrt{n}|\boldsymbol{x}^T Rem_\beta^{(j)}| &\le \sqrt{n}\|\boldsymbol{x}\|_2 \|Rem_\beta^{(j)}\|_2 \\ &= o_P(n^{1/2} b(n,\lambda,J)) = o_P(\sqrt{N} s^{-1/2} b(n,\lambda,J)), \end{aligned} \tag{7.12}$$

where we used the boundedness of $\boldsymbol{x}$. Also,

$$\begin{aligned} \sqrt{N/d(\lambda)}\Big|s^{-1}\sum_{j=1}^s Rem_f^{(j)}(z_0)\Big| &\le \sqrt{N/d(\lambda)}\|\widetilde{K}_{z_0}\|_{\mathcal{C}} \Big\|s^{-1}\sum_{j=1}^s Rem_f^{(j)}\Big\|_{\mathcal{C}} \\ &\lesssim \sqrt{N}\Big\|s^{-1}\sum_{j=1}^s Rem_f^{(j)}\Big\|_{\mathcal{C}} \\ &= o_P(\sqrt{N} s^{-1/2} b(n,\lambda,J)\log N), \end{aligned} \tag{7.13}$$

where the second inequality follows from Lemma A.4. Therefore by (7.12) and (7.13), we have

$$(II) = o_P(\sqrt{N} s^{-1/2} b(n,\lambda,J)\log N). \tag{7.14}$$

Now by definition of $b(n,\lambda,J)$ and condition (3.13), we have $(II) = o_P(1)$.

Combining (7.10) and (7.14), it follows that if $s \to \infty$, then

$$(\boldsymbol{x}^T, 1)\begin{pmatrix} \sqrt{n}(\widehat{\boldsymbol{\beta}}^{(j)} - \boldsymbol{\beta}_0^{(j)*}) \\ \sqrt{N/d(\lambda)}\big(\bar{f}(z_0) - f_0^*(z_0)\big) \end{pmatrix} \rightsquigarrow N\big(0, \sigma^2(\boldsymbol{x}^T\boldsymbol{\Omega}^{-1}\boldsymbol{x} + \Sigma_{22})\big).$$

Combining (7.11) and (7.14), it follows that if $s$ is fixed, then

$$\begin{aligned} &(\boldsymbol{x}^T, 1)\begin{pmatrix} \sqrt{n}(\widehat{\boldsymbol{\beta}}^{(j)} - \boldsymbol{\beta}_0^{(j)*}) \\ \sqrt{N/d(\lambda)}\big(\bar{f}(z_0) - f_0^*(z_0)\big) \end{pmatrix} \\ &\qquad \rightsquigarrow N\big(0, \sigma^2(\boldsymbol{x}^T\boldsymbol{\Omega}^{-1}\boldsymbol{x} + \Sigma_{22} + 2s^{-1/2}\boldsymbol{x}^T\boldsymbol{\Sigma}_{12})\big). \end{aligned}$$

By the arbitrariness of $\boldsymbol{x}$, we reach the conclusion of the theorem using Wold device. □

7.4. *Proof of Lemma 7.2: Controlling the Remainder Term.*

PROOF. (i) Recall that $Rem^{(j)} = (Rem_\beta^{(j)}, Rem_f^{(j)}) \in \mathcal{A}$. We first derive the bound of $\|Rem^{(j)}\|_{\mathcal{A}}$. Recall

$$Rem^{(j)} = \frac{1}{n}\sum_{i\in S_j} \Delta m^{(j)}(U_i)R_{U_i} - \mathbb{E}_U[\Delta m^{(j)}(U)R_U].$$

Let $Z_n(m) = c_r^{-1}d(\lambda)^{-1/2}n^{-1/2}\sum_{i\in S_j}\big\{m(U_i)R_{U_i} - \mathbb{E}[m(U)R_U]\big\}$, where $c_r$ is the constant specified in Lemma A.4. Note that $Z_n(m)$ is implicitly related to $j$ but we omit the superscript of $(j)$. We have $Rem^{(j)} = c_r^{-1}\sqrt{n/d(\lambda)}Z_n(\Delta m^{(j)})$. We apply Lemma F.1 to obtain an exponential inequality for $\sup_{m\in\mathcal{F}}\|Z_n(m)\|_{\mathcal{A}}$. The first step is to show that $Z_n(m)$ is a sub-Gaussian process by Lemma G.1. Let $g(U_i, m) = c_r^{-1}\sqrt{n/d(\lambda)}\big(m(U_i)R_{U_i} - \mathbb{E}[m(U)R_U]\big)$. Now for any $m_1$ and $m_2$,

$$\begin{aligned}
&\|g(U_i, m_1) - g(U_i, m_2)\|_{\mathcal{A}}\\
&= c_r^{-1}\sqrt{n/d(\lambda)}\Big\{\|\big(m_1(U_i) - m_2(U_i)\big)R_{U_i}\|_{\mathcal{A}} + \|\mathbb{E}\big[\big(m_1(U) - m_2(U)\big)R_U\big]\|_{\mathcal{A}}\Big\}\\
&\le 2\sqrt{n}\|m_1 - m_2\|_{\sup},
\end{aligned}$$

where we used the fact that $\|R_u\|_{\mathcal{A}} \le c_r d(\lambda)^{1/2}$ by Lemma A.4. Note that $Z_n(m) = \frac{1}{n}\sum_{i\in S_j} g(U_i, m)$.Therefore by Lemma G.1, we have for any $t > 0$,

$$\begin{aligned}
\mathbb{P}\Big(\|Z_n(m_1) - Z_n(m_2)\|_{\mathcal{A}} \ge t\Big) &= \mathbb{P}\bigg(\Big\|\frac{1}{n}\sum_{i=1}^n \big\{g(U_i, m_1) - g(U_i, m_2)\big\}\Big\|_{\mathcal{A}} \ge t\bigg)\\
&\le 2\exp\bigg(-\frac{t^2}{8\|m_1 - m_2\|_{\sup}^2}\bigg) \qquad (7.15)
\end{aligned}$$

Then by Lemma F.1, we have

$$\mathbb{P}\Big(\sup_{m\in\mathcal{F}}\|Z_n(m)\|_{\mathcal{A}} \ge CJ(\mathcal{F}, \mathrm{diam}(\mathcal{F})) + x\Big) \le C\exp\bigg(\frac{-x^2}{C\mathrm{diam}(\mathcal{F})^2}\bigg), \qquad (7.16)$$

where $\mathrm{diam}(\mathcal{F}) = \sup_{m_1,m_2\in\mathcal{F}}\|m_1 - m_2\|_{\sup}$.

Define $q_{n,\lambda} = c_r r_{n,\lambda} d(\lambda)^{1/2}$ and $\widetilde{m} = q_{n,\lambda}^{-1}\Delta m^{(j)}/2$. Again we do not specify its relationship with $j$. Define the event $\mathcal{E} = \big\{\|\Delta m^{(j)}\|_{\mathcal{A}} \le r_{n,\lambda}\big\}$. On the event $\mathcal{E}$, we have

$$\|\widetilde{m}\|_{\sup} \le c_r d(\lambda)^{1/2}(2q_{n,\lambda})^{-1}\|\Delta m^{(j)}\|_{\mathcal{A}} \le 1/2,$$

where we used the fact that $\|\widetilde{m}\|_{\sup} \le c_r d(\lambda)^{1/2}\|\widetilde{m}\|_{\mathcal{A}}$ by Lemma A.4. This implies $|\boldsymbol{x}^T\widetilde{\boldsymbol{\beta}} + \widetilde{f}(z)| \le 1/2$ for any $(\boldsymbol{x}, z)$. Letting $\boldsymbol{x} = 0$, one gets $\|\widetilde{f}\|_{\sup} \le 1/2$, which further implies $|\boldsymbol{x}^T\widetilde{\boldsymbol{\beta}}| \le 1$ for all $\boldsymbol{x}$ by triangular inequality.

Moreover, on the even $\mathcal{E}$ we have

$$\|\widetilde{f}\|_{\mathcal{H}} \le \lambda^{-1/2}\|\widetilde{m}\|_{\mathcal{A}} \le \lambda^{-1/2}/(2q_{n,\lambda})\|\Delta m^{(j)}\|_{\mathcal{A}} \le c_r^{-1} d(\lambda)^{-1/2}\lambda^{-1/2}$$

by the definition of $\|\cdot\|_{\mathcal{A}}$. Hence, we have shown that $\mathcal{E} \subset \{\widetilde{m} \in \mathcal{F}\}$. Combining this fact with (7.16), and noting that $\mathrm{diam}(\mathcal{F}) \le 1$, we have

$$\mathbb{P}\Big(\big\{\|Z_n(\widetilde{m})\|_{\mathcal{A}} \ge CJ(\mathcal{F},1)+x\big\} \cap \mathcal{E}\Big) \le C\exp\big(-x^2/C\big), \tag{7.17}$$

by Lemma F.1. Using the definition of $\widetilde{m}$, and the relationship that $Rem^{(j)} = c_r^{-1}\sqrt{n/d(\lambda)}Z_n(\Delta m^{(j)})$, we calculate that

$$Z_n(\widetilde{m}) = (1/2)d(\lambda)^{-1/2}n^{1/2}q_{n,\lambda}^{-1}Rem^{(j)} = (1/2)c_r^{-1}d(\lambda)^{-1}n^{1/2}r_{n,\lambda}^{-1}Rem^{(j)}.$$

Plugging the above form of $Z_n(\widetilde{m})$ into (7.17) and letting $x = \log n$ in (7.17), we have

$$\mathbb{P}\Big(\big\{\|Rem^{(j)}\|_{\mathcal{A}} \ge b(n,\lambda,J)\big\} \cap \mathcal{E}\Big) \le C\exp\big(-\log^2 n/C\big), \tag{7.18}$$

where we used the definition that $b(n,\lambda,J) = Cd(\lambda)n^{-1/2}r_{n,\lambda}(J(\mathcal{F},1)+\log n)$. Therefore we have

$$\begin{aligned}\mathbb{P}\Big(\|Rem^{(j)}\|_{\mathcal{A}} \ge b(n,\lambda,J)\Big) &\le \mathbb{P}\Big(\big\{\|Rem^{(j)}\|_{\mathcal{A}} \ge b(n,\lambda,J)\big\} \cap \mathcal{E}\Big) + \mathbb{P}(\mathcal{E}^c)\\ &\le C\exp\big(-\log^2 n/C\big) + \mathbb{P}(\mathcal{E}^c).\end{aligned} \tag{7.19}$$

We have the following lemma that controls $\mathbb{P}(\mathcal{E}^c)$.

**Lemma 7.4.** Suppose the conditions in Lemma 3.1 hold. There exist a constant $c$ such that

$$\mathbb{P}(\mathcal{E}^c) = \mathbb{P}\Big(\|\Delta m^{(j)}\|_{\mathcal{A}} \ge r_{n,\lambda}\Big) \lesssim n\exp(-c\log^2 n).$$

for all $j = 1, \dots, s$.

By Lemma 7.4 and (7.19) we have

$$\mathbb{P}\Big(\|Rem^{(j)}\|_{\mathcal{A}} \ge b(n,\lambda,J)\Big) \lesssim n\exp(-c\log^2 n). \tag{7.20}$$

We can apply similar arguments as above to bound $\|Rem_f^{(j)}\|_{\mathcal{C}}$, by changing $\omega(\mathcal{F},1)$ to $\omega(\mathcal{F}_2,1)$, which is dominated by $\omega(\mathcal{F},1)$. The bound of $\|Rem_\beta^{(j)}\|_2$ then follows from triangular inequality.

(ii) We will use an Azuma-type inequality in Hilbert space to control the averaging remainder term $s^{-1}\sum_{j=1}^s Rem^{(j)}$, as all $Rem^{(j)}$ are independent and have zero mean. Define the event $\mathcal{A}_j = \big\{\|Rem^{(j)}\|_{\mathcal{A}} \le b(n,\lambda,J)\big\}$. By

Lemma G.1, we have

$$\mathbb{P}\Big(\Big\{\cap_j \mathcal{A}_j\Big\} \cap \Big\{\|s^{-1}\sum_{j=1}^{s} Rem^{(j)}\|_{\mathcal{A}} > s^{-1/2} b(n,\lambda,J)\log N\Big\}\Big) \le 2\exp(-\log^2 N/2). \tag{7.21}$$

Moreover, by (7.20),

$$\mathbb{P}(\mathcal{A}_j^c) \lesssim n\exp(-c\log^2 n). \tag{7.22}$$

Hence it follows that

$$\begin{aligned}
&\mathbb{P}\Big(\|s^{-1}\sum_{j=1}^{s} Rem^{(j)}\|_{\mathcal{A}} > s^{-1/2} b(n,\lambda,J)\log N\Big)\\
&\quad \le \mathbb{P}\Big(\Big\{\cap_{j=1}^{s} \mathcal{A}_j\Big\} \cap \Big\{\|s^{-1}\sum_{j=1}^{s} Rem^{(j)}\|_{\mathcal{A}} > s^{-1/2} b(n,\lambda,J)\log N\Big\}\Big) + \mathbb{P}(\cup_j \mathcal{A}_j^c)\\
&\quad \lesssim 2\exp(-\log^2 N/2) + ns\exp(-c\log^2 n) \lesssim N\exp(-c\log^2 n),
\end{aligned}$$

where the second inequality follows from (7.21), (7.22) and union bound. By our technical assumption that $s \lesssim N^{\psi}$ (stated before Assumption 3.1), we have $N\exp(-c\log^2 n) \asymp N\exp(-c'\log^2 N) \to 0$ as $N \to \infty$. This completes the proof of Part (ii).

Applying similar arguments as in (i), we get the similar inequalities for $\|1/s\sum_{j=1}^{s} Rem_{\beta}^{(j)}\|_2$ and $\|1/s\sum_{j=1}^{s} Rem_{f}^{(j)}\|_{\mathcal{C}}$. □

### 7.5. *Proof of Theorem 3.6.*

PROOF. In view of Theorem 3.6, we first prove

$$\begin{pmatrix} \sqrt{n}(\boldsymbol{\beta}_0^{(j)*} - \boldsymbol{\beta}_0^{(j)}) \\ \sqrt{N/d(\lambda)}(f_0^*(z_0) - f_0(z_0) - W_\lambda f_0(z_0)) \end{pmatrix} \to \mathbf{0} \tag{7.23}$$

for both (i) and (ii). By Proposition 2.3, we have

$$\begin{pmatrix} \boldsymbol{\beta}_0^{(j)*} - \boldsymbol{\beta}_0^{(j)} \\ f_0^*(z_0) - f_0(z_0) \end{pmatrix} = \begin{pmatrix} L_\lambda f_0 \\ W_\lambda f_0(z_0) + \boldsymbol{A}(z_0)^T L_\lambda f_0 \end{pmatrix}. \tag{7.24}$$

By Lemma A.5, it follows that under Assumption 3.3, $\|L_\lambda f_0\|_2 \lesssim \lambda$. Now we turn to $f_0^*(z_0) - f_0(z_0)$. Observe that

$$\boldsymbol{A}(z) = \langle \boldsymbol{A}, \widetilde{K}_z\rangle_{\mathcal{C}} = \langle \boldsymbol{B}, \widetilde{K}_z\rangle_{L_2(\mathbb{P}_Z)} = \sum_{\ell=1}^{\infty} \frac{\langle \boldsymbol{B}, \phi_\ell\rangle_{L_2(\mathbb{P}_Z)}}{1+\lambda/\mu_\ell}\phi_\ell(z), \tag{7.25}$$

Applying Cauchy-Schwarz, we obtain

$$
\begin{aligned}
A_k(z_0)^2 &\le \Big(\sum_{\ell=1}^{\infty}\frac{\langle B_k,\phi_\ell\rangle_{L_2(\mathbb{P}_Z)}^2}{\mu_\ell}\phi_\ell^2(z_0)\Big)\Big(\sum_{\ell=1}^{\infty}\frac{\mu_\ell}{(1+\lambda/\mu_\ell)^2}\Big)\\
&\le c_\phi^2\|B_k\|_{\mathcal{H}}^2\,\mathrm{Tr}(K),
\end{aligned}
$$

where the last inequality follows from the uniform boundedness of $\phi_\ell$. Hence we have that $A_k(z_0)$ is uniformly bounded, which implies

$$
\boldsymbol{A}(z_0)^T L_\lambda f_0 \le \|\boldsymbol{A}(z_0)\|_2\|L_\lambda f_0\|_2 \lesssim \lambda.
$$

Therefore, if we choose $\lambda = o\big(\sqrt{d(\lambda)/N}\wedge n^{-1/2}\big)$, then we get (7.23), which eliminates the estimation bias for $\boldsymbol{\beta}_0^{(j)}$.

Now we consider the asymptotic variance for cases (i) and (ii). It suffices to show that $\boldsymbol{\alpha}_{z_0} = \mathbf{0}$ under Assumption 3.3. Recall that $\boldsymbol{\alpha}_{z_0} = \lim_{N\to\infty} d(\lambda)^{-1/2}W_\lambda \boldsymbol{A}(z_0)$. By Lemma A.2 and (7.25), we have

$$
\begin{aligned}
W_\lambda A_k(z_0) &= \sum_{\ell=1}^{\infty}\frac{\langle B_k,\phi_\ell\rangle_{L_2(\mathbb{P}_Z)}}{1+\lambda/\mu_\ell}\frac{\lambda}{\lambda+\mu_\ell}\phi_\ell(z_0)\\
&\le \Big(\sum_{\ell=1}^{\infty}\frac{\langle B_k,\phi_\ell\rangle_{L_2(\mathbb{P}_Z)}^2}{\mu_\ell}\phi_\ell^2(z_0)\Big)\Big(\sum_{\ell=1}^{\infty}\frac{\mu_\ell}{(1+\lambda/\mu_\ell)^2}\Big)\\
&\le c_\phi^2\|B_k\|_{\mathcal{H}}^2\,\mathrm{Tr}(K).
\end{aligned}
$$

Hence by dominated conference theorem, as $\lambda\to 0$ we have $W_\lambda A_k(z_0)\to 0$. As $d(\lambda)^{-1} = O(1)$, it follows that $\boldsymbol{\alpha}_{z_0} = \lim_{N\to\infty} d(\lambda)^{-1/2}W_\lambda\boldsymbol{A}(z_0) = \mathbf{0}$.

When $d(\lambda)\to\infty$, we have $\boldsymbol{\gamma}_{z_0} = -\lim_{N\to\infty}\boldsymbol{A}(z_0)/d(\lambda)^{1/2} = \mathbf{0}$, as $A_k(z_0)$ is uniformly bounded. Hence $\boldsymbol{\Sigma}_{12}^* = \boldsymbol{\Sigma}_{21}^* = \mathbf{0}$ and $\Sigma_{22}^* = \sigma^2\sigma_{z_0}^2$. □

**Acknowledgements.** We thank Co-Editor Runze Li, an Associate Editor, and two referees for helpful comments that lead to important improvements on the paper.

## SUPPLEMENTARY MATERIAL

**Supplementary material for: A Partially Linear Framework for Massive Heterogenous Data** (DOI: To Be Assigned; .pdf). We provide the detailed proofs in the supplement.

## References.

Department of operations research
and financial engineering
Princeton University
Princeton, New Jersey 08544
USA
E-mail: tianqi@princeton.edu
hanliu@princeton.edu

Department of Statistics
Purdue University
West Lafayette, IN 47906
USA
E-mail: chengg@purdue.edu

# SUPPLEMENTARY MATERIAL FOR: A PARTIALLY LINEAR FRAMEWORK FOR MASSIVE HETEROGENOUS DATA

By Tianqi Zhao[§], Guang Cheng[††] and Han Liu[§]

In this supplemental material, we provide the detailed proofs of results presented in the main text. Appendix A contains theoretical justification of RKHS extension to the partially linear function space, as discussed in Section 2. Appendix B, C and D present the proofs of results in Section 3, 4 and 5 respectively. Appendix E contains the proofs of lemmas used in Section 7. Appendix F proves an exponential inequality for empirical processes in a Banach space, and Appendix G provides the proofs of auxiliary lemmas which are used in Appendix E.

## APPENDIX A: RKHS EXTENSION TO PARTIALLY LINEAR FUNCTION SPACE

In this section, we provide detailed theoretical justifications for the RKHS extension to the partially linear space. We first study the properties of the inner product $\langle\cdot,\cdot\rangle_{\mathcal{C}}$ and its induced kernel $\widetilde{K}$, and then prove Proposition 2.3. In the end we provide technical lemmas for the properties of $R_u$ and $P_\lambda$ .

**A.1. A Collection of Lemmas.** The following lemmas are direct consequences of defining the new inner product $\langle\cdot,\cdot\rangle_{\mathcal{C}}$. Lemma A.1 proves the existence of kernel $\widetilde{K}$ under the new inner product $\langle\cdot,\cdot\rangle_{\mathcal{C}}$, and derives its closed form. Lemma A.2 justifies the existence of the linear operator $W_\lambda$ and derives its closed form. Lemma A.3 studies the limit of the $\boldsymbol{B}-\boldsymbol{A}$, where recall $\boldsymbol{A}$ is the Reisz representer of $\boldsymbol{B}=\mathbb{E}[\boldsymbol{X}\mid Z]$ under the inner product $\langle\cdot,\cdot\rangle_{\mathcal{C}}$.

**Lemma A.1.** The linear evaluation functional $E_z$ of $\mathcal{H}$ under the inner product $\langle\cdot,\cdot\rangle_{\mathcal{C}}$ is bounded. Hence $\langle\cdot,\cdot\rangle_{\mathcal{C}}$ induces a new kernel $\widetilde{K}(z,z)$ with the form

$$\widetilde{K}_z(\cdot)=\sum_{\ell=1}^{\infty}\frac{\phi_\ell(z)}{1+\lambda/\mu_\ell}\phi_\ell(\cdot). \tag{A.1}$$

Moreover, under Assumption 3.2, we have that $\|\widetilde{K}_z\|_{\mathcal{C}}\le c_\phi d(\lambda)^{1/2}$, where $c_\phi$ is the constant specified in Assumption 3.2. This implies that $\|f\|_{\sup}\le c_\phi d(\lambda)^{1/2}\|f\|_{\mathcal{C}}$ for all $f\in\mathcal{H}$.

PROOF. (i) Boundedness of $E_z$: We have for any $f \in \mathcal{H}$,

$$|E_z f| = |f(z)| = |\langle f, K_z\rangle_{\mathcal{H}}| \le \|K_z\|_{\mathcal{H}}\|f\|_{\mathcal{H}} \le \lambda^{-1/2} c_k \|f\|_{\mathcal{C}},$$

where the last inequality follows from the relationship $\lambda\|f\|_{\mathcal{H}}^2 \le \|f\|_{\mathcal{C}}^2$ implied by the definition of $\langle\cdot,\cdot\rangle_{\mathcal{C}}$. It follows that $E_z$ is bounded.

(ii) Existence and exact form of $\widetilde{K}$: By Definition 2.1, we have that $\langle\cdot,\cdot\rangle_{\mathcal{C}}$ induces a new kernel $\widetilde{K}(z,z)$. As $\widetilde{K}_z \in \mathcal{H}$ for all $z \in \mathcal{Z}$, by Fourier expansion, $\widetilde{K}_z = \sum_{\ell=1}^{\infty} \kappa_\ell \phi_\ell$. Then we have

$$\begin{aligned}\kappa_\ell = \langle \widetilde{K}_z, \phi_\ell\rangle_{L_2(\mathbb{P}_Z)} &= \langle \widetilde{K}_z, \phi_\ell\rangle_{\mathcal{C}} - \lambda\langle \widetilde{K}_z, \phi_\ell\rangle_{\mathcal{H}} \\ &= \phi_\ell(z) - \lambda\kappa_\ell/\mu_\ell.\end{aligned}$$

Solving for $\kappa_\ell$, we have $\kappa_\ell = \phi_\ell(z)/(1+\lambda/\mu_\ell)$. Hence we get the formula for $\widetilde{K}_z(\cdot)$ in (A.1).

(iii) Uniform bound of $\widetilde{K}_z$: By (A.1) and reproducing property, we have that

$$\|\widetilde{K}_z\|_{\mathcal{C}}^2 = \langle \widetilde{K}_z, \widetilde{K}_z\rangle_{\mathcal{C}} = \widetilde{K}(z,z) = \sum_{\ell=1}^{\infty} \frac{\phi_\ell^2(z)}{1+\lambda/\mu_\ell} \le c_\phi^2 d(\lambda).$$

as desired. Hence for all $z \in \mathcal{Z}$,

$$|f(z)| \le \|f\|_{\mathcal{C}}\|\widetilde{K}_z\|_{\mathcal{C}} \le c_\phi d(\lambda)^{1/2}\|f\|_{\mathcal{C}},$$

by Cauchy-Schwarz. This implies that $\|f\|_{\sup} \le c_\phi d(\lambda)^{1/2}\|f\|_{\mathcal{C}}$ for all $f \in \mathcal{H}$. □

**Lemma A.2.** There exists a bounded linear operator $W_\lambda : \mathcal{H} \to \mathcal{H}$ such that for any $f, \widetilde{f} \in \mathcal{H}$, we have

$$\langle W_\lambda f, \widetilde{f}\rangle_{\mathcal{C}} = \lambda\langle f, \widetilde{f}\rangle_{\mathcal{H}}. \tag{A.2}$$

Moreover, we have for all eigenfunctions $\phi_\ell$, $\ell = 1, 2, \ldots$

$$W_\lambda \phi_\ell(\cdot) = \frac{\lambda}{\lambda+\mu_\ell}\phi_\ell(\cdot). \tag{A.3}$$

PROOF. The proof for the existence of $W_\lambda$ uses Riesz representation theorem. Define the bilinear form $V(f,\widetilde{f}) := \lambda\langle f, \widetilde{f}\rangle_{\mathcal{H}}$, for any $f, \widetilde{f} \in \mathcal{H}$. For any fixed $f$, this defines a functional $V_f(\cdot) = V(f,\cdot)$. It is easy to verify that $V_f$ is linear. Moreover, $V_f$ is bounded under the inner product $\langle\cdot,\cdot\rangle_{\mathcal{C}}$, as

$$|V_f(\widetilde{f})| = |\lambda\langle f, \widetilde{f}\rangle_{\mathcal{H}}| \le \lambda\|f\|_{\mathcal{H}}\|\widetilde{f}\|_{\mathcal{H}} \le \lambda^{1/2}\|f\|_{\mathcal{H}}\|\widetilde{f}\|_{\mathcal{C}},$$

for all $\widetilde{f} \in \mathcal{H}$. Hence by Riesz representation theorem, there exists an unique element $f_1 \in \mathcal{H}$ such that $V_f(\widetilde{f}) = \langle f_1, \widetilde{f}\rangle_{\mathcal{C}}$ for all $\widetilde{f} \in \mathcal{H}$. We let $W_\lambda f = f_1$, and it follows that $\langle W_\lambda f, \widetilde{f}\rangle_{\mathcal{C}} = \lambda\langle f, \widetilde{f}\rangle_{\mathcal{H}}$ for all $f, \widetilde{f} \in \mathcal{H}$.

We next prove linearity of $W_\lambda$. By definition, for any $f_1, f_2, \widetilde{f} \in \mathcal{H}$ and $a, b \in \mathbb{R}$, we have

$$\begin{aligned}\langle W_\lambda(af_1 + bf_2), \widetilde{f}\rangle_{\mathcal{C}} &= \lambda\langle af_1 + bf_2, \widetilde{f}\rangle_{\mathcal{H}} \\ &= a\lambda\langle f_1, \widetilde{f}\rangle_{\mathcal{H}} + b\lambda\langle f_2, \widetilde{f}\rangle_{\mathcal{H}} \\ &= \langle aW_\lambda f_1 + bW_\lambda f_2, \widetilde{f}\rangle_{\mathcal{C}},\end{aligned}$$

which implies $W_\lambda(af_1 + bf_2) = aW_\lambda f_1 + bW_\lambda f_2$.

Furthermore,$W_\lambda$ is a bounded operator under $\|\cdot\|_{\mathcal{C}}$, as for any $f \in \mathcal{H}$

$$\begin{aligned}\|W_\lambda f\|_{\mathcal{C}} &= \sup_{\|\widetilde{f}\|_{\mathcal{C}}\le 1} \langle W_\lambda f, \widetilde{f}\rangle_{\mathcal{C}} = \sup_{\|\widetilde{f}\|_{\mathcal{C}}\le 1} \lambda\langle f, \widetilde{f}\rangle_{\mathcal{H}} \\ &\le \lambda\|f\|_{\mathcal{H}} \sup_{\|\widetilde{f}\|_{\mathcal{C}}\le 1} \|\widetilde{f}\|_{\mathcal{H}} \le \|f\|_{\mathcal{C}},\end{aligned}$$

where the last inequality follows from the fact that $\lambda^{1/2}\|f\|_{\mathcal{H}} \le \|f\|_{\mathcal{C}}$ implied by the definition of $\langle\cdot,\cdot\rangle_{\mathcal{C}}$. Hence, the operator norm of $W_\lambda$ is bounded by 1.

To prove the second half of the lemma, we have that $\langle W_\lambda f, \widetilde{f}\rangle_{\mathcal{C}} = \lambda\langle f, \widetilde{f}\rangle_{\mathcal{H}}$. Also, by the definition of $\langle\cdot,\cdot\rangle_{\mathcal{C}}$ we have $\langle W_\lambda f, \widetilde{f}\rangle_{\mathcal{C}} = \langle W_\lambda f, \widetilde{f}\rangle_{L_2(P_Z)} + \lambda\langle W_\lambda f, \widetilde{f}\rangle_{\mathcal{H}}$. It follows from the two equations that

$$\langle W_\lambda f, \widetilde{f}\rangle_{L_2(\mathbb{P}_Z)} = \lambda\langle (id - W_\lambda)f, \widetilde{f}\rangle_{\mathcal{H}}, \tag{A.4}$$

for any $f, \widetilde{f} \in \mathcal{H}$. $W_\lambda\phi_\ell$ has a Fourier expansion: $W_\lambda\phi_\ell = \sum_{k=1}^\infty w_k\phi_k$. Letting $f = \widetilde{f} = \phi_\ell$ in (A.4) yields $w_\ell = \lambda/(\lambda + \mu_\ell)$, and letting $f = \phi_\ell$ and $\widetilde{f} = \phi_r$ in (A.4) yields $w_r = 0$ for $r \ne \ell$. Hence the conclusion follows. □

**Lemma A.3.** Under Assumptions 3.1 and 3.2, the following equations hold:

$$\lim_{\lambda\to 0} \mathbb{E}\Big[\boldsymbol{X}\big(\boldsymbol{B}(Z) - \boldsymbol{A}(Z)\big)^T\Big] = \mathbf{0}, \tag{A.5}$$

$$\lim_{\lambda\to 0} \mathbb{E}\Big[\boldsymbol{B}(Z)\big(\boldsymbol{B}(Z) - \boldsymbol{A}(Z)\big)^T\Big] = \mathbf{0}, \tag{A.6}$$

$$\lim_{\lambda\to 0} \mathbb{E}\Big[\big(\boldsymbol{B}(Z) - \boldsymbol{A}(Z)(\boldsymbol{B}(Z) - \boldsymbol{A}(Z))^T\big)\Big] = \mathbf{0}. \tag{A.7}$$

The lemma shows that the difference $\boldsymbol{B} - \boldsymbol{A}$ goes to zero as $\lambda \to 0$. Intuititively, as $\lambda \to 0$, the inner product $\langle\cdot,\cdot\rangle_{\mathcal{C}}$ converges to $\langle\cdot,\cdot\rangle_{L_2(\mathbb{P}_Z)}$ by its definition, hence the the representer $\boldsymbol{A}$ of $\boldsymbol{B}$ converges to $\boldsymbol{B}$ itself. The following is the formal proof of this lemma.

PROOF. By reproducing property of $\widetilde{K}_z$, the definition of $A_k$ and (A.1), we have

$$A_k(z) = \langle A_k, \widetilde{K}_z\rangle_{\mathcal{C}} = \langle B_k, \widetilde{K}_z\rangle_{L_2(\mathbb{P}_Z)} = \sum_{i=1}^\infty \frac{\langle B_k, \phi_i\rangle_{L_2(\mathbb{P}_Z)}}{1 + \lambda/\mu_i}\phi_i(z). \tag{A.8}$$

By Fourier expansion of $B_k$, it follows from the above equation that

$$B_k(z) - A_k(z) = \sum_{i=1}^{\infty} \frac{\langle B_k, \phi_i \rangle_{L_2(\mathbb{P}_Z)} \lambda/\mu_i}{1 + \lambda/\mu_i} \phi_i(z), \tag{A.9}$$

Therefore, for any $j, k \in \{1, \dots, p\}$, we have

$$\begin{aligned} \mathbb{E}\big[X_j(B_k(Z) - A_k(Z))\big] &= \sum_{i=1}^{\infty} \frac{\lambda/\mu_i}{1 + \lambda/\mu_i} \langle B_k, \phi_i \rangle_{L_2(\mathbb{P}_Z)} \mathbb{E}[X_j \phi_i(Z)] \\ &= \sum_{i=1}^{\infty} \frac{\lambda/\mu_i}{1 + \lambda/\mu_i} \langle B_k, \phi_i \rangle_{L_2(\mathbb{P}_Z)} \langle B_j, \phi_i \rangle_{L_2(\mathbb{P}_Z)} \leq \infty, \end{aligned}$$

where the second equality is by the fact that $\mathbb{E}[X_j \phi_i(Z)] = \mathbb{E}[B_j \phi_i(Z)] = \langle B_j, \phi_i \rangle_{L_2(\mathbb{P}_Z)}$, and the inequality is by Cauchy-Schwarz and the fact that $\sum_{i=1}^{\infty} \langle B_j, \phi_i \rangle^2_{L_2(\mathbb{P}_Z)} \leq \infty$ as $B_j \in L_2(\mathbb{P}_Z)$. By dominated convergence theorem, we have (A.5) holds. Moreover,

$$\begin{aligned} \mathbb{E}\big[B_j(Z)(B_k(Z) - A_k(Z))\big] &= \langle B_j, B_k - A_k \rangle_{L_2(\mathbb{P}_Z)} \\ &= \sum_{i=1}^{\infty} \frac{\lambda/\mu_i}{1 + \lambda/\mu_i} \langle B_j, \phi_i \rangle_{L_2(\mathbb{P}_Z)} \langle B_k, \phi_i \rangle_{L_2(\mathbb{P}_Z)} \\ &\to 0, \end{aligned} \tag{A.10}$$

where the second equality is by (A.9) and the limit is by dominated convergence theorem. Hence (A.6) holds. Applying similar arguments, we can show that (A.7) also holds. □

**A.2. Proof of Proposition 2.3.** With the theoretical foundation established in the previous section, we are now ready to construct $R_u$ and $P_\lambda$, whose exact forms are presented in Proposition 2.3.

Proof. The proof follows similarly as Proposition 2.1 in Cheng and Shang (2015). We first want to construct $R_u \in \mathcal{A}$ such that it possess the following reproducing property:

$$\langle R_u, m \rangle_{\mathcal{A}} = \boldsymbol{\beta}^T \boldsymbol{x} + f(z), \tag{A.11}$$

for any $u = (\boldsymbol{x}, z)$ and $m = (\boldsymbol{\beta}, f) \in \mathcal{A}$. As $R_u \in \mathcal{A}$, it has two components:

$R_u = (L_u, N_u)$. Hence the L.H.S. of (A.11) can be written as

$$\begin{aligned}\langle R_u, m\rangle_{\mathcal{A}} &= \mathbb{E}\big[(\boldsymbol{X}^T L_u + N_u(Z))(\boldsymbol{X}^T\boldsymbol{\beta} + f(Z))\big] + \lambda\langle N_u, f\rangle_{\mathcal{H}} \\ &= \boldsymbol{\beta}^T\mathbb{E}[\boldsymbol{X}\boldsymbol{X}^T]L_u + \boldsymbol{\beta}^T\mathbb{E}[\boldsymbol{X}N_u(Z)] + L_u^T\mathbb{E}[\boldsymbol{X}f(Z)] \\ &\qquad + \mathbb{E}[N_u(Z)f(Z)] + \lambda\langle N_u, f\rangle_{\mathcal{H}} \\ &= \boldsymbol{\beta}^T\big(\mathbb{E}[\boldsymbol{X}\boldsymbol{X}^T]L_u + \mathbb{E}[\boldsymbol{B}(Z)N_u(Z)]\big) + L_u^T\mathbb{E}[\boldsymbol{B}(Z)f(Z)] + \langle N_u, f\rangle_{\mathcal{C}} \\ &= \boldsymbol{\beta}^T\big(\mathbb{E}[\boldsymbol{X}\boldsymbol{X}^T]L_u + \mathbb{E}[\boldsymbol{B}(Z)N_u(Z)]\big) + \langle \boldsymbol{A}^T L_u + N_u, f\rangle_{\mathcal{C}} \end{aligned} \tag{A.12}$$

where in the second last inequality we used the definition of $\langle\cdot,\cdot\rangle_{\mathcal{C}}$ and in the last equality we used $\mathbb{E}[\boldsymbol{B}(Z)f(Z)] = \langle \boldsymbol{B}, f\rangle_{L_2(\mathbb{P}_Z)} = \langle \boldsymbol{A}, f\rangle_{\mathcal{C}}$. On the other hand, the R.H.S of (A.11) is

$$\boldsymbol{\beta}^T\boldsymbol{x} + f(z) = \boldsymbol{\beta}^T\boldsymbol{x} + \langle \widetilde{K}_z, f\rangle. \tag{A.13}$$

Comparing (A.12) and (A.13), we have the following set of equations:

$$\begin{aligned}\boldsymbol{x} &= \mathbb{E}[\boldsymbol{X}\boldsymbol{X}^T]L_u + \mathbb{E}[\boldsymbol{B}(Z)N_u(Z)] \\ \widetilde{K}_z &= \boldsymbol{A}^T L_u + N_u.\end{aligned}$$

From the second equation we get $N_u = \widetilde{K}_z - \boldsymbol{A}^T L_u$. Substitute it into the first equation, we get

$$\begin{aligned}\boldsymbol{x} &= \mathbb{E}[\boldsymbol{X}\boldsymbol{X}^T]L_u + \mathbb{E}[\boldsymbol{B}(Z)(\widetilde{K}_z(Z) - \boldsymbol{A}(Z)^T L_u)] \\ &= \big(\boldsymbol{\Omega} + \mathbb{E}[\boldsymbol{B}(Z)\boldsymbol{B}^T(Z)]\big)L_u + \mathbb{E}[\boldsymbol{B}(Z)\widetilde{K}_z(Z)] - \mathbb{E}[\boldsymbol{B}(Z)\boldsymbol{A}(Z)^T L_u] \\ &= (\boldsymbol{\Omega} + \boldsymbol{\Sigma}_\lambda)L_u + \langle \boldsymbol{B}, \widetilde{K}_z\rangle_{L_2(\mathbb{P}_Z)} \\ &= (\boldsymbol{\Omega} + \boldsymbol{\Sigma}_\lambda)L_u + \langle \boldsymbol{A}, \widetilde{K}_z\rangle_{\mathcal{C}} = (\boldsymbol{\Omega} + \boldsymbol{\Sigma}_\lambda)L_u + \boldsymbol{A}(z),\end{aligned}$$

where in the second inequality we used the fact that $\boldsymbol{B}(Z)$ and $\boldsymbol{X} - \boldsymbol{B}(Z)$ are orthogonal. Therefore it follows that

$$L_u = (\boldsymbol{\Omega} + \boldsymbol{\Sigma}_\lambda)^{-1}(\boldsymbol{x} - \boldsymbol{A}(z)).$$

This finishes the proof for the construction of $R_u$.

We next construct $P_\lambda$ such that

$$\langle P_\lambda m, \widetilde{m}\rangle_{\mathcal{A}} = \lambda\langle f, \widetilde{f}\rangle_{\mathcal{H}}, \tag{A.14}$$

for any $m = (\boldsymbol{\beta}, f), \widetilde{m} = (\widetilde{\boldsymbol{\beta}}, \widetilde{f}) \in \mathcal{A}$. As $P_\lambda m \in \mathcal{A}$, it has two components: $P_\lambda m = (L_\lambda f, N_\lambda f)$. Similar to the derivation of (A.12), the L.H.S. of (A.14) can be written as

$$\begin{aligned}\langle P_\lambda m, \widetilde{m}\rangle_{\mathcal{A}} &= \mathbb{E}\big[(\boldsymbol{X}^T L_\lambda f + N_\lambda f(Z))(\boldsymbol{X}^T\widetilde{\boldsymbol{\beta}} + \widetilde{f})\big] + \lambda\langle N_\lambda f, \widetilde{f}\rangle_{\mathcal{H}} \\ &= \widetilde{\boldsymbol{\beta}}^T\big(\mathbb{E}[\boldsymbol{X}\boldsymbol{X}^T]L_\lambda f + \mathbb{E}[\boldsymbol{B}(Z)N_\lambda f(Z)]\big) + \langle \boldsymbol{A}^T L_\lambda f + N_\lambda f, \widetilde{f}\rangle_{\mathcal{C}}\end{aligned} \tag{A.15}$$

The R.H.S. of (A.14) is

$$\lambda\langle f, \widetilde{f}\rangle_{\mathcal{H}} = \langle W_\lambda f, \widetilde{f}\rangle_{\mathcal{C}}. \tag{A.16}$$

Comparing (A.15) and (A.16), we obtain the following set of equations:

$$\mathbf{0} = \mathbb{E}[\boldsymbol{X}\boldsymbol{X}^T]L_\lambda f + \mathbb{E}[\boldsymbol{B}(Z)N_\lambda f(Z)]$$
$$W_\lambda f = \boldsymbol{A}^T L_\lambda f + N_\lambda f.$$

Solving the above two equations, we get

$$L_\lambda f = -(\boldsymbol{\Omega} + \boldsymbol{\Sigma}_\lambda)^{-1}\langle \boldsymbol{B}, W_\lambda f\rangle_{L_2(\mathbb{P}_Z)} \text{ and } N_\lambda f = W_\lambda f - \boldsymbol{A}^T L_\lambda f,$$

as desired. □

**A.3. Properties of $\boldsymbol{R_u}$ and $\boldsymbol{P_\lambda}$.** We first present a lemma that bounds the $\mathcal{A}$-norm of $R_u$.

**Lemma A.4.** Under Assumptions 3.1 and 3.2, there exists a constant $c_r > 0$ independent to $u$ such that $\|R_u\|_{\mathcal{A}} \le c_r d(\lambda)^{1/2}$. It follows that $\|m\|_{\sup} \le c_r d(\lambda)^{1/2}\|m\|_{\mathcal{A}}$ for all $m$.

Proof. By Proposition 2.3, we have that for $u = (\boldsymbol{x}, z)$,

$$\begin{aligned}\langle R_u, R_u\rangle_{\mathcal{A}} &= \boldsymbol{x}^T L_u + N_u(z)\\ &= \widetilde{K}_z(z) + (\boldsymbol{x} - \boldsymbol{A}(z))^T L_u\\ &= \widetilde{K}_z(z) + (\boldsymbol{x} - \boldsymbol{A}(z))^T(\boldsymbol{\Omega} + \boldsymbol{\Sigma}_\lambda)^{-1}(\boldsymbol{x} - \boldsymbol{A}(z)).\end{aligned} \tag{A.17}$$

From Lemma A.1, we have $\widetilde{K}_z(z) = \|\widetilde{K}\|_{\mathcal{C}}^2 \le c_\phi^2 d(\lambda)$. For the second term in (A.17), we first show that $\boldsymbol{\Sigma}_\lambda$ is positive definite. Recall the definition of $\boldsymbol{\Sigma}_\lambda = \mathbb{E}\big[\boldsymbol{B}(Z)(\boldsymbol{B}(Z) - \boldsymbol{A}(Z))\big]$. (A.10) shows that

$$[\boldsymbol{\Sigma}_\lambda]_{jk} = \mathbb{E}\big[B_j(Z)(B_k(Z) - A_k(Z))\big] = \sum_{\ell=1}^{\infty} \frac{\lambda/\mu_\ell}{1 + \lambda/\mu_\ell}\langle B_j, \phi_\ell\rangle_{L_2}\langle B_k, \phi_\ell\rangle_{L_2},$$

which implies that $\boldsymbol{\Sigma}_\lambda$ is positive definite. Indeed, for any $\boldsymbol{x} \in \mathbb{R}^p$ and $\boldsymbol{x} \neq 0$,

$$\boldsymbol{x}^T\boldsymbol{\Sigma}_\lambda\boldsymbol{x} = \sum_{j=1}^{p}\sum_{k=1}^{p} x_j x_k[\boldsymbol{\Sigma}_\lambda]_{jk} = \sum_{\ell=1}^{\infty} \frac{\lambda/\mu_\ell}{1 + \lambda/\mu_\ell}\Big(\sum_{j=1}^{p} x_j\langle B_j, \phi_\ell\rangle_{L_2}\Big)^2 > 0.$$

Therefore, it follows that the second term in (A.17) is bounded by

$$\begin{aligned}(\boldsymbol{x} - \boldsymbol{A}(z))^T(\boldsymbol{\Omega} + \boldsymbol{\Sigma}_\lambda)^{-1}(\boldsymbol{x} - \boldsymbol{A}(z)) &\le \|\boldsymbol{\Omega} + \boldsymbol{\Sigma}_\lambda\|^{-1}\|\boldsymbol{x} - \boldsymbol{A}(z)\|_2^2\\ &\le \tau_{\min}(\boldsymbol{\Omega})^{-1}\|\boldsymbol{x} - \boldsymbol{A}(z)\|_2^2.\end{aligned} \tag{A.18}$$

where recall $\tau_{\min}(\boldsymbol{\Omega})$ is the minimum eigenvalue of $\boldsymbol{\Omega}$. As $\boldsymbol{x}$ is uniformly bounded, we are left to bound $\boldsymbol{A}(z)$. By (A.8), we have

$$A_k(z) = \sum_{\ell=1}^{\infty} \frac{\langle B_k, \phi_\ell \rangle_{L_2(\mathbb{P}_Z)}}{1 + \lambda/\mu_\ell} \phi_\ell(z), \tag{A.19}$$

Hence by the assumption that $B_j \in L_2(\mathbb{P}_Z)$ for $j = 1, \ldots, p$ and uniform boundedness of $\phi_\ell$, we have

$$A_k^2(z) \leq c_\phi \sum_{\ell=1}^{\infty} \langle B_k, \phi_\ell \rangle_{L_2(\mathbb{P}_Z)}^2 \sum_{\ell=1}^{\infty} (1 + \lambda/\mu_\ell)^{-2} \leq c_\phi \|B_k\|_{L_2(\mathbb{P}_Z)}^2 d(\lambda). \tag{A.20}$$

Therefore, by (A.17) and (A.18), we have $\|R_u\|_{\mathcal{A}} \leq c_r d(\lambda)^{1/2}$, where $c_r$ is determined by $p, c_\phi, c_x, \tau_{\min}$ and $\|B_k\|_{L_2(\mathbb{P}_Z)}$. Therefore, for any $u = (\boldsymbol{x}, z)$

$$|m(u)| = |\langle m, R_u \rangle_{\mathcal{A}}| \leq \|m\|_{\mathcal{A}} \|R_u\|_{\mathcal{A}} \leq c_r d(\lambda)^{1/2} \|m\|_{\mathcal{A}},$$

which implies that $\|m\|_{\sup} \leq c_r d(\lambda)^{1/2} \|m\|_{\mathcal{A}}$. □

Based on the above lemma, if we have an extra condition that $B_k(z)$ is smooth, then we can bound the parametric and nonparametric components of $R_u$ and $P_\lambda m_0$ more precisely.

**Lemma A.5.** Suppose Assumptions 3.1 - 3.3 hold. Then we have

(i) $\|L_u\|_2^2 \leq C_1'$ where $C_1' = 2\tau_{\min}^{-2}\Big(c_x^2 p + c_\phi^2 \operatorname{Tr}(K) \sum_{k=1}^p \|B_k\|_{\mathcal{H}}^2\Big)$, and $\|N_u\|_{\mathcal{C}}^2 \leq C_1 d(\lambda)$ where $C_1 = 2c_\phi^2\Big(1 + C_1' \sum_{k=1}^p \|B_k\|_{L_2(\mathbb{P}_Z)}^2\Big)$;

(ii) Moreover, $\|L_\lambda f_0\|_2^2 \leq C_2' \lambda^2$ and $\|N_\lambda f_0\|_{\mathcal{C}}^2 \leq 2\|f_0\|_{\mathcal{H}}^2 \lambda + C_2 \lambda^2$, where $C_2' = \tau_{\min}^{-2} \|f_0\|_{\mathcal{H}}^2 \sum_{k=1}^p \|B_k\|_{\mathcal{H}}^2$ and $C_2 = 2C_2' \sum_{k=1}^p \|B_k\|_{L_2(\mathbb{P}_Z)}^2$.

Proof. (i) By Proposition 2.3,

$$L_u = (\boldsymbol{\Omega} + \boldsymbol{\Sigma}_\lambda)^{-1}(\boldsymbol{x} - \boldsymbol{A}(z)) \text{ and } N_u = \widetilde{K}_z - \boldsymbol{A}^T L_u,$$

By the first equation we have $\|L_u\|_2^2 \leq \|(\boldsymbol{\Omega} + \boldsymbol{\Sigma}_\lambda)^{-1}\|_2^2 \|\boldsymbol{x} - \boldsymbol{A}(z)\|_2^2$. Recall

$$A_k(z) = \sum_{\ell=1}^{\infty} \frac{\langle B_k, \phi_\ell \rangle_{L_2}}{1 + \lambda/\mu_\ell} \phi_\ell(z).$$

Hence by Assumption 3.3 it follows that

$$A_k(z)^2 \leq \sum_{\ell=1}^{\infty} \frac{\langle B_k, \phi_\ell \rangle_{L_2}^2}{\mu_\ell} \sum_{\ell=1}^{\infty} \mu_\ell \frac{\phi_\ell(z)^2}{(1 + \lambda/\mu_\ell)^2} \leq c_\phi^2 \|B_k\|_{\mathcal{H}}^2 \operatorname{Tr}(K),$$

where the first inequality is by Cauchy-Schwarz and the second is by As-sumption 3.2 that $\phi_\ell$ are uniformly bounded. Hence for all $z \in \mathcal{Z}$,

$$\|\boldsymbol{A}(z)\|_2^2 \le c_\phi^2 \operatorname{Tr}(K) \sum_{k=1}^{p} \|B_k\|_{\mathcal{H}}^2. \tag{A.21}$$

Also, we showed in the proof of Lemma A.4 that $\|(\boldsymbol{\Omega} + \boldsymbol{\Sigma}_\lambda)^{-1}\|$ is bounded by $\tau_{\min}^{-1}$. Finally, by the boundedness of the support of $\mathcal{X}$, we have

$$\|L_U\|_2^2 \le 2\tau_{\min}^{-2}\Big(c_x^2 p + c_\phi^2 \operatorname{Tr}(K) \sum_{k=1}^{p} \|B_k\|_{\mathcal{H}}^2\Big) = C_1'.$$

To control $N_u$, we have

$$\|N_u\|_{\mathcal{C}}^2 \le 2\big(\|\widetilde{K}_z\|_{\mathcal{C}}^2 + \|L_u^T \boldsymbol{A}\|_{\mathcal{C}}^2\big) \tag{A.22}$$

For the first term in (A.22), by Lemma A.1, we have $\|\widetilde{K}_z\|_{\mathcal{C}}^2 \le c_\phi^2 d(\lambda)$. For the second term, by (A.19) we have

$$\begin{aligned} \|A_k\|_{\mathcal{C}}^2 &= \langle A_k, A_k\rangle_{\mathcal{C}} = \langle B_k, A_k\rangle_{L_2(\mathbb{P}_Z)} = \sum_{\ell=1}^{\infty} \frac{\langle B_k, \phi_\ell\rangle_{L_2(\mathbb{P}_Z)}^2}{1+\lambda/\mu_\ell} \\ &\le c_\phi^2 \|B_k\|_{L_2(\mathbb{P}_Z)}^2 \sum_{\ell=1}^{\infty} \frac{1}{1+\lambda/\mu_\ell} = c_\phi^2 \|B_k\|_{L_2(\mathbb{P}_Z)}^2 d(\lambda). \end{aligned} \tag{A.23}$$

where the inequality is by Cauchy-Schwartz and uniform boundedness of $\phi_\ell$. Hence it follows that

$$\begin{aligned} \|L_u^T \boldsymbol{A}\|_{\mathcal{C}} &= \|\sum_{k=1}^{p} (L_u)_k A_k\|_{\mathcal{C}} \le \sum_{k=1}^{p} |(L_u)_k| \|A_k\|_{\mathcal{C}} \\ &\le \|L_u\|_2 \Big(\sum_{k=1}^{p} \|A_k\|_{\mathcal{C}}^2\Big)^{1/2} \le (C_1')^{1/2} c_\phi d(\lambda)^{1/2} \Big(\sum_{k=1}^{p} \|B_k\|_{L_2(\mathbb{P}_Z)}^2\Big)^{1/2} \end{aligned} \tag{A.24}$$

Therefore, by (A.22) we obtain

$$\|N_u\|_{\mathcal{C}}^2 \le 2c_\phi^2\Big(1 + C_1' \sum_{k=1}^{p} \|B_k\|_{L_2(\mathbb{P}_Z)}^2\Big) d(\lambda) = C_1 d(\lambda)$$

(ii) Let $\{\theta_\ell\}_{\ell=1}^\infty$ be the Fourier coefficient of $f_0$ under the basis $\{\phi_\ell\}_{\ell=1}^\infty$ (given the $L_2(\mathbb{P}_Z)$-inner product $\langle\cdot,\cdot\rangle_{L_2(\mathbb{P}_Z)}$). We have $\|f_0\|_{L_2(\mathbb{P}_Z)}^2 = \sum_{i=1}^\infty \theta_i^2$ and $\|f_0\|_{\mathcal{H}}^2 = \sum_{i=1}^\infty \theta_i^2/\mu_i$, where the second equality is by Mercer's theorem. By Fourier expansion of $f_0$ and (A.3), we have

$$W_\lambda f_0 = \sum_{\ell=1}^{\infty} \frac{\lambda \theta_\ell}{\mu_\ell + \lambda} \phi_\ell \tag{A.25}$$

By Proposition 2.3, we have $L_\lambda f_0 = -(\boldsymbol{\Omega} + \boldsymbol{\Sigma}_\lambda)^{-1} \langle W_\lambda f_0, \boldsymbol{B}\rangle_{L_2(\mathbb{P}_Z)}$ and

$N_\lambda f_0 = (L_\lambda f_0)^T \boldsymbol{A} + W_\lambda f_0$. To control $L_\lambda f_0$, we calculate

$$\begin{aligned}
\langle B_j, W_\lambda f_0 \rangle^2_{L_2(\mathbb{P}_Z)} &= \Big(\sum_{\ell=1}^\infty \langle B_k, \phi_\ell\rangle_{L_2} \frac{\lambda\theta_\ell}{\mu_\ell+\lambda}\Big)^2 \\
&\le \lambda^2 \Big(\sum_{i=1}^\infty \frac{\langle B_k, \phi_\ell\rangle^2_{L_2}}{\mu_\ell+\lambda}\Big)\Big(\sum_{\ell=1}^\infty \frac{\theta_\ell^2}{\mu_\ell+\lambda}\Big) \\
&\le \lambda^2 \|B_k\|^2_{\mathcal{H}} \|f_0\|^2_{\mathcal{H}}.
\end{aligned}$$

Hence by the positive definiteness of $\boldsymbol{\Sigma}_\lambda$, we have

$$\begin{aligned}
\|L_{f_0}\|_2^2 &\le \|(\boldsymbol{\Omega}+\boldsymbol{\Sigma}_\lambda)^{-1}\| \|\langle W_\lambda f_0, \boldsymbol{B}\rangle_{L_2(\mathbb{P}_Z)}\|_2^2 \\
&\le \tau_{\min}^{-2} \|f_0\|^2_{\mathcal{H}} \sum_{k=1}^p \|B_k\|^2_{\mathcal{H}} \lambda^2 = C_2' \lambda^2.
\end{aligned}$$

For $N_\lambda f_0$, we have

$$\|N_\lambda f_0\|^2_{\mathcal{C}} = 2\|(L_\lambda f_0)^T \boldsymbol{A}\|^2_{\mathcal{C}} + 2\|W_\lambda f_0\|^2_{\mathcal{C}} \tag{A.26}$$

For the first term (A.26), by the fact that $\sum_{\ell=1}^\infty \langle B_k, \phi_\ell\rangle^2_{L_2(\mathbb{P}_Z)} = \|B_k\|^2_{L_2(\mathbb{P}_Z)}$, we first get an inequality for $\|A_k\|^2_{\mathcal{C}}$ that is different than (A.23):

$$\|A_k\|^2_{\mathcal{C}} = \langle A_k, A_k\rangle_{\mathcal{C}} = \langle B_k, A_k\rangle_{L_2(\mathbb{P}_Z)} = \sum_{\ell=1}^\infty \frac{\langle B_k, \phi_\ell\rangle^2_{L_2(\mathbb{P}_Z)}}{1+\lambda/\mu_\ell} \le \|B_k\|^2_{L_2(\mathbb{P}_Z)}.$$

Therefore, following the same derivation as (A.24), we have

$$\|(L_\lambda f_0)^T \boldsymbol{A}\|^2_{\mathcal{C}} \le \|L_\lambda f_0\|_2^2 \sum_{k=1}^p \|A_k\|^2_{\mathcal{C}} \le C_2' \lambda^2 \sum_{k=1}^p \|B_k\|^2_{L_2(\mathbb{P}_Z)}.$$

For the second term in (A.26) we have

$$\begin{aligned}
\|W_\lambda f_0\|_{\mathcal{C}} &= \sup_{\|f\|_{\mathcal{C}}=1} |\langle W_\lambda f_0, f\rangle_{\mathcal{C}}| = \sup_{\|f\|_{\mathcal{C}}=1} \lambda |\langle f_0, f\rangle_{\mathcal{H}}| \\
&\le \sup_{\|f\|_{\mathcal{C}}=1} \sqrt{\lambda\|f_0\|^2_{\mathcal{H}}}\sqrt{\lambda\|f\|^2_{\mathcal{H}}} \le \lambda^{1/2}\|f_0\|_{\mathcal{H}},
\end{aligned}$$

where the last equality follows from the fact that $\lambda\|f\|^2_{\mathcal{H}} \le \|f\|^2_{\mathcal{C}}$. Therefore, by (A.26), we obtain

$$\|N_\lambda f_0\|^2_{\mathcal{C}} \le 2\lambda\|f_0\|^2_{\mathcal{H}} + 2C_2'\lambda^2 \sum_{k=1}^p \|B_k\|^2_{L_2(\mathbb{P}_Z)} = 2\|f_0\|^2_{\mathcal{H}}\lambda + C_2\lambda^2,$$

as desired. □

## APPENDIX B: PROOFS IN SECTION 3

### B.1. Proof of Theorem 3.3.

PROOF. Similar to the arguments in the proof of Theorem 3.2, the bound on $\mathrm{MSE}(\widehat{\boldsymbol{\beta}}^{(j)})$ follows from Eq. (3.5), (3.7) and Lemma A.5. □

### B.2. Proof of Theorem 3.5.

PROOF. We first relabel the data as $U_i^{(j)}$ and $\varepsilon_i^{(j)}$, which denotes the $i$-th data in subsample $j$, for $i = 1, \ldots, n$, $j = 1, \ldots, s$. By (3.5), we obtain

$$\sqrt{n} V_s^{-1} \underset{\sim}{\boldsymbol{u}}^T (\underset{\sim}{\widehat{\boldsymbol{\beta}}} - \underset{\sim}{\boldsymbol{\beta}}_0) = \frac{1}{\sqrt{n}} \sum_{i=1}^n V_s^{-1} \sum_{j=1}^s \varepsilon_i^{(j)} \boldsymbol{u}_j^T L_{U_i^{(j)}} + \sqrt{n} \sum_{j=1}^s V_s^{-1} \boldsymbol{u}_j^T L_\lambda f_0$$

$$+ \sqrt{n} \sum_{j=1}^s V_s^{-1} \boldsymbol{u}_j^T Rem_\beta^{(j)}. \tag{B.1}$$

We first prove that the first term is asymptotically normal. It is a sum of independent random variables. Moreover, $\sum_{j=1}^s \varepsilon_i^{(j)} \boldsymbol{u}_j^T L_{U_i^{(j)}}$ is $\sigma^2 C_1'$-sub-Gaussian, where $C_1'$ is the constant defined in Lemma A.5. Indeed, for any $t \in \mathbb{R}$,

$$\begin{aligned}
\mathbb{E}\left[\exp\left(t \sum_{j=1}^s \varepsilon_i^{(j)} \boldsymbol{u}_j^T L_{U_i^{(j)}}\right)\right] &= \prod_{j=1}^s \mathbb{E}\left[\exp(t\varepsilon_i^{(j)} \boldsymbol{u}_j^T L_{U_i^{(j)}})\right] \\
&= \prod_{j=1}^s \mathbb{E}\left[\mathbb{E}\left[\exp(t\varepsilon_i^{(j)} \boldsymbol{u}_j^T L_{U_i^{(j)}}) \,|\, L_{U_i^{(1)}}, \ldots, L_{U_i^{(s)}}\right]\right] \\
&\le \prod_{j=1}^s \mathbb{E}\left[\exp\left(\frac{t^2 \sigma^2 \|\boldsymbol{u}_j\|_2^2 \|L_{U_i^{(j)}}\|_2^2}{2}\right)\right] \\
&\le \exp\left(\frac{t^2 \sigma^2 C_1'}{2}\right),
\end{aligned}$$

where the first inequality is by Hoeffding, and the last one is by Lemma A.5. Moreover, we have

$$V_s = \sum_{j=1}^s \boldsymbol{u}_j^T \boldsymbol{\Omega}^{-1} \boldsymbol{u}_j \ge \sum_{j=1}^s \|\boldsymbol{u}_j\|_2^2 \tau_{\min}(\boldsymbol{\Omega}^{-1}) = \tau_{\max}(\boldsymbol{\Omega})^{-1},$$

where $\tau_{\min}(\cdot), \tau_{\max}(\cdot)$ denote the minimum and maximum eigenvalues. Therefore, we have $V_s^{-1} \le \tau_{\max}(\boldsymbol{\Omega})$. We now compute the variance of the summands in the first summation. By independence among subpopulations, we have

$$\begin{aligned}\mathrm{Var}\left(V_s^{-1}\sum_{j=1}^s \varepsilon_i^{(j)}\boldsymbol{u}_j^T L_{U_i^{(j)}}\right) &= V_s^{-2}\sum_{j=1}^s \mathrm{Var}\left(\varepsilon_i^{(j)}\boldsymbol{u}_j^T L_{U_i^{(j)}}\right)\\ &= V_s^{-2}\sigma^2\sum_{j=1}^s \boldsymbol{u}_j^T\mathbb{E}[L_U L_U^T]\boldsymbol{u}_j \qquad \text{(B.2)}\end{aligned}$$

Plugging the formula for $L_U$ in Proposition 2.3, we have

$$\mathbb{E}[L_U L_U^T] \;=\; (\boldsymbol{\Omega}+\boldsymbol{\Sigma}_\lambda)^{-1}\mathbb{E}[(\boldsymbol{X}-\boldsymbol{A}(Z))(\boldsymbol{X}-\boldsymbol{A}(Z))^T](\boldsymbol{\Omega}+\boldsymbol{\Sigma}_\lambda)^{-1}.$$

By Lemma A.3 in Section A.1, we have that $\boldsymbol{\Sigma}_\lambda = \mathbb{E}_Z\big[\boldsymbol{B}(Z)\big(\boldsymbol{B}(Z)-\boldsymbol{A}(Z)\big)^T\big]\to 0$, and also

$$\begin{aligned}&\mathbb{E}\big[(\boldsymbol{X}-\boldsymbol{A}(Z))(\boldsymbol{X}-\boldsymbol{A}(Z))^T\big]\\ &= \mathbb{E}[(\boldsymbol{X}-\boldsymbol{B}(Z))(\boldsymbol{X}-\boldsymbol{B}(Z))^T]+\mathbb{E}[(\boldsymbol{B}(Z)-\boldsymbol{A}(Z))(\boldsymbol{B}(Z)-\boldsymbol{A}(Z))^T]\\ &\quad +2\mathbb{E}[(\boldsymbol{X}-\boldsymbol{A}(Z))(\boldsymbol{B}(Z)-\boldsymbol{A}(Z))^T]\to\boldsymbol{\Omega}^{-1}\end{aligned}$$

This implies $\mathbb{E}[L_U L_U^T]\to\boldsymbol{\Omega}^{-1}$. By (B.2), we have $\mathrm{Var}\left(V_s^{-1}\sum_{j=1}^s \varepsilon_i^{(j)}\boldsymbol{u}_j^T L_{U_i^{(j)}}\right)\to\sigma^2$. By sub-Gaussianity, we have $V_s^{-1}\sum_{j=1}^s\varepsilon_i^{(j)}\boldsymbol{u}_j^T L_{U_i^{(j)}}$ has bounded third moment, hence it satisfies Lyapunov condition. Therefore, applying central limit theorem, we have

$$\frac{1}{\sqrt{n}}\sum_{i=1}^n V_s^{-1}\sum_{j=1}^s\varepsilon_i^{(j)}\boldsymbol{u}_j^T L_{U_i^{(j)}}\rightsquigarrow N(0,\sigma^2).$$

For the second term in B.1, by Lemma A.5, we have $\|L_\lambda f_0\|_2^2\le C_2'\lambda^2$. Hence, it holds that

$$\left|\sqrt{n}\sum_{j=1}^s V_s^{-1}\boldsymbol{u}_j^T L_\lambda f_0\right|\le\tau_{\max}(\boldsymbol{\Omega})\sqrt{n}\sum_{j=1}^s\|\boldsymbol{u}_j\|_2\sqrt{C_2'}\lambda\le\sqrt{ns}\sqrt{C_2'}\tau_{\max}(\boldsymbol{\Omega})\lambda,$$

where we used Cauchy-Schwarz in both inequalities. Therefore, we have $\sqrt{n}\sum_{j=1}^s V_s^{-1}\boldsymbol{u}_j^T L_\lambda f_0 = o(1)$ if $\lambda = o(N^{-1/2})$.

For the third term in B.1, by Cauchy-Schwarz, we have

$$\left|\sqrt{n}\sum_{j=1}^s\boldsymbol{u}_j^T Rem_\beta^{(j)}\right|\le\sqrt{n}\sum_{j=1}^s\|\boldsymbol{u}_j\|_2\sup_j\|Rem_\beta^{(j)}\|_2\le\sqrt{ns}\sup_j\|Rem_\beta^{(j)}\|_2.$$

Moreover, if (3.12) holds, we can show that $d(\lambda)n^{-1/2}(J(\mathcal{F},1)+\log n) = o(1)$. Hence by Lemma 3.1 and union bound, we have

$$\mathbb{P}\big(\sup_j\|Rem_\beta^{(j)}\|_2\ge b(n,\lambda,J)\big)\le ns\exp(-\log^2 n). \qquad \text{(B.3)}$$

By explicit calculation, we can show that the conditions $\lambda = o(N^{-1/2})$ and (3.16) imply $\sqrt{ns}b(n, \lambda, J) = o(1)$. Moreover, by the assumption $s \lesssim N^{\psi}$ (stated before Assumption 3.1), we have $ns \exp(-\log^2 n) \asymp N \exp(-c \log^2 N) = o(1)$. Hence it follows from (B.3) that $\left|\sqrt{n} \sum_{j=1}^{s} V_s^{-1} \boldsymbol{u}_j^T Rem_{\beta}^{(j)}\right| = o_P(1)$. This completes the proof □

### B.3. Proof of Theorem 3.7.

PROOF. By first order optimality condition, we have

$$\sum_{i \in S_j} \boldsymbol{X}_i \big(Y_i - \boldsymbol{X}_i^T \check{\boldsymbol{\beta}}^{(j)} - \bar{f}(Z_i)\big) = 0.$$

Hence we have

$$\check{\boldsymbol{\beta}}^{(j)} = \Big(\sum_{i \in S_j} \boldsymbol{X}_i \boldsymbol{X}_i^T\Big)^{-1} \sum_{i \in S_j} \boldsymbol{X}_i \big(Y_i - \bar{f}(Z_i)\big). \tag{B.4}$$

As for $i \in S_j$ we have $Y_i = \boldsymbol{X}_i^T \boldsymbol{\beta}_0^{(j)} + f_0(Z_i) + \varepsilon_i$, hence

$$\begin{aligned} \sqrt{n}(\check{\boldsymbol{\beta}}^{(j)} - \boldsymbol{\beta}_0) = n^{-1/2} & \sum_{i \in S_j} (\widehat{\boldsymbol{\Sigma}}^{(j)})^{-1} \boldsymbol{X}_i \varepsilon_i \\ & + n^{-1/2} \sum_{i \in S_j} (\widehat{\boldsymbol{\Sigma}}^{(j)})^{-1} \boldsymbol{X}_i \big(f_0(Z_i) - \bar{f}(Z_i)\big), \end{aligned} \tag{B.5}$$

where the $\widehat{\boldsymbol{\Sigma}}^{(j)} = \frac{1}{n} \sum_{i \in S_j} \boldsymbol{X}_i \boldsymbol{X}_i^T$ is sample covariance of $\boldsymbol{X}$ based on data from the $j$-th subpopulation. For the first term on the R.H.S. of (B.5), it is the same as the one for ordinary least squares, and so

$$n^{-1/2} \sum_{i \in S_j} (\widehat{\boldsymbol{\Sigma}}^{(j)})^{-1} \boldsymbol{X}_i \varepsilon_i \rightsquigarrow N(0, \sigma^2 \boldsymbol{\Sigma}^{-1}). \tag{B.6}$$

For the second term on the R.H.S. of (B.5), by triangular inequality, we have for all $1 \le j \le s$ that

$$\begin{aligned} \Big\| \frac{1}{\sqrt{n}} \sum_{i \in S_j} \big(\widehat{\boldsymbol{\Sigma}}^{(j)}\big)^{-1} \boldsymbol{X}_i \big(f_0(Z_i) - \bar{f}(Z_i)\big) \Big\|_2 & \le \|\bar{f} - f_0\|_{\sup} n^{-1/2} \sum_{i \in S_j} \|\widehat{\boldsymbol{\Sigma}}^{(j)} \boldsymbol{X}_i\|_2 \\ & \le C n^{1/2} \|\bar{f} - f_0\|_{\sup}, \end{aligned}$$

for some constant $C$ that is related to the dimension $p$ and boundedness of the support of $\mathcal{X}$. As $\|\bar{f} - f_0\|_{\sup} \le c_{\phi} d(\lambda)^{1/2} \|\bar{f} - f_0\|_{\mathcal{C}}$ by Lemma A.1, it suffices to bound $\|\bar{f} - f_0\|_{\mathcal{C}}$. By (7.3) and triangular inequality,

$$\|\bar{f} - f_0\|_{\mathcal{C}} \le \Big\| \frac{1}{N} \sum_{i=1}^{N} N_{U_i} \varepsilon_i \Big\|_{\mathcal{C}} + \|N_{\lambda} f_0\|_{\mathcal{C}} + \Big\| \frac{1}{s} \sum_{j=1}^{s} Rem_f^{(j)} \Big\|_{\mathcal{C}}. \tag{B.7}$$

For the first term on the R.H.S., define $\mathcal{Q}_i = \{|\varepsilon_i| \le \log N\}$. Since $\|N_U\|_{\mathcal{C}} \le c_\phi d(\lambda)^{1/2}$, we have that $\{\varepsilon_i N_{U_i} I_{\mathcal{Q}_i}\}_{i=1}^N$ is a sequence of random variables in Hilbert space $\mathcal{F}$ that are i.i.d. with mean zero and bounded by $c_\phi d(\lambda)^{1/2}\log N$. Therefore by Lemma G.1 we have

$$\begin{aligned}
&\mathbb{P}\Big(\|\frac{1}{N}\sum_{i=1}^N \varepsilon_i N_{U_i}\|_{\mathcal{C}} > c_\phi \log^2 N\sqrt{d(\lambda)/N}\Big)\\
&\quad\le \mathbb{P}\Big(\{\cap_i \mathcal{Q}_i\}\cap\{\|\frac{1}{N}\sum_{i=1}^N \varepsilon_i N_{U_i}\|_{\mathcal{C}} > c_\phi \log^2 N\sqrt{d(\lambda)/N}\}\Big) + \mathbb{P}\big((\cap_i \mathcal{Q}_i)^c\big)\\
&\quad\le 2\exp(-\log^2 N) + 2N\exp(-\log^2 N) \to 0.
\end{aligned} \tag{B.8}$$

Therefore, we have $\|N^{-1}\sum_{i=1}^N \varepsilon_i N_{U_i}\|_{\mathcal{C}} = o_P\big(\log^2 N\sqrt{d(\lambda)/N}\big)$. Moreover, by (A.27), we have $\|N_\lambda f_0\|_{\mathcal{C}} = O(\lambda^{-1/2})$. Furthermore, by Lemma 7.2, when $s$ satisfies Condition (3.12) , we have $\|\frac{1}{s}\sum_{j=1}^s Rem_f^{(j)}\|_{\mathcal{C}} = o_P(s^{-1/2}b(n,\lambda,J)\log N)$. By the definition of $b(n,\lambda,J)$, we have if $s$ satisfies (3.13), then $\|\frac{1}{s}\sum_{j=1}^s Rem_f^{(j)}\|_{\mathcal{C}} = o_P(\sqrt{d(\lambda)/N})$. Therefore, by (B.7), it holds

$$\|\bar f - f_0\|_{\mathcal{C}} = o_P\big(\log^2 N\sqrt{d(\lambda)/N} + \lambda^{1/2}\big).$$

Hence we have

$$\begin{aligned}
&\max_{1\le j\le s}\|\frac{1}{\sqrt{n}}\sum_{i\in S_j}\big(\widehat{\mathbf{\Sigma}}^{(j)}\big)^{-1}\boldsymbol{X}_i\big(f_0(Z_i) - \bar f(Z_i)\big)\|_2\\
&\qquad = o_P\big(n^{1/2}d(\lambda)^{1/2}(\log^2 N\sqrt{d(\lambda)/N} + \lambda^{1/2})\big).
\end{aligned} \tag{B.9}$$

Plugging the relationship $n = N/s$, we have when $s \gtrsim d(\lambda)^2\log^4 N$ and $\lambda = O(d(\lambda)/N)$,

$$\|\frac{1}{\sqrt{n}}\sum_{i\in S_j}\big(\widehat{\mathbf{\Sigma}}^{(j)}\big)^{-1}\boldsymbol{X}_i\big(f_0(Z_i) - \bar f(Z_i)\big)\|_2 = o_P(1).$$

Hence by (B.5) and (B.6),

$$\sqrt{n}(\check{\boldsymbol{\beta}} - \boldsymbol{\beta}_0) \rightsquigarrow N(0, \sigma^2\mathbf{\Sigma}^{-1}),$$

as desired. □

### B.4. Proof of Theorem 3.8.

Proof. (i) By (3.5), we have for the $j$th sub-population

$$\widehat{\boldsymbol{\beta}}^{(j)} - \boldsymbol{\beta}_0^{(j)} = \frac{1}{n}\sum_{i\in S_j} L_{U_i}\varepsilon_i - L_\lambda f_0 - Rem_\beta^{(j)}, \tag{B.10}$$

where $Rem_{\beta}^{(j)} = 1/n \sum_{i \in S_j} \Big( L_{U_i} \Delta m^{(j)}(U_i) - \mathbb{E}_U[L_U \Delta m^{(j)}(U)] \Big)$. Equation (B.10) also holds for the $k$-th sub-population. Hence under $H_0 : \boldsymbol{\beta}_0^{(j)} = \boldsymbol{\beta}_0^{(k)}$, we have

$$\widehat{\boldsymbol{\beta}}^{(j)} - \widehat{\boldsymbol{\beta}}^{(k)} = \frac{1}{n} \sum_{i \in S_j} L_{U_i} \varepsilon_i - \frac{1}{n} \sum_{i \in L_k} L_{U_i} \varepsilon_i - (Rem_{\beta}^{(j)} - Rem_{\beta}^{(k)}), \tag{B.11}$$

By independence between two sub-populations, we have $\frac{1}{n} \sum_{i \in S_j} L_{U_i} \varepsilon_i - \frac{1}{n} \sum_{i \in L_k} L_{U_i} \varepsilon_i \rightsquigarrow N(\mathbf{0}, 2\sigma^2 \boldsymbol{\Omega}^{-1})$. Moreover, when the conditions in Theorem 3.6 are satisfied, we have $\sqrt{n} \|Q\| \|Rem_{\beta}^{(j)} - Rem_{\beta}^{(j)}\|_2 = o_P(1)$ by triangular inequality. Therefore the result follows.

(ii) By (B.5), we have

$$\sqrt{n}(\check{\boldsymbol{\beta}}^{(j)} - \boldsymbol{\beta}_0^{(j)}) = \frac{1}{\sqrt{n}} \big(\widehat{\boldsymbol{\Sigma}}^{(j)}\big)^{-1} \sum_{i \in S_j} \boldsymbol{X}_i \varepsilon_i + \frac{1}{\sqrt{n}} \big(\widehat{\boldsymbol{\Sigma}}^{(j)}\big)^{-1} \sum_{i \in S_j} \boldsymbol{X}_i \big( f_0(Z_i) - \bar{f}(Z_i) \big), \tag{B.12}$$

where $\widehat{\boldsymbol{\Sigma}}^{(j)} = \frac{1}{n} \sum_{i \in S_j} \boldsymbol{X}_i \boldsymbol{X}_i^T$. The above equation is also true for $k$-th sub-population. So if $s$ satisfies Condition (3.24), (3.12) and (3.13), we have

$$\Big\| \frac{1}{\sqrt{n}} \big(\widehat{\boldsymbol{\Sigma}}^{(j)}\big)^{-1} \sum_{i \in S_j} \boldsymbol{X}_i \big( f_0(Z_i) - \bar{f}(Z_i) \big) \Big\|_2 = o_P(1).$$

We have another equation that is same as (B.12) with $j$ replaced by $k$. Hence subtracting the two equations, we have under $H_0 : \boldsymbol{\beta}_0^{(j)} = \boldsymbol{\beta}_0^{(k)}$,

$$\sqrt{n}(\check{\boldsymbol{\beta}}^{(j)} - \check{\boldsymbol{\beta}}^{(k)}) = \frac{1}{\sqrt{n}} \big(\widehat{\boldsymbol{\Sigma}}^{(j)}\big)^{-1} \sum_{i \in S_j} \boldsymbol{X}_i \varepsilon_i - \frac{1}{\sqrt{n}} \big(\widehat{\boldsymbol{\Sigma}}^{(k)}\big)^{-1} \sum_{i \in L_k} \boldsymbol{X}_i \varepsilon_i + o_P(1).$$

Hence the conclusion follows from CLT and independence of sub-populations $j$ and $k$. □

**B.5. Proof of Theorem 3.9.** Before presenting the proof, we define the following preliminaries: for any $\mathcal{G} \subset \{1, 2, \ldots, s\}$ with $|\mathcal{G}| = d$, let

$$T_{0,\mathcal{G}} := \max_{j \in \mathcal{G}, 1 \le k \le p} \frac{1}{\sqrt{n}} \sum_{i \in S_j} \big(\boldsymbol{\Sigma}^{-1}\big)_k \boldsymbol{X}_i \varepsilon_i,$$

where $\big(\boldsymbol{\Sigma}^{-1}\big)_k$ denotes the $k$-th row of the precision matrix $\boldsymbol{\Sigma}^{-1}$ of $\boldsymbol{X}$. Furthermore, let

$$W_{0,\mathcal{G}} := \max_{j \in \mathcal{G}, 1 \le k \le p} n^{-1/2} \sum_{i \in S_j} \Gamma_{i,k},$$

where $\{\mathbf{\Gamma}_i = (\Gamma_{i,1}, \dots \Gamma_{i,p})\}$ for each $i \in S_j$, $j \in \mathcal{G}$ is a sequence of mean zero independent Gaussian vector with $\mathbb{E}[\mathbf{\Gamma}_i \mathbf{\Gamma}_i^T] = (\mathbf{\Sigma})^{-1}\sigma^2$. Lastly, it is useful to recall

$$W_{\mathcal{G}} := \max_{j \in \mathcal{G}, 1 \le k \le p} \frac{1}{\sqrt{n}} \sum_{i \in S_j} \big(\widehat{\mathbf{\Sigma}}^{(j)}\big)_k^{-1} \boldsymbol{X}_i e_i,$$

and $c_{\mathcal{G}}(\alpha) = \inf\{t \in \mathbb{R} : \mathbb{P}(W_{\mathcal{G}} \le t \,|\, \mathbb{X}) \ge 1 - \alpha\}$.

The proof strategy is similar to that of Theorem 3.2 in Chernozhukov et al. (2013). Specifically, we first approximate $T_{\mathcal{G}}$ by $T_{0,\mathcal{G}}$, and then apply Gaussian approximation to $T_{0,\mathcal{G}}$ and $W_{0,\mathcal{G}}$. Then, we argue that $W_{0,\mathcal{G}}$ and $W_{\mathcal{G}}$ are close. Hence we can approximate the quantiles of $T_{\mathcal{G}}$ by those of $W_{\mathcal{G}}$. The detailed proof is presented as follows.

PROOF. By (B.5), we have

$$\sqrt{n}(\check{\boldsymbol{\beta}}^{(j)} - \widetilde{\boldsymbol{\beta}}^{(j)}) = n^{-1/2} \sum_{i \in S_j} (\widehat{\mathbf{\Sigma}}^{(j)})^{-1} \boldsymbol{X}_i \varepsilon_i + \mathbf{\Delta}^{(j)}, \tag{B.13}$$

where

$$\mathbf{\Delta}^{(j)} = n^{-1/2} \sum_{i \in S_j} (\widehat{\mathbf{\Sigma}}^{(j)})^{-1} \boldsymbol{X}_i \big(f_0(Z_i) - \bar{f}(Z_i)\big).$$

By (B.9) in the proof of Theorem 3.7, we have

$$\max_{j \in \mathcal{G}} \|\mathbf{\Delta}^{(j)}\|_\infty = o_P\big(n^{1/2} d(\lambda)^{1/2} (\log^2 N \sqrt{d(\lambda)/N} + \lambda^{1/2})\big).$$

and when $s \gtrsim d(\lambda)^2 \log(pd) \log^4 N$ and $\lambda = O(d(\lambda)/N)$, we have

$$\max_{j \in \mathcal{G}} \|\mathbf{\Delta}^{(j)}\|_\infty = o_P(\log^{-1/2}(pd)).$$

By the definitions of $T_{\mathcal{G}}$ and $T_{0,\mathcal{G}}$ and (B.13), we have

$$|T_{\mathcal{G}} - T_{0,\mathcal{G}}| \le \max_{j \in \mathcal{G}} \frac{1}{\sqrt{n}} \Big\| \sum_{i \in S_j} (\widehat{\mathbf{\Sigma}}^{(j)})^{-1} \boldsymbol{X}_i \varepsilon_i - \mathbf{\Sigma}^{-1} \boldsymbol{X}_i \varepsilon_i \Big\|_\infty + \max_{j \in \mathcal{G}} \|\mathbf{\Delta}^{(j)}\|_\infty,$$

where we used the fact that $\max_j a_j - \max_j b_j \le \max_j |a_j - b_j|$ for any two finite sequences $\{a_j\}, \{b_j\}$. By the above inequality and Lemma B.2, there exist $\zeta_1$ and $\zeta_2$ such that

$$\mathbb{P}\big(|T_{\mathcal{G}} - T_{0,\mathcal{G}}| \ge \zeta_1\big) \le \zeta_2, \tag{B.14}$$

where $\zeta_1 \sqrt{1 \vee \log(pd/\zeta_1)} = o(1)$ and $\zeta_2 = o(1)$.

We next turn to bound the distance between quantiles of $W_{\mathcal{G}}$ and $W_{0,\mathcal{G}}$. Let $c_{0,\mathcal{G}}(\alpha) := \inf\{t \in \mathbb{R} : \mathbb{P}(W_{0,\mathcal{G}} \le t) \ge 1 - \alpha\}$, and let $\pi(\nu) := C_2 \nu^{1/3}(1 \vee \log(pd/\nu))^{2/3}$ with $C_2 > 0$, and

$$\Psi := \max_{\substack{1 \le k, \ell \le p \\ j \in \mathcal{G}}} \sigma^2 |(\widehat{\mathbf{\Sigma}}^{(j)} - \mathbf{\Sigma})_{k\ell}|.$$

As the data size in each subpopulation is the same, we can relabel $\{\boldsymbol{X}_i \in \mathbb{R}^p\}_{i\in S_j, j\in\mathcal{G}}$ as $\{\boldsymbol{X}_i^{(j)} \in \mathbb{R}^p\}_{1\le i\le n, j\in\mathcal{G}}$ and $\{\Gamma_i \in \mathbb{R}^p\}_{1\le i\le n, j\in\mathcal{G}}$ as $\{\Gamma_i^{(j)} \in \mathbb{R}^p\}_{1\le i\le n, j\in\mathcal{G}}$. Then we can re-write $W_{0,\mathcal{G}} = \max_{j\in\mathcal{G}, 1\le k\le p} U_k^{(j)}$, and $W_{\mathcal{G}} = \max_{j\in\mathcal{G}, 1\le k\le p} V_k^{(j)}$, where

$$U_k^{(j)} = \frac{1}{\sqrt{n}}\sum_{i=1}^n \Gamma_{ik}^{(j)} \quad \text{and} \quad V_k^{(j)} = \frac{1}{\sqrt{n}}\sum_{i=1}^n (\widehat{\boldsymbol{\Sigma}}^{(j)})_k^{-1}\boldsymbol{X}_i^{(j)} e_i^{(j)}.$$

($e_i^{(j)}$ is defined in the similar way). Notice that $\boldsymbol{U} = \{U_k^{(j)}\}_{1\le k\le p, j\in\mathcal{G}}$ can be viewed as an $(p\cdot d)$-dimensional Gaussian random vector with mean zero and covariance

$$\begin{pmatrix} \sigma^2\boldsymbol{\Sigma}^{-1} & 0 & \dots & 0 \\ 0 & \sigma^2\boldsymbol{\Sigma}^{-1} & \dots & 0 \\ \vdots & \vdots & \ddots & \vdots \\ 0 & 0 & \dots & \sigma^2\boldsymbol{\Sigma}^{-1} \end{pmatrix} \in \mathbb{R}^{(p\cdot d)\times(p\cdot d)}.$$

Conditioned on $\mathbb{X}$, $\boldsymbol{V} = \{V_k^{(j)}\}_{1\le k\le p, j\in\mathcal{G}}$ can be viewed as an $(p\cdot d)$-dimensional Gaussian random vector with mean zero and covariance

$$\begin{pmatrix} \sigma^2\big(\widehat{\boldsymbol{\Sigma}}^{(1)}\big)^{-1} & 0 & \dots & 0 \\ 0 & \sigma^2\big(\widehat{\boldsymbol{\Sigma}}^{(2)}\big)^{-1} & \dots & 0 \\ \vdots & \vdots & \ddots & \vdots \\ 0 & 0 & \dots & \sigma^2\big(\widehat{\boldsymbol{\Sigma}}^{(d)}\big)^{-1} \end{pmatrix} \in \mathbb{R}^{(p\cdot d)\times(p\cdot d)}.$$

Using Gaussian comparison (Lemma 3.1 in Chernozhukov et al. (2013)) and applying the same argument as in the proof of Lemma 3.2 in Chernozhukov et al. (2013), we obtain for any $\nu > 0$

$$\mathbb{P}\big(c_{0,\mathcal{G}}(\alpha) \le c_{\mathcal{G}}(\alpha + \pi(\nu))\big) \ge 1 - \mathbb{P}(\Psi > \nu), \tag{B.15}$$

$$\mathbb{P}\big(c_{\mathcal{G}}(\alpha) \le c_{0,\mathcal{G}}(\alpha + \pi(\nu))\big) \ge 1 - \mathbb{P}(\Psi > \nu). \tag{B.16}$$

By Lemma B.1, we have

$$\sup_{\alpha\in(0,1)} \big|\mathbb{P}\big(T_{0,\mathcal{G}} > c_{\mathcal{G}}(\alpha)\big) - \alpha\big| \le \sup_{\alpha\in(0,1)} \big|\mathbb{P}\big(W_{0,\mathcal{G}} > c_{\mathcal{G}}(\alpha)\big) - \alpha\big| + n^{-c}. \tag{B.17}$$

To further control $\mathbb{P}\big(W_{0,\mathcal{G}} > c_{\mathcal{G}}(\alpha)\big)$, we define $\mathcal{E}_1 = \{c_{0,\mathcal{G}}(\alpha - \pi(\nu)) \le c_{\mathcal{G}}(\alpha)\}$, $\mathcal{E}_2 = \{c_{\mathcal{G}}(\alpha) \le c_{0,\mathcal{G}}(\alpha + \pi(\nu))\}$. We have

$$\begin{aligned} \mathbb{P}\big(W_{0,\mathcal{G}} > c_{\mathcal{G}}(\alpha)\big) &= \mathbb{P}\big(W_{0,\mathcal{G}} > c_{\mathcal{G}}(\alpha), \mathcal{E}_1\big) + \mathbb{P}\big(W_{0,\mathcal{G}} > c_{\mathcal{G}}(\alpha), \mathcal{E}_1^c\big) \\ &\le \mathbb{P}\big(W_{0,\mathcal{G}} > c_{0,\mathcal{G}}(\alpha - \pi(\nu))\big) + \mathbb{P}(\mathcal{E}_1^c) \\ &\le \alpha - \pi(\nu) + \mathbb{P}(\Psi > \nu), \end{aligned}$$

where the last inequality is by the definition of $c_{0,\mathcal{G}}(\alpha)$ and (B.15). Similarly,

we have

$$\begin{aligned}\mathbb{P}\big(W_{0,\mathcal{G}} > c_{\mathcal{G}}(\alpha)\big) &= 1 - \mathbb{P}\big(W_{0,\mathcal{G}} \le c_{\mathcal{G}}(\alpha)\big) \\ &= 1 - \mathbb{P}\big(W_{0,\mathcal{G}} \le c_{\mathcal{G}}(\alpha), \mathcal{E}_2\big) - \mathbb{P}\big(W_{0,\mathcal{G}} \le c_{\mathcal{G}}(\alpha), \mathcal{E}_2^c\big) \\ &\ge 1 - \mathbb{P}\big(W_{0,\mathcal{G}} \le c_{0,\mathcal{G}}(\alpha + \pi(\nu))\big) - \mathbb{P}(\mathcal{E}_2^c) \\ &\ge \alpha + \pi(\nu) - \mathbb{P}(\Psi > \nu),\end{aligned}$$

where the last inequality is by the definition of $c_{0,\mathcal{G}}(\alpha)$ and (B.16). Hence it follows from (B.17) that

$$\sup_{\alpha\in(0,1)} \big|\mathbb{P}\big(T_{0,\mathcal{G}} > c_{\mathcal{G}}(\alpha)\big) - \alpha\big| \le \pi(\nu) + \mathbb{P}(\Psi > \nu) + n^{-c}. \tag{B.18}$$

Define the event $\mathcal{E}_3 = \{|T_{0,\mathcal{G}} - T_{\mathcal{G}}| \le \zeta_1\}$. By (B.14), we have $\mathbb{P}(\mathcal{E}_3^c) \le \zeta_2$. Hence, we deduce that for any $\alpha$

$$\begin{aligned}\mathbb{P}(T_{\mathcal{G}} \ge c_G(\alpha)) - \alpha &\le \mathbb{P}(T_{\mathcal{G}} \ge c_G(\alpha), \mathcal{E}_3) + \mathbb{P}(\mathcal{E}_3^c) - \alpha \\ &\le \mathbb{P}(T_{0,\mathcal{G}} \ge c_G(\alpha) - \zeta_1) + \zeta_2 - \alpha, \\ &\le \mathbb{P}(T_{0,\mathcal{G}} \ge c_G(\alpha)) + C\zeta_1\sqrt{1 \vee \log(ps/\zeta_1)} + \zeta_2 - \alpha \\ &\le \pi(\nu) + \mathbb{P}(\Psi > \nu) + n^{-c} + C\zeta_1\sqrt{1 \vee \log(ps/\zeta_1)} + \zeta_2,\end{aligned}$$

where the second last inequality is by Corollary 16 of Wasserman (2014) (Gaussian anti-concentration). By similar arguments, we get the same bound for $\alpha - \mathbb{P}(T_{\mathcal{G}} \ge c_G(\alpha))$, so we have

$$\sup_{\alpha} \big|\mathbb{P}(T_{\mathcal{G}} \ge c_G(\alpha)) - \alpha\big| \le \pi(\nu) + \mathbb{P}(\Psi > \nu) + n^{-c} + C\zeta_1\sqrt{1 \vee \log(pd/\zeta_1)} + \zeta_2.$$

Lastly, we bound $\Psi$. By (B.21) in the proof of Lemma B.2 and the fact that elementwise infinity norm is bounded by spectral norm, we obtain

$$\Psi \le \max_{j\in\mathcal{G}} \|\widehat{\mathbf{\Sigma}}^{(j)} - \mathbf{\Sigma}\|_\infty = o_P\Big(p\sqrt{(\log d)/n}\Big).$$

Hence, choosing $\nu = p\sqrt{(\log d)/n}$, we get

$$\sup_{\alpha} \big|\mathbb{P}(T_{\mathcal{G}} \ge c_G(\alpha)) - \alpha\big| = o(1),$$

which concludes the proof. □

**Lemma B.1.** Suppose Assumption 3.1 holds. For any $\mathcal{G} \subset \{1, 2, \ldots, s\}$ with $|\mathcal{G}| = d$, if $(\log(pdn))^7/n \le C_1 n^{-c_1}$ for some constants $c_1, C_1 > 0$, then we have

$$\sup_{x\in\mathbb{R}} \Big|\mathbb{P}\big(T_{0,\mathcal{G}} \le x\big) - \mathbb{P}\big(W_{0,\mathcal{G}} \le x\big)\Big| \le n^{-c},$$

for some constant $c > 0$.

Proof. As the data size in each subpopulation is the same, we can relabel

$\{\boldsymbol{X}_i \in \mathbb{R}^p\}_{i\in S_j, j\in\mathcal{G}}$ as $\{\boldsymbol{X}_i^{(j)} \in \mathbb{R}^p\}_{1\le i\le n, j\in\mathcal{G}}$. Under such notation, we have $T_{0,\mathcal{G}} = \max_{j\in\mathcal{G},1\le k\le p} n^{-1/2}\sum_{i=1}^n \xi_{ik}^{(j)}$, where $\xi_{ik}^{(j)} = (\boldsymbol{\Sigma}^{-1})_k \boldsymbol{X}_i^{(j)}\varepsilon_i$. For each $i$, $\{\xi_{ik}^{(j)}\}_{1\le k\le p, j\in\mathcal{G}}$ can be viewed as a $(p\cdot d)$-dimensional vector with covariance matrix

$$\begin{pmatrix} \sigma^2\boldsymbol{\Sigma}^{-1} & 0 & \dots & 0 \\ 0 & \sigma^2\boldsymbol{\Sigma}^{-1} & \dots & 0 \\ \vdots & \vdots & \ddots & \vdots \\ 0 & 0 & \dots & \sigma^2\boldsymbol{\Sigma}^{-1} \end{pmatrix} \in \mathbb{R}^{(p\cdot d)\times(p\cdot d)}.$$

The same thing can be done for $\boldsymbol{\Gamma}$ which results in the same covariance matrix. Then we apply Corollary 2.1 of Chernozhukov et al. (2013) to prove the Gaussian approximation result stated in the lemma. It suffices to verify Condition (E1) therein. We have $\mathbb{E}\big[(\xi_{ik}^{(j)})^2\big] = (\boldsymbol{\Sigma}^{-1})_{kk}$ is a constant, and $\max_{\ell=1,2}\mathbb{E}\big[|\xi_{ik}^{(j)}|^{2+\ell}/B^\ell\big] + \mathbb{E}\big[\exp(|\xi_{ik}^{(j)}|/B)\big] \le 4$ for some large enough constant $B$, by the sub-Gaussianity of $\varepsilon_i^{(j)}$ and the boundedness of $\boldsymbol{X}_i^{(j)}$. Hence Condition (E1) is verified, and by the assumption that $(\log(pdn))^7/n \le C_1 n^{-c_1}$, we get the desired result. □

**Lemma B.2.** Suppose Assumption 3.1 holds. For any $\mathcal{G} \subset \{1, 2, \dots, s\}$ with $|\mathcal{G}| = d$, suppose $p^2\log(pd)/\sqrt{n} = o(1)$. Then there exist $\zeta_1$ and $\zeta_2$ such that

$$\mathbb{P}\bigg(\max_{j\in\mathcal{G}}\frac{1}{\sqrt{n}}\Big\|\sum_{i\in S_j}(\widehat{\boldsymbol{\Sigma}}^{(j)})^{-1}\boldsymbol{X}_i\varepsilon_i - \boldsymbol{\Sigma}^{-1}\boldsymbol{X}_i\varepsilon_i\Big\|_\infty > \zeta_1\bigg) \le \zeta_2,$$

where $\zeta_1\sqrt{1\vee\log(pd/\zeta_1)} = o(1)$ and $\zeta_2 = o(1)$.

Proof. We have

$$\begin{aligned} &\max_{j\in\mathcal{G}}\frac{1}{\sqrt{n}}\Big\|\sum_{i\in S_j}(\widehat{\boldsymbol{\Sigma}}^{(j)})^{-1}\boldsymbol{X}_i\varepsilon_i - \boldsymbol{\Sigma}^{-1}\boldsymbol{X}_i\varepsilon_i\Big\|_\infty \\ &\qquad\le \max_{j\in\mathcal{G}}\big\|(\widehat{\boldsymbol{\Sigma}}^{(j)})^{-1} - \boldsymbol{\Sigma}^{-1}\big\|_1 \max_{j\in\mathcal{G}}\Big\|\frac{1}{\sqrt{n}}\sum_{i\in S_j}\boldsymbol{X}_i\varepsilon_i\Big\|_\infty \\ &\qquad\le \max_{j\in\mathcal{G}} p\big\|(\widehat{\boldsymbol{\Sigma}}^{(j)})^{-1} - \boldsymbol{\Sigma}^{-1}\big\| \max_{j\in\mathcal{G}}\Big\|\frac{1}{\sqrt{n}}\sum_{i\in S_j}\boldsymbol{X}_i\varepsilon_i\Big\|_\infty, \end{aligned} \tag{B.19}$$

where $\|\cdot\|_1$ denotes the elementwise $L_1$ norm of matrices. As $\varepsilon_i$ are i.i.d. sub-Gaussian random variables, we have by Hoeffding's inequality that for any $j \in \mathcal{G}$ and $1 \le k \le p$

$$\mathbb{P}\Big(\frac{1}{\sqrt{n}}\sum_{i\in S_j} X_{ik}\varepsilon_i > t \,|\, \mathbb{X}\Big) \le \exp\Big(-\frac{nt^2}{\sum_{i\in S_j} X_{ik}^2\sigma^2}\Big) \le \exp\Big(-\frac{t^2}{c_x^2\sigma^2}\Big),$$

where the second inequality is by the boundedness of $X_{ik}$. By law of iterated expectation and union bound we have

$$\mathbb{P}\Big(\max_{j\in\mathcal{G}}\Big\|\frac{1}{\sqrt{n}}\sum_{i\in S_j}\boldsymbol{X}_i\varepsilon_i\Big\|_\infty > t\Big) \le pd\exp\Big(-\frac{t^2}{c_x^2\sigma^2}\Big).$$

Letting $t = 2c_x\sigma\sqrt{\log(pd)}$, we get with probability at least $1-(pd)^{-1}$ that

$$\max_{j\in\mathcal{G}}\Big\|\frac{1}{\sqrt{n}}\sum_{i\in S_j}\boldsymbol{X}_i\varepsilon_i\Big\|_\infty \le 2c_x\sigma\sqrt{\log(pd)}. \tag{B.20}$$

By the boundedness of $\boldsymbol{X}$, we have $\|\boldsymbol{X}_i\boldsymbol{X}_i^T - \mathbb{E}[\boldsymbol{X}_i\boldsymbol{X}_i^T]\| \le 2\|\boldsymbol{X}_i\boldsymbol{X}_i^T\| \le 2\|\boldsymbol{X}_i\|_2^2 \le 2pc_x^2$. Therefore, by Lemma G.3, we have for all $j \in \mathcal{G}$ that

$$\mathbb{P}\Big(\|\widehat{\boldsymbol{\Sigma}}^{(j)} - \boldsymbol{\Sigma}\| \ge t\Big) \le \mathbb{P}\Big(\Big\|\frac{1}{n}\sum_{i\in S_j}\boldsymbol{X}_i\boldsymbol{X}_i^T - \mathbb{E}[\boldsymbol{X}_i\boldsymbol{X}_i^T]\Big\| \ge t\Big) \le p\exp\Big(-\frac{nt^2}{32p^2c_x^4}\Big).$$

and so it follows from union bound that

$$\mathbb{P}\Big(\max_{j\in\mathcal{G}}\|\widehat{\boldsymbol{\Sigma}}^{(j)} - \boldsymbol{\Sigma}\| \ge t\Big) \le pd\exp\Big(-\frac{nt^2}{32p^2c_x^4}\Big).$$

Choosing $t = 64p\sqrt{(\log d)/n}$, we obtain

$$\max_{j\in\mathcal{G}}\|\widehat{\boldsymbol{\Sigma}}^{(j)} - \boldsymbol{\Sigma}\| = o_P\Big(p\sqrt{\frac{\log d}{n}}\Big). \tag{B.21}$$

Thus, by Lemma G.4, we get

$$\max_{j\in\mathcal{G}}\|(\widehat{\boldsymbol{\Sigma}}^{(j)})^{-1} - \boldsymbol{\Sigma}^{-1}\| = o_P\Big(p\sqrt{\frac{\log d}{n}}\Big). \tag{B.22}$$

Combining (B.19), (B.20) and (B.22), we have

$$\max_{j\in\mathcal{G}}\frac{1}{\sqrt{n}}\Big\|\sum_{i\in S_j}(\widehat{\boldsymbol{\Sigma}}^{(j)})^{-1}\boldsymbol{X}_i\varepsilon_i - \boldsymbol{\Sigma}^{-1}\boldsymbol{X}_i\varepsilon_i\Big\|_\infty = o_P\Big(p^2\frac{\log(pd)}{\sqrt{n}}\Big)$$

We choose $\zeta_1$ such that $p^2\log(pd)/(\sqrt{n}\zeta_1) = o(1)$ and $\zeta_1\sqrt{1\vee\log(pd/\zeta_1)} =$

$o(1)$, e.g., $\zeta_1^2 = p^2 \log(pd)/\sqrt{n}$. Then by the above equation we have

$$\mathbb{P}\bigg(\max_{j\in\mathcal{G}} \frac{1}{\sqrt{n}}\bigg\|\sum_{i\in S_j}(\widehat{\mathbf{\Sigma}}^{(j)})^{-1}\boldsymbol{X}_i\varepsilon_i - \mathbf{\Sigma}^{-1}\boldsymbol{X}_i\varepsilon_i\bigg\|_\infty \geq \zeta_1\bigg) < \zeta_2,$$

where $\zeta_2 = o(1)$. □

## APPENDIX C: PROOFS IN SECTION 4

### C.1. Proof of Corollary 4.1.

PROOF. We begin by computing $d(\lambda)$. As $\mu_i = 0$ for $i > r$, we have that $d(\lambda) = \sum_{i=1}^r \frac{1}{1+\lambda/\mu_i} \asymp r$. Therefore by Theorem 3.6, $\lambda = o\big(\sqrt{d(\lambda)/N} \wedge n^{-1/2}\big) = o(N^{-1/2})$. We next calculate the asymptotic covariance.

$$\begin{aligned} A_k(z_0) &= \langle A_k, \widetilde{K}_{z_0}\rangle_{\mathcal{C}} = \langle B_k, \widetilde{K}_{z_0}\rangle_{L_2(\mathbb{P}_Z)} \\ &= \sum_{i=1}^r \frac{\langle B_k, \phi_i\rangle_{L_2(\mathbb{P}_Z)}}{1+\lambda/\mu_i}\phi_i(z_0) \to \sum_{i=1}^r \langle B_k, \phi_i\rangle_{L_2(\mathbb{P}_Z)}\phi_i(z_0). \end{aligned}$$

Hence $\boldsymbol{\gamma}_{z_0} = d(\lambda)^{-1/2}\sum_{i=1}^r \langle B_k, \phi_i\rangle_{L_2(\mathbb{P}_Z)}\phi_i(z_0)$. The formula for $\mathbf{\Sigma}_{12}^*$ and $\Sigma_{22}^*$ then follows from Theorem 3.6.

We next calculate the entropy integral $J(\mathcal{F}, \delta)$ for finite rank RKHS and the upper bound for $s$. Define $\widetilde{\mathcal{F}}_2 = \{f \in \mathcal{H} : \|f\|_{\sup} \leq 1, \|f\|_{\mathcal{H}} \leq 1\}$. By Carl and Triebel (1980), for finite rank RKHS,

$$\log \mathcal{N}(\widetilde{\mathcal{F}}_2, \|\cdot\|_{\sup}, \delta) \asymp r \log \delta^{-1}.$$

We have that $\mathcal{N}(\mathcal{F}, \|\cdot\|_{\sup}, \delta) \leq \mathcal{N}(\mathcal{F}_1, \|\cdot\|_{\sup}, \delta)\mathcal{N}(\mathcal{F}_2, \|\cdot\|_{\sup}, \delta)$. As $\mathcal{N}(\mathcal{F}_1, \|\cdot\|_{\sup}, \delta)$ is dominated by $\mathcal{N}(\mathcal{F}_2, \|\cdot\|_{\sup}, \delta)$, it suffices to bound $\mathcal{N}(\mathcal{F}_2, \|\cdot\|_{\sup}, \delta)$. Now by Van Der Vaart and Wellner (1996), we have that

$$\begin{aligned} \mathcal{N}(\mathcal{F}_2, \|\cdot\|_{\sup}, \delta) &\leq \mathcal{N}(d(\lambda)^{-1/2}\lambda^{-1/2}\widetilde{\mathcal{F}}_2, \|\cdot\|_{\sup}, \delta) \\ &= \mathcal{N}(\widetilde{\mathcal{F}}_2, \|\cdot\|_{\sup}, d(\lambda)^{1/2}\lambda^{1/2}\delta). \end{aligned}$$

Hence

$$\begin{aligned} J(\mathcal{F}, \delta) &\leq \int_0^\delta \sqrt{\log \mathcal{N}(\widetilde{\mathcal{F}}_2, \|\cdot\|_{\sup}, d(\lambda)^{1/2}\lambda^{1/2}\varepsilon)}d\varepsilon \\ &\asymp \int_0^\delta \sqrt{r\log(d(\lambda)^{-1/2}\lambda^{-1/2}\varepsilon^{-1})}d\varepsilon \\ &\asymp \sqrt{r}\delta\sqrt{\log(d(\lambda)^{-1/2}\lambda^{-1/2}\delta^{-1})} \end{aligned}$$

Now we are ready to calculate the upper bound for $s$. We plug in $n = N/s$ and $d(\lambda) \asymp r$ into (3.12) and (3.13), and by the condition $\lambda = o\Big(\frac{1}{\sqrt{N}}\Big)$, we get $s = o\Big(\frac{N}{\sqrt{\log \lambda^{-1}} \log^6 N}\Big)$. This upper bound needs to allow the case that $s = 1$, which yields the lower bound for $\lambda$: $\sqrt{\log(\lambda^{-1})} = o\big(N \log^{-6} N\big)$. □

### C.2. Proof of Corollary 4.2.

PROOF. Recall that we have $d(\lambda) \asymp r$. To optimize the rate, we choose $\lambda$ such that $d(\lambda)/N \asymp \lambda$, which yields $\lambda = r/N$. By Theorem 3.2 we have

$$\mathbb{E}\big[\|\bar{f}_{N,\lambda} - f_0\|^2_{L_2(\mathbb{P}_Z)}\big] \le Cr/N + s^{-1} a(n, \lambda, J).$$

For the remainder term to be small, we need $s^{-1} a(n, \lambda, J) \lesssim N^{-1}$. Plugging in $a(n, \lambda, J), d(\lambda)$ and $\lambda$, we get the upper bound for $s$. □

### C.3. A Lemma for Exponentially Decaying RKHS.

**Lemma C.1.** Let $d(\lambda) = (-\log \lambda)^{1/p}$. For all $t > 0$, $p \ge 1$ and some positive constants $c, \alpha$, we have

$$\lim_{\lambda \to 0} \frac{1}{d(\lambda)} \sum_{\ell=1}^{\infty} \frac{1}{(1 + \lambda c \exp(\alpha \ell^p))^t} = \alpha^{-1/p}.$$

PROOF. We have by convexity that

$$\sum_{\ell=1}^{\infty} \frac{1}{(1 + \lambda c \exp(\alpha \ell^p))^t} \le \int_0^{\infty} \frac{1}{(1 + \lambda c \exp(\alpha x^p))^t} dx.$$

We then approximate the integral by

$$\begin{aligned}
&\int_0^{\infty} \frac{dx}{(1 + \lambda c \exp(\alpha x^p))^t} \\
&= \int_0^{(\alpha^{-1}\log(1/\lambda))^{1/p}} \frac{dx}{(1 + \lambda c \exp(\alpha x^p))^t} \\
&\qquad + \int_{(\alpha^{-1}\log(1/\lambda))^{1/p}}^{\infty} \frac{dx}{(1 + \lambda c \exp(\alpha x^p))^t} \\
&\le (\alpha^{-1}\log(1/\lambda))^{1/p} + \int_{(\alpha^{-1}\log(1/\lambda))^{1/p}}^{\infty} (c\lambda)^{-t} \exp(-t\alpha x^p) dx \\
&= (\alpha^{-1}\log(1/\lambda))^{1/p} + o(1), \qquad \text{(C.1)}
\end{aligned}$$

where the last equality is by L'Hospital's Rule for $\lambda \to 0$.

Moreover, we have for any $\epsilon \in (0,1)$ that

$$\begin{aligned}\sum_{\ell=1}^{\infty} \frac{1}{(1+\lambda c \exp(\alpha \ell^p))^t} &\geq \int_1^{\infty} \frac{1}{(1+\lambda c \exp(\alpha x^p))^t} dx \\ &\geq \int_1^{(\epsilon\alpha^{-1}\log(1/\lambda))^{1/p}} \frac{1}{(1+\lambda c \exp(\alpha x^p))^t} dx \\ &\geq \frac{1}{(1+c\lambda^{1-\epsilon})^t}\Big((\epsilon\alpha^{-1}\log(1/\lambda))^{1/p}-1\Big) \\ &= \frac{1}{(1+c\lambda^{1-\epsilon})^t}(\epsilon\alpha^{-1}\log(1/\lambda))^{1/p}+O(1). \qquad \text{(C.2)}\end{aligned}$$

Combining (C.1) and (C.2), we get

$$\left(\frac{\epsilon}{\alpha}\right)^{1/p} \leq \lim_{\lambda\to 0}\frac{1}{d(\lambda)}\sum_{\ell=1}^{\infty}\frac{1}{(1+\lambda c\exp(\alpha\ell^p))^t} \leq \left(\frac{1}{\alpha}\right)^{1/p}.$$

for any $\epsilon \in (0,1)$. Lastly, letting $\epsilon \to 1$, we get the desired result. □

### C.4. Proof of Corollary 4.3.

Proof. As before, we start by calculating $d(\lambda)$. By Lemma C.1 with $t=1$, we have $d(\lambda) \asymp (-\log\lambda)^{1/p}$. As $d(\lambda)\to\infty$, Theorem 3.6 shows that $\boldsymbol{\alpha}_{z_0} = \boldsymbol{\gamma}_{z_0} = 0$. Moreover, by (A.25), it holds that

$$\begin{aligned}|W_\lambda f_0(z_0)| &= \lambda\left|\sum_{\ell=1}^{\infty}\frac{\theta_\ell}{\lambda+\mu_\ell}\phi_\ell(z_0)\right| \\ &\leq \lambda\sum_{\ell=0}^{\infty}|\phi_\ell(z_0)\langle f_0,\phi_\ell\rangle_{\mathcal{H}}| = O(\lambda).\end{aligned}$$

Therefore, by Theorem 3.6, we can completely remove the asymptotic bias by choosing $\lambda = o(\sqrt{d(\lambda)/N}\wedge n^{-1/2}) = o\big(N^{-1/2}\log^{1/(2p)}N\wedge n^{-1/2}\big)$. We next calculate the entropy integral. We have that for RKHS with exponentially decaying eigenvalues, by Proposition 17 in Williamson et al. (2001) with $p=2$,

$$\log\mathcal{N}(\widetilde{\mathcal{F}}_2,\|\cdot\|_{\sup},\delta) \asymp \left(\log\frac{1}{\delta}\right)^{\frac{p+1}{p}}.$$

Then following the deduction in the proof of Corollary 4.1, we have

$$\begin{aligned} J(\delta) &\leq \int_0^\delta \sqrt{\log\big(1+\mathcal{N}(\widetilde{\mathcal{F}}_2, \|\cdot\|_{\sup}, c_r d(\lambda)^{1/2}\lambda^{1/2}\varepsilon)\big)}d\varepsilon \\ &\asymp \int_0^\delta \sqrt{\Big(\log\frac{1}{c_r d(\lambda)^{1/2}\lambda^{1/2}\varepsilon}\Big)^{\frac{p+1}{p}}}d\varepsilon \\ &\asymp \delta \log^{\frac{p+1}{2p}}(d(\lambda)^{-1/2}\lambda^{-1/2}\delta^{-1}). \end{aligned}$$

For the range on $s$, we plug in $n = N/s$ and $d(\lambda) \asymp (-\log\lambda)^{1/p}$ into (3.12) and (3.13), and we get that it suffices to take

$$s = o\Big(\frac{N}{\log^6 N \log^{(p+4)/p}(\lambda^{-1})}\Big).$$

Again the upper bound must allow the case that $s = 1$, which yields the lower bound for the choice of $\lambda$. □

### C.5. Proof of Corollary 4.4.

Proof. Recall that we have $d(\lambda) \asymp (-\log\lambda)^{1/p}$. To balance variance and bias, we choose $\lambda = \frac{(\log N)^{1/p}}{N}$. By Theorem 3.2 we have

$$\mathbb{E}\big[\|\bar{f}_{N,\lambda} - f_0\|^2_{L_2(\mathbb{P}_Z)}\big] \leq C(\log N)^{1/p}/N + s^{-1}a(n,\lambda,J).$$

For the remainder term to be small, we need $s^{-1}a(n,\lambda,J) \lesssim (\log N)^{1/p}/N$. Plugging in $d(\lambda)$, $\lambda$ and $J(\mathcal{F},1)$, we get the upper bound for $s$. □

### C.6. Proof of Corollary 4.5.

Proof. Again, we begin by calculating $d(\lambda)$. As $\mu_j \leq cj^{-2\nu}$, we approximate $d(\lambda)$ using integration. For simplicity, let $c = 1$ here. We have

$$\begin{aligned} d(\lambda) \leq \int_0^\infty \frac{1}{1+\lambda x^{2\nu}}dx &= \int_0^{\lambda^{-\frac{1}{2\nu}}} \frac{1}{1+\lambda x^{2\nu}}dx + \int_{\lambda^{-\frac{1}{2\nu}}}^\infty \frac{1}{1+\lambda x^{2\nu}}dx \\ &\leq \Big(1+\frac{1}{1-2\nu}\Big)\lambda^{-\frac{1}{2\nu}}. \end{aligned}$$

On the other hand, we also have

$$\begin{aligned} d(\lambda) \geq \int_1^\infty \frac{1}{1+\lambda x^{2\nu}}dx &= \int_1^{\lambda^{-\frac{1}{2\nu}}} \frac{1}{1+\lambda x^{2\nu}}dx + \int_{\lambda^{-\frac{1}{2\nu}}}^\infty \frac{1}{1+\lambda x^{2\nu}}dx \\ &\geq \Big(2+\frac{1}{2-4\nu}\Big)\lambda^{-\frac{1}{2\nu}}. \end{aligned}$$

Hence we conclude that $d(\lambda) \asymp \lambda^{-\frac{1}{2\nu}}$.

As $d(\lambda) \to \infty$, Theorem 3.6 shows that $\boldsymbol{\alpha}_{z_0} = \boldsymbol{\gamma}_{z_0} = 0$. Similar to proof of Corollary 4.3, we get $|W_\lambda f_0(z_0)| = o(\lambda)$, and by Theorem 3.6, we can remove the asymptotic bias by choosing $\lambda = o\big(\sqrt{d(\lambda)/N} \wedge n^{-1/2}\big) = o\big(N^{-\frac{2\nu}{4\nu+1}} \wedge n^{-1/2}\big)$. We next calculate the entropy integral. We have that for RKHS with polynomially decaying eigenvalues, by Proposition 16 in Williamson et al. (2001),

$$\log \mathcal{N}(\widetilde{\mathcal{F}}_2, \|\cdot\|_{\sup}, \delta) \asymp \left(\frac{1}{\delta}\right)^{\frac{1}{\nu}}.$$

Then following the deduction in the proof of Corollary 4.1, we have

$$\begin{aligned} J(\mathcal{F}, \delta) &\leq \int_0^\delta \sqrt{\log \mathcal{N}(\widetilde{\mathcal{F}}_2, \|\cdot\|_{\sup}, d(\lambda)^{1/2}\lambda^{1/2}\varepsilon)}\, d\varepsilon \\ &\asymp \int_0^\delta \sqrt{\Big(\frac{1}{(d(\lambda)\lambda)^{1/2}\varepsilon}\Big)^{\frac{1}{\nu}}}\, d\varepsilon \\ &\asymp (d(\lambda)\lambda)^{-\frac{1}{4\nu}} \delta^{1-\frac{1}{2\nu}}. \end{aligned}$$

For the range on $s$, we plug in $n = N/s$ and $d(\lambda) \asymp \lambda^{-\frac{1}{2\nu}}$ into (3.12) and (3.13), and it follows that $s$ needs to satisfy

$$s = o\Big(\lambda^{\frac{10\nu-1}{4\nu^2}} N \log^{-6} N\Big).$$

Again the upper bound must allow the case that $s = 1$, which yields the lower bound for the choice of $\lambda$: $\lambda^{-1} = o\big(N^{\frac{4\nu^2}{10\nu-1}}\big)$. □

## C.7. Proof of Corollary 4.6 .

Proof. Recall that we have $d(\lambda) \asymp \lambda^{-1/2\nu}$. To optimize the rate, we choose $\lambda$ such that $d(\lambda)/N \asymp \lambda$, which yields $\lambda = N^{-\frac{2\nu}{2\nu+1}}$. By Theorem 3.2 we have

$$\mathbb{E}\big[\|\bar{f}_{N,\lambda} - f_0\|^2_{L_2(\mathbb{P}_Z)}\big] \leq C N^{-\frac{2\nu}{2\nu+1}} + s^{-1} a(n, \lambda, J)$$

For the remainder term to be small, we need $s^{-1}a(n,\lambda,J) \lesssim N^{-\frac{2\nu}{2\nu+1}}$. Plugging in $a(n,\lambda,J), d(\lambda)$ and $\lambda$, we get the upper bound for $s$. □

**C.8. Proof of Lemma 4.1.** Recall that $\lambda = d(\lambda)^{-2\nu}$. By Theorem 3.4 for the asymptotic variance, we compute that

$$
\begin{aligned}
d(\lambda)^{-1}\|\widetilde{K}_{z_0}\|^2_{L_2(\mathbb{P}_Z)} &= d(\lambda)^{-1}\sum_{\ell=1}^{\infty}\Big(\frac{\phi_\ell(z_0)}{1+\lambda/\mu_\ell}\Big)^2 \\
&= d(\lambda)^{-1}\Big(1+\sum_{\ell=1}^{\infty}\frac{2cos^2(2\ell\pi z_0)+2sin^2(2\ell\pi z_0)}{(1+\lambda(2\ell\pi)^{2\nu})^2}\Big) \\
&= d(\lambda)^{-1}\Big(1+\sum_{\ell=1}^{\infty}\frac{2}{(1+(2\ell\pi/d(\lambda))^{2\nu})^2}\Big)
\end{aligned}
$$

And we have that

$$
\begin{aligned}
\sum_{\ell=0}^{\infty}\frac{2\pi/d(\lambda)}{(1+(2\ell\pi/d(\lambda))^{2\nu})^2} &\le \sum_{\ell=1}^{\infty}\int_{2\pi(\ell-1)/d(\lambda)}^{2\pi\ell/d(\lambda)}\frac{1}{(1+x^{2\nu})^2}dx \\
&\to \int_0^{\infty}\frac{1}{(1+x^{2\nu})^2}dx
\end{aligned}
$$

and similarly

$$
\begin{aligned}
\sum_{\ell=1}^{\infty}\frac{2\pi/d(\lambda)}{(1+(2\ell\pi/d(\lambda))^{2\nu})^2} &\ge \sum_{\ell=1}^{\infty}\int_{2\pi(\ell-1)/d(\lambda)}^{2\pi\ell/d(\lambda)}\frac{1}{(1+x^{2\nu})^2}dx \\
&\to \int_0^{\infty}\frac{1}{(1+x^{2\nu})^2}dx
\end{aligned}
$$

The two inequalities yield

$$
d(\lambda)^{-1}\|\widetilde{K}_{z_0}\|^2_{L_2(\mathbb{P}_Z)} \to \int_0^{\infty}\frac{1}{\pi(1+x^{2\nu})^2}dx
$$

and so

$$
\sigma^2_{z_0} = \int_0^{\infty}\frac{1}{\pi(1+x^{2\nu})^2}dx. \tag{C.3}
$$

## APPENDIX D: PROOF OF RESULTS IN SECTION 5

**D.1. Proof of Proposition 5.1.**

Proof. Recall (7.2) from proof of Theorem 3.2 in Section 7.3

$$
\widehat{m}^{(j)} - m_0 = \frac{1}{n}\sum_{i\in S_j} R_{U_i}\varepsilon_i - P_\lambda m_0 - Rem^{(j)}.
$$

Also recall $m_0^* = m_0 - P_\lambda m_0$. Taking average of the above equation for all $j$

over $s$, we have

$$\bar{m} - m_0^* = \frac{1}{N}\sum_{i=1}^{N} \varepsilon_i R_{U_i} + \frac{1}{s}\sum_{j=1}^{s} Rem^{(j)}, \tag{D.1}$$

which decomposes into

$$\bar{\boldsymbol{\beta}} - \boldsymbol{\beta}_0^* = \frac{1}{N}\sum_{i=1}^{N} L_{U_i}\varepsilon_i - \frac{1}{s}\sum_{j=1}^{s} Rem_{\beta}^{(j)}, \tag{D.2}$$

and

$$\bar{f} - f_0^* = \frac{1}{N}\sum_{i=1}^{N} N_{U_i}\varepsilon_i - \frac{1}{s}\sum_{j=1}^{s} Rem_{f}^{(j)}. \tag{D.3}$$

Similar to proof in Section 7.3, we can show that the first term weakly converges to a normal distribution, and the remainder term is asymptotically ignorable. Recall the definition of $m_0^* = (id - P_\lambda)f_0$.

Therefore, we deduct that

$$\begin{aligned}
&\big(\boldsymbol{x}^T, 1\big)\begin{pmatrix} \sqrt{N}(\bar{\boldsymbol{\beta}} - \boldsymbol{\beta}_0^*) \\ \sqrt{N/d(\lambda)}\big(\bar{f}(z_0) - f_0^*(z_0)\big) \end{pmatrix} \\
&= \sqrt{N}\boldsymbol{x}^T(\bar{\boldsymbol{\beta}} - \boldsymbol{\beta}_0^*) + (N/d(\lambda))^{\frac{1}{2}}(\bar{f}(z_0) - f_0^*(z_0)) \\
&\leq \frac{1}{\sqrt{N}}\sum_{i=1}^{N}(\varepsilon_i \boldsymbol{x}^T L_{U_i} + d(\lambda)^{-1/2}\varepsilon_i N_{U_i}) \\
&\qquad \frac{1}{s}\sum_{j=1}^{s}\sqrt{N}\boldsymbol{x}^T Rem_{\beta}^{(j)} - \frac{1}{s}\sum_{j=1}^{s}\sqrt{N/d(\lambda)}Rem_f^{(j)}(z_0).
\end{aligned} \tag{D.4}$$

We can show that the first term is asymptotic normal by central limit theorem: first note that the summands are i.i.d. and with mean zero. Moreover, by Proposition 2.3,

$$\begin{aligned}
\boldsymbol{x}^T L_U + d(\lambda)^{-1/2} N_U(z_0) &= \boldsymbol{x}^T L_U + d(\lambda)^{-1/2}(\widetilde{K}_Z(z_0) - \boldsymbol{A}(z_0)^T L_U) \\
&= (\boldsymbol{x} - d(\lambda)^{-1/2}\boldsymbol{A}(z_0))^T L_U + d(\lambda)^{-1/2}\widetilde{K}_Z(z_0)
\end{aligned}$$

We compute

$$\mathbb{E}\big[(\boldsymbol{x}^T L_{U_i} + d(\lambda)^{-1/2} N_{U_i}(z_0))^2\big]$$

(D.5)
$$= \underbrace{\mathbb{E}\big[(\boldsymbol{x} - d(\lambda)^{-1/2}\boldsymbol{A}(z_0))^T L_U L_U^T (\boldsymbol{x} - d(\lambda)^{-1/2}\boldsymbol{A}(z_0))\big]}_{(I)} + \underbrace{\mathbb{E}\big[d(\lambda)^{-1}\widetilde{K}_Z(z_0)^2\big]}_{(II)} + \underbrace{\mathbb{E}\big[2d(\lambda)^{-1/2}\widetilde{K}_Z(z_0)(\boldsymbol{x} - d(\lambda)^{-1/2}\boldsymbol{A}(z_0))^T L_U\big]}_{(III)}.$$

For the second term in (D.5), we have

$$(II) = \mathbb{E}\big[d(\lambda)^{-1}\widetilde{K}_Z(z_0)^2\big] = d(\lambda)^{-1}\|\widetilde{K}_{z_0}\|^2_{L_2(\mathbb{P}_Z)} \to \sigma^2_{z_0}. \tag{D.6}$$

For the first term in (D.5), recall from the proof of Theorem 3.5 that $\mathbb{E}[L_U L_U^T] \to \boldsymbol{\Omega}^{-1}$. Also, $\boldsymbol{A}(z_0)/d(\lambda)^{1/2} \to -\boldsymbol{\gamma}_{z_0}$. Hence we have

$$(I) \to (\boldsymbol{x} + \boldsymbol{\gamma}_{z_0})^T \boldsymbol{\Omega}^{-1}(\boldsymbol{x} + \boldsymbol{\gamma}_{z_0}). \tag{D.7}$$

Moreover, recall from Section 7.3 that

$$d(\lambda)^{-1/2}\mathbb{E}\big[\widetilde{K}_Z(z_0)L_U\big] \to \boldsymbol{\Omega}^{-1}\boldsymbol{\alpha}_{z_0}.$$

Therefore it follows that

$$(III) \to 2(\boldsymbol{x} + \boldsymbol{\gamma}_{z_0})^T \boldsymbol{\Omega}^{-1}\boldsymbol{\alpha}_{z_0}. \tag{D.8}$$

Hence combining (D.7), (D.6) and (D.8), the limit of (D.5) is

$$\mathbb{E}\big[(\boldsymbol{x}^T L_{U_i} + d(\lambda)^{-1/2} N_{U_i}(z_0))^2\big] \to \boldsymbol{x}^T \boldsymbol{\Omega}^{-1}\boldsymbol{x} + 2\boldsymbol{x}^T \boldsymbol{\Sigma}_{12} + \Sigma_{22},$$

for any $\boldsymbol{x} \in \mathbb{R}^p$. Therefore the limit distribution follows by central limit theorem. Now for the remainder terms, by Lemma 7.2, if Condition (3.12) is satisfied, we have

$$\begin{aligned}|\frac{1}{s}\sum_{j=1}^s \sqrt{N}\boldsymbol{x}^T Rem_\beta^{(j)}| &\le C\sqrt{N}\|\frac{1}{s}\sum_{j=1}^s Rem_\beta^{(j)}\|_2 \\ &= o_P(N^{1/2}s^{-1/2}b(n,\lambda,J)\log N).\end{aligned}$$

and

$$\begin{aligned}|\frac{1}{s}\sum_{j=1}^s \sqrt{N/d(\lambda)}\boldsymbol{x}^T Rem_f^{(j)}(z_0)| &\le C\sqrt{N/d(\lambda)}\|\frac{1}{s}\sum_{j=1}^s Rem_f^{(j)}\|_{\sup} \\ &\le C' N^{1/2}\|\frac{1}{s}\sum_{j=1}^s Rem_f^{(j)}\|_{\mathcal{C}} \\ &= o_P(N^{1/2}s^{-1/2}b(n,\lambda,J)\log N).\end{aligned}$$

where in the second inequality we used Lemma A.1. Then if Condition (3.13)

is satisfied, we have $N^{1/2}s^{-1/2}b(n,\lambda,J)\log N \to 0$. Hence by (D.4), it follows that

$$\big(\boldsymbol{x}^T, 1\big)\left(\begin{array}{c}\sqrt{N}(\bar{\boldsymbol{\beta}} - \boldsymbol{\beta}_0^*)\\ \sqrt{N/d(\lambda)}\big(\bar{f}(z_0) - f_0^*(z_0)\big)\end{array}\right) \to N\big(0, \sigma^2(\boldsymbol{x}^T\boldsymbol{\Omega}^{-1}\boldsymbol{x} + 2\boldsymbol{x}^T\boldsymbol{\Sigma}_{12} + \Sigma_{22})\big).$$

Hence the conclusion follows by the arbitrariness of $\boldsymbol{x}$ using Wold device. □

## APPENDIX E: PROOFS OF LEMMAS IN SECTION 7

### E.1. Proof of Lemma 7.1.

PROOF. Recall from Section 7.4 that

$$Z_n(\widetilde{m}) = \frac{1}{2}d(\lambda)^{-1/2}n^{1/2}q_{n,\lambda}^{-1}Rem^{(j)} = \frac{1}{2}c_r^{-1}d(\lambda)^{-1}n^{1/2}r_{n,\lambda}^{-1}Rem^{(j)}$$

We showed in Section 7.4 that $Z_n(m)$ is a sub-Gaussian process. Letting $\mathbb{U}^{(j)} = (\mathbb{X}^{(j)}, \mathbb{Z}^{(j)})$, where $\mathbb{X}^{(j)}$ and $\mathbb{Z}^{(j)}$ are designs on $j$-th sub-population. Without causing any confusion, we can remove the the superscript $(j)$. We have

$$\begin{aligned}\mathbb{E}[\|Rem^{(j)}\|_{\mathcal{A}}^2] &= \mathbb{E}\big[\mathbb{E}[\|Rem^{(j)}\|_{\mathcal{A}}^2 \,|\, \mathbb{U}]\big]\\ &= \mathbb{E}\big[\mathbb{E}[\|Rem^{(j)}\|_{\mathcal{A}}^2 \,|\, \mathbb{U}]I_{\mathcal{E}}\big] + \mathbb{E}\big[\mathbb{E}[\|Rem^{(j)}\|_{\mathcal{A}}^2 \,|\, \mathbb{U}]I_{\mathcal{E}^c}\big],\end{aligned} \tag{E.1}$$

where $\mathcal{E}$ is the event defined in Section 7.4. For the first term in (E.1), we have

$$\begin{aligned}&\mathbb{E}\big[\mathbb{E}[\|Rem^{(j)}\|_{\mathcal{A}}^2 \,|\, \mathbb{U}]I_{\mathcal{E}}\big]\\ &\quad= 4c_r^2 d(\lambda)^2 n^{-1} r_{n,\lambda}^2 \mathbb{E}\big[\mathbb{E}[\|Z_n(\widetilde{m})\|_{\mathcal{A}}^2 \,|\, \mathbb{U}]I_{\mathcal{E}}\big]\\ &\quad= 4c_r^2 d(\lambda)^2 n^{-1} r_{n,\lambda}^2 \int_0^\infty \mathbb{P}\big(\mathbb{E}[\|Z_n(\widetilde{m})\|_{\mathcal{A}}^2 \,|\, \mathbb{U}]I_{\mathcal{E}} \ge x\big)dx\\ &\quad= 4c_r^2 d(\lambda)^2 n^{-1} r_{n,\lambda}^2 \Big\{\int_0^{J(\mathcal{F},1)^2} \mathbb{P}\big(\mathbb{E}[\|Z_n(\widetilde{m})\|_{\mathcal{A}}^2 \,|\, \mathbb{U}]I_{\mathcal{E}} \ge x\big)dx\\ &\qquad\qquad + \int_{J(\mathcal{F},1)^2}^\infty \mathbb{P}\big(\mathbb{E}[\|Z_n(\widetilde{m})\|_{\mathcal{A}}^2 \,|\, \mathbb{U}]I_{\mathcal{E}} \ge x\big)dx\Big\}\\ &\quad\le 4c_r^2 d(\lambda)^2 n^{-1} r_{n,\lambda}^2 \Big\{J(\mathcal{F},1)^2 + \int_0^\infty \mathbb{P}\big(\mathbb{E}[\|Z_n(\widetilde{m})\|_{\mathcal{A}}^2 \,|\, \mathbb{U}]I_{\mathcal{E}} \ge x + J(\mathcal{F},1)^2\big)dx\Big\}.\end{aligned}$$

In Section 7.4 we proved that $\mathcal{E} \subset \{\widetilde{m} \in \mathcal{F}\}$. Therefore we have

$$\mathbb{E}[\|Z_n(\widetilde{m})\|_{\mathcal{A}}^2 \,|\, \mathbb{U}]I_{\mathcal{E}} \le \sup_{m\in\mathcal{F}} \|Z_n(m)\|_{\mathcal{A}}^2.$$

and by Lemma F.1 and the fact that $\mathrm{diam}(\mathcal{F}) \le 1$, we have

$$\mathbb{P}\Big(\sup_{m\in\mathcal{F}} \|Z_n(m)\|_{\mathcal{A}}^2 \ge x + J(\mathcal{F},1)^2\Big) \le \mathbb{P}\Big(\sup_{m\in\mathcal{F}} \|Z_n(m)\|_{\mathcal{A}} \ge (\sqrt{x} + J(\mathcal{F},1))/2\Big) \le C\exp\big(-x/C\big).$$

Hence it follows that

$$\begin{aligned}\mathbb{E}\big[\mathbb{E}[\|Rem^{(j)}\|_{\mathcal{A}}^2 \,|\, \mathbb{U}]I_{\mathcal{E}}\big] &\le 2c_r^2 d(\lambda)^2 n^{-1} r_{n,\lambda}^2\Big(J(\mathcal{F},1)^2 + \int_0^\infty C\exp(-x/C)dx\Big)\\ &= 2c_r^2 d(\lambda)^2 n^{-1} r_{n,\lambda}^2\big(J(\mathcal{F},1)^2 + C^2\big). \qquad \text{(E.2)}\end{aligned}$$

We now turn to control the second term in (E.1). By Lemma G.2, we get that $\mathbb{E}\big[\|\Delta f^{(j)}\|_{\mathcal{H}}^2 \,|\, \mathbb{U}\big] \le 2\sigma^2/\lambda + 4\|f_0\|_{\mathcal{H}}^2$. Also by first order optimality condition with respect to $\widehat{\boldsymbol\beta}$,

$$\widehat{\boldsymbol\beta}^{(j)} - \boldsymbol\beta_0^{(j)} = (\mathbb{X}^T\mathbb{X})^{-1}\mathbb{X}^T(f_0(\mathbb{Z}) - \widehat{f}^{(j)}(\mathbb{Z}) + \boldsymbol\varepsilon^{(j)}),$$

where we omitted the superscript of $(j)$ for the designs $\mathbb{X}^{(j)}$ and $\mathbb{Z}^{(j)}$. Hence

$$\|\widehat{\boldsymbol\beta} - \boldsymbol\beta\|_2^2 \le 2\|(\mathbb{X}^T\mathbb{X})^{-1}\mathbb{X}^T\boldsymbol\varepsilon^{(j)}\|_2^2 + 2\|(\mathbb{X}^T\mathbb{X})^{-1}\mathbb{X}^T(f_0(\mathbb{Z}) - \widehat{f}^{(j)}(\mathbb{Z}))\|_2^2,$$

Taking conditional expectation yields

$$\begin{aligned}\mathbb{E}[\|\widehat{\boldsymbol\beta}^{(j)} - \boldsymbol\beta_0^{(j)}\|_2^2 \,|\, \mathbb{U}] &\le 2\mathbb{E}\big[\|(\mathbb{X}^T\mathbb{X})^{-1}\mathbb{X}^T\boldsymbol\varepsilon^{(j)}\|_2^2 \,|\, \mathbb{U}\big]\\ &\quad + 2\mathbb{E}\big[\|(\mathbb{X}^T\mathbb{X})^{-1}\mathbb{X}^T(f_0(\mathbb{Z}) - \widehat{f}^{(j)}(\mathbb{Z}))\|_2^2 \,|\, \mathbb{U}\big]. \qquad \text{(E.3)}\end{aligned}$$

Denote $\widehat{\boldsymbol\Sigma}^{(j)} = \mathbb{X}^T\mathbb{X}$, we first control the first term (E.3). Note that

$$(\mathbb{X}^T\mathbb{X})^{-1}\mathbb{X}^T\boldsymbol\varepsilon^{(j)} = \frac{1}{n}\sum_{i\in S_j}(\widehat{\boldsymbol\Sigma}^{(j)})^{-1}\boldsymbol X_i\varepsilon_i.$$

Taking conditional expectation and by independence, we have

$$\begin{aligned}\mathbb{E}\big[\|(\mathbb{X}^T\mathbb{X})^{-1}\mathbb{X}^T\boldsymbol\varepsilon^{(j)}\|_2^2 \,|\, \mathbb{U}\big] &= \mathbb{E}\big[\|\frac{1}{n}\sum_{i\in S_j}(\widehat{\boldsymbol\Sigma}^{(j)})^{-1}\boldsymbol X_i\varepsilon_i\|_2^2 \,|\, \mathbb{U}\big]\\ &= \frac{1}{n^2}\sum_{i\in S_j}\mathbb{E}\Big[(\widehat{\boldsymbol\Sigma}^{(j)})^{-1}\boldsymbol X_i\varepsilon_i\|_2^2 \,|\, \mathbb{U}\Big]\\ &= \frac{1}{n^2}\sum_{i\in S_j}\sigma^2\|(\widehat{\boldsymbol\Sigma}^{(j)})^{-1}\boldsymbol X_i\|_2^2. \qquad \text{(E.4)}\end{aligned}$$

For the second term in (E.3), we have similar to above that

$$(\mathbb{X}^T\mathbb{X})^{-1}\mathbb{X}^T(f_0(\mathbb{Z}) - \widehat{f}^{(j)}(\mathbb{Z})) = \frac{1}{n}\sum_{i\in S_j}(\widehat{\boldsymbol\Sigma}^{(j)})^{-1}\boldsymbol X_i(f_0(Z_i) - \widehat{f}^{(j)}(Z_i)).$$

Taking conditional expectation, we have

$$\begin{aligned}\mathbb{E}\big[\|(\mathbb{X}^T\mathbb{X})^{-1}\mathbb{X}^T(f_0(\mathbb{Z})-\widehat{f}^{(j)}(\mathbb{Z}))\|_2^2\,|\,\mathbb{U}\big]\\ &=\mathbb{E}\big[\|\frac{1}{n}\sum_{i\in S_j}(\widehat{\boldsymbol{\Sigma}}^{(j)})^{-1}\boldsymbol{X}_i(f_0(Z_i)-\widehat{f}^{(j)}(Z_i))\|_2^2\,|\,\mathbb{U}\big]\\ &\le\frac{1}{n}\sum_{i\in S_j}\mathbb{E}\big[\|(\widehat{\boldsymbol{\Sigma}}^{(j)})^{-1}\boldsymbol{X}_i(f_0(Z_i)-\widehat{f}^{(j)}(Z_i))\|_2^2\,|\,\mathbb{U}\big],\end{aligned}$$

where the inequality is by $(\sum_{i=1}^n a_i)^2\le n\sum_{i=1}^n a_i^2$. Hence we have

$$\begin{aligned}\mathbb{E}\big[\|(\mathbb{X}^T\mathbb{X})^{-1}\mathbb{X}^T(f_0(\mathbb{Z})-\widehat{f}^{(j)}(\mathbb{Z}))\|_2^2\,|\,\mathbb{U}\big]\\ &\le\frac{1}{n}\sum_{i\in S_j}\mathbb{E}\big[(f_0(Z_i)-\widehat{f}^{(j)}(Z_i))^2\,|\,\mathbb{U}\big]\|(\widehat{\boldsymbol{\Sigma}}^{(j)})^{-1}\boldsymbol{X}_i\|_2^2\\ &\le\frac{1}{n}\sum_{i\in S_j}\mathbb{E}\big[\|f_0-\widehat{f}^{(j)}\|_{\sup}^2\,|\,\mathbb{U}\big]\|(\widehat{\boldsymbol{\Sigma}}^{(j)})^{-1}\boldsymbol{X}_i\|_2^2\\ &\le\frac{1}{n}\sum_{i\in S_j}c_k\mathbb{E}\big[\|f_0-\widehat{f}^{(j)}\|_{\mathcal{H}}^2\,|\,\mathbb{U}\big]\|(\widehat{\boldsymbol{\Sigma}}^{(j)})^{-1}\boldsymbol{X}_i\|_2^2\\ &\le\frac{1}{n}\sum_{i\in S_j}c_k(2\sigma^2/\lambda+4\|f_0\|_{\mathcal{H}}^2)\|(\widehat{\boldsymbol{\Sigma}}^{(j)})^{-1}\boldsymbol{X}_i\|_2^2,\end{aligned}\tag{E.5}$$

where $c_k=\sup_z K(z,z)$. The second last inequality follows from the fact that $\|f\|_{\sup}\le\sup_z\|K_z\|_{\mathcal{H}}\|f\|_{\mathcal{H}}=c_k^{1/2}\|f\|_{\mathcal{H}}$ and the last inequality is by Lemma G.2. Combing (E.4) and (E.5), we have by (E.3) that

$$\mathbb{E}[\|\widehat{\boldsymbol{\beta}}^{(j)}-\boldsymbol{\beta}_0^{(j)}\|_2^2\,|\,\mathbb{U}]\le C\lambda^{-1}\frac{1}{n}\sum_{i\in S_j}\|(\widehat{\boldsymbol{\Sigma}}^{(j)})^{-1}\boldsymbol{X}_i\|_2^2$$

holds almost surely for some constant $C$ As we have

$$\|\widehat{\boldsymbol{\beta}}^{(j)}-\boldsymbol{\beta}_0^{(j)}\|_{L_2(\mathbb{P}_X)}^2=(\widehat{\boldsymbol{\beta}}^{(j)}-\boldsymbol{\beta}_0^{(j)})^T\boldsymbol{\Sigma}(\widehat{\boldsymbol{\beta}}^{(j)}-\boldsymbol{\beta}_0^{(j)})\le\|\boldsymbol{\Sigma}^{-1/2}\|\|\widehat{\boldsymbol{\beta}}^{(j)}-\boldsymbol{\beta}_0^{(j)}\|_2^2,$$

it follows that

$$\mathbb{E}[\|\widehat{\boldsymbol{\beta}}^{(j)}-\boldsymbol{\beta}_0^{(j)}\|_{L_2(\mathbb{P}_X)}^2\,|\,\mathbb{U}]\le C\lambda^{-1}\frac{1}{n}\sum_{i\in S_j}\|(\widehat{\boldsymbol{\Sigma}}^{(j)})^{-1}\boldsymbol{X}_i\|_2^2,\tag{E.6}$$

for a constant $C$ that is different from above. Lastly, we have $\|\widehat{f}^{(j)}-f_0\|_{L_2(\mathbb{P}_Z)}\le\|\widehat{f}^{(j)}-f_0\|_{\sup}\le c_k^{1/2}\|\widehat{f}^{(j)}-f_0\|_{\mathcal{H}}$, and it follows that

$$\mathbb{E}\big[\|\widehat{f}^{(j)}-f_0\|_{L_2(\mathbb{P}_Z)}^2\,|\,\mathbb{U}\big]\lesssim\lambda^{-1}.\tag{E.7}$$

Note that for any $m=(\boldsymbol{\beta},f)$, $\|m\|_{\mathcal{A}}^2=\|\boldsymbol{X}^T\boldsymbol{\beta}+f(Z)\|_{L_2(\mathbb{P}_U)}^2+\lambda\|f\|_{\mathcal{H}}^2\le$

$2\|\boldsymbol{\beta}\|^2_{L_2(\mathbb{P}_X)} + 2\|f\|^2_{L_2(\mathbb{P}_Z)} + \lambda\|f\|^2_{\mathcal{H}}$. Hence by (E.6), (E.7) and (G.1), we have

$$\mathbb{E}\big[\|\Delta m^{(j)}\|^2_{\mathcal{A}} \,|\, \mathbb{U}\big] \le C\lambda^{-1}n^{-1}\sum_{i\in S_j}\|(\widehat{\boldsymbol{\Sigma}}^{(j)})^{-1}\boldsymbol{X}_i\|_2^2.$$

Moreover,

$$\begin{aligned}\|Rem^{(j)}\|_{\mathcal{A}} &\le \|\frac{1}{n}\sum_{i\in S_j}\Delta m^{(j)}(U_i)R_{U_i}\|_{\mathcal{A}} + \|\mathbb{E}_U[\Delta m^{(j)}(U)R_U]\| \\ &= \frac{1}{n}\sum_{i\in S_j}\|\Delta m^{(j)}(U_i)R_{U_i}\|_{\mathcal{A}} + \mathbb{E}_U\big[\|\Delta m^{(j)}(U)R_U\|\big] \\ &\le 2d(\lambda)\|\Delta m^{(j)}\|_{\mathcal{A}},\end{aligned}$$

where in the last inequality we used $|\Delta m^{(j)}(U)| \le d(\lambda)^{1/2}\|\Delta m^{(j)}\|_{\mathcal{A}}$ and $\|R_U\|_{\mathcal{A}} \le d(\lambda)^{1/2}$. Hence

$$\mathbb{E}\big[\|Rem^{(j)}\|^2_{\mathcal{A}} \,|\, \mathbb{U}\big] \le 4d(\lambda)^2\mathbb{E}\big[\|\Delta m^{(j)}\|^2_{\mathcal{A}} \,|\, \mathbb{U}\big] \le Cd(\lambda)^2\lambda^{-1}n^{-1}\sum_{i\in S_j}\|(\widehat{\boldsymbol{\Sigma}}^{(j)})^{-1}\boldsymbol{X}_i\|_2^2.$$

Hence, we have that the second term in (E.1)

$$\begin{aligned}\mathbb{E}\big[\mathbb{E}[\|Rem^{(j)}\|^2_{\mathcal{A}} \,|\, \mathbb{U}]I_{\mathcal{E}^c}\big] &\le Cd(\lambda)^2\lambda^{-1}n^{-1}\sum_{i\in S_j}\mathbb{E}\big[\|(\widehat{\boldsymbol{\Sigma}}^{(j)})^{-1}\boldsymbol{X}_i\|_2^2 I_{\mathcal{E}^c}\big] \\ &\le Cd(\lambda)^2\lambda^{-1}n^{-1}\mathbb{P}(\mathcal{E}^c)\sum_{i\in S_j}\mathbb{E}\big[\|(\widehat{\boldsymbol{\Sigma}}^{(j)})^{-1}\boldsymbol{X}_i\|_2^4\big] \\ &\le C'd(\lambda)^2\lambda^{-1}\mathbb{P}(\mathcal{E}^c), \end{aligned} \tag{E.8}$$

where the second last inequality is by Holder's inequality and the last one by assumption on the design. By (E.8) and Lemma 7.4, we obtain

$$\mathbb{E}\big[\mathbb{E}[\|Rem^{(j)}\|^2_{\mathcal{A}} \,|\, \mathbb{U}]I_{\mathcal{E}^c}\big] \lesssim d(\lambda)^2\lambda^{-1}n\exp(-c\log^2 N). \tag{E.9}$$

Finally, plugging (E.2) and (E.9) into (E.1), we have for sufficiently large $n$,

$$\mathbb{E}\big[\|Rem^{(j)}\|^2_{\mathcal{A}}\big] \le 2c_r^2 d(\lambda)^2 n^{-1} r_{n,\lambda}^2\big(J(\mathcal{F},1)^2 + C\big) + C'd(\lambda)^2\lambda^{-1}n\exp(-c\log^2 N), \tag{E.10}$$

as desired.

We can apply similar arguments as above to bound $\|Rem_f^{(j)}\|_{\mathcal{C}}$ and $\|1/s\sum_{j=1}^s Rem_f^{(j)}\|_{\mathcal{C}}$, by changing $J(\mathcal{F},1)$ to $J(\mathcal{F}_2,1)$, which is dominated by $J(\mathcal{F},1)$. The bounds of $\|Rem_\beta^{(j)}\|_2$ and $\|1/s\sum_{j=1}^s Rem_\beta^{(j)}\|_2$ then follow from triangular inequality. □

**E.2. Proof of Lemma 7.3.**

PROOF. The main term (I) can be rearranged as follows:

$$(I) = \underbrace{\frac{1}{\sqrt{n}} \sum_{i \in S_j} \varepsilon_i \Big( \boldsymbol{x}^T L_{U_i} + s^{-1/2} d(\lambda)^{-1/2} N_{U_i}(z_0) \Big)}_{(III)} + \underbrace{\frac{1}{\sqrt{N}} \sum_{i \notin S_j} d(\lambda)^{-1/2} N_{U_i}(z_0) \varepsilon_i}_{(IV)}$$

When analyzing (I), we consider two cases: (1) $s \to \infty$ and (2) $s$ is fixed.

Case 1: $s \to \infty$. We first apply CLT to the first component of term (III), i.e., $\frac{1}{\sqrt{n}} \sum_{i \in S_j} \varepsilon_i \boldsymbol{x}^T L_{U_i}$. The summands are i.i.d. with mean zero. Moreover,

$$\mathbb{E}\big[(\varepsilon \boldsymbol{x}^T L_U)^2\big] = \sigma^2 \boldsymbol{x}^T \mathbb{E}[L_U L_U^T] \boldsymbol{x}.$$

By the proof of Theorem 3.5, we have $\mathbb{E}[L_U L_U^T] \to \boldsymbol{\Omega}^{-1}$. Therefore by CLT, we have

$$\frac{1}{\sqrt{n}} \sum_{i \in S_j} \varepsilon_i \boldsymbol{x}^T L_{U_i} \rightsquigarrow N(0, \sigma^2 \boldsymbol{x}^T \boldsymbol{\Omega}^{-1} \boldsymbol{x}). \tag{E.11}$$

We next consider $\frac{1}{\sqrt{n}} \sum_{i \in S_j} s^{-1/2} d(\lambda)^{-1/2} N_{U_i}(z_0) \varepsilon_i$ which is the second component in (III). Again the summands are i.i.d. with mean zero. By Proposition 2.3 we have

$$\begin{aligned}
\text{(E.12)}\quad & \mathbb{E}\big[(d(\lambda)^{-1/2} \varepsilon N_U(z_0))^2\big] \\
&= \sigma^2 d(\lambda)^{-1} \mathbb{E}\big[(\widetilde{K}_Z(z_0) - L_U^T \boldsymbol{A}(z_0))^2\big] \\
&= \sigma^2 d(\lambda)^{-1} \mathbb{E}\big[\widetilde{K}_Z(z_0)^2\big] + \sigma^2 d(\lambda)^{-1} \boldsymbol{A}(z_0)^T \mathbb{E}[L_U L_U^T] \boldsymbol{A}(z_0) \\
&\quad - 2\sigma^2 d(\lambda)^{-1} \mathbb{E}[\widetilde{K}_Z(z_0) L_U^T \boldsymbol{A}(z_0)].
\end{aligned}$$

For the first term in (E.12), by condition in the lemma, we have

$$\sigma^2 d(\lambda)^{-1} \mathbb{E}\big[\widetilde{K}_Z(z_0)^2\big] = \sigma^2 d(\lambda)^{-1} \|\widetilde{K}_{z_0}\|^2_{L_2(\mathbb{P}_Z)} \to \sigma^2 \sigma^2_{z_0}$$

For the second term in (E.12), as $\boldsymbol{A}(z_0)/d(\lambda)^{1/2} \to -\boldsymbol{\gamma}_{z_0}$, and $\mathbb{E}[L_U L_U^T] \to \boldsymbol{\Omega}^{-1}$, we have

$$\sigma^2 d(\lambda)^{-1} \boldsymbol{A}(z_0)^T \mathbb{E}[L_U L_U^T] \boldsymbol{A}(z_0) \to \sigma^2 \boldsymbol{\gamma}_{z_0}^T \boldsymbol{\Omega}^{-1} \boldsymbol{\gamma}_{z_0}.$$

For the last term in (E.12), we consider

$$d(\lambda)^{-1/2} \mathbb{E}\big[\widetilde{K}_Z(z_0) L_U\big] = d(\lambda)^{-1/2} (\boldsymbol{\Omega} + \boldsymbol{\Sigma}_\lambda)^{-1} \mathbb{E}[\widetilde{K}_Z(z_0)(\boldsymbol{X} - \boldsymbol{A}(Z))].$$

We have $(\boldsymbol{\Omega}+\boldsymbol{\Sigma}_\lambda)^{-1} \to \boldsymbol{\Omega}^{-1}$ and

$$\begin{aligned}
\mathbb{E}[\widetilde{K}_Z(z_0)(\boldsymbol{X}-\boldsymbol{A}(Z))] &= d(\lambda)^{-1/2}(\langle \boldsymbol{B}, \widetilde{K}_{z_0}\rangle_{L_2(\mathbb{P}_Z)} - \langle \boldsymbol{A}, \widetilde{K}_{z_0}\rangle_{L_2(\mathbb{P}_Z)}) \\
&= d(\lambda)^{-1/2}(\langle \boldsymbol{A}, \widetilde{K}_{z_0}\rangle_{\mathcal{C}} - \langle \boldsymbol{A}, \widetilde{K}_{z_0}\rangle_{L_2(\mathbb{P}_Z)}) \\
&= d(\lambda)^{-1/2}\lambda\langle \boldsymbol{A}, \widetilde{K}_{z_0}\rangle_{\mathcal{H}} \\
&= d(\lambda)^{-1/2}\langle W_\lambda \boldsymbol{A}, \widetilde{K}_{z_0}\rangle_{\mathcal{C}} \\
&= d(\lambda)^{-1/2} W_\lambda \boldsymbol{A}(z_0) \to \boldsymbol{\alpha}_{z_0}
\end{aligned}$$

Hence $d(\lambda)^{-1/2}\mathbb{E}\big[\widetilde{K}_Z(z_0)L_U\big] \to \boldsymbol{\Omega}^{-1}\boldsymbol{\alpha}_{z_0}$ and so $d(\lambda)^{-1}\mathbb{E}[\widetilde{K}_Z(z_0)L_U^T\boldsymbol{A}(z_0)] \to \boldsymbol{\gamma}_{z_0}^T\boldsymbol{\Omega}^{-1}\boldsymbol{\alpha}_{z_0}$. In summary, we have

$$\mathbb{E}\big[(d(\lambda)^{-1/2}\varepsilon N_U(z_0))^2\big] \to \sigma^2\big(\sigma_{z_0}^2 + \boldsymbol{\gamma}_{z_0}^T\boldsymbol{\Omega}^{-1}\boldsymbol{\gamma}_{z_0} + 2\boldsymbol{\gamma}_{z_0}^T\boldsymbol{\Omega}^{-1}\boldsymbol{\alpha}_{z_0}\big) = \Sigma_{22}.$$

By central limit theorem, it follows that

$$\frac{1}{\sqrt{n}}\sum_{i\in S_j}\varepsilon_i d(\lambda)^{-1/2}N_{U_i}(z_0) \rightsquigarrow N(0,\sigma^2\Sigma_{22}). \tag{E.13}$$

As $s\to\infty$, we have $\frac{1}{\sqrt{n}}\sum_{i=1}^n \varepsilon_i s^{-1/2}d(\lambda)^{-1/2}N_{U_i}(z_0) \to 0$. So the second component in (III) is asymptotically ignorable. Therefore by (E.11), we obtain

$$(III) \rightsquigarrow N\left(0, \sigma^2\boldsymbol{x}^T\boldsymbol{\Omega}^{-1}\boldsymbol{x}\right).$$

As for (IV), we apply similar arguments as in the previous paragraph and consider $s\to\infty$. It follows that

$$(IV) = \sqrt{1-s^{-1}}\Big\{\frac{1}{\sqrt{N-n}}\sum_{i\notin S_j} d(\lambda)^{-1/2}N_{U_i}(z_0)\varepsilon_i\Big\} \rightsquigarrow N(0,\sigma^2\Sigma_{22}).$$

Lastly, note that $(III)$ and $(IV)$ are independent, so are their limits. Therefore, it follows that

$$(I) \rightsquigarrow N(0,\sigma^2(\boldsymbol{x}^T\boldsymbol{\Omega}^{-1}\boldsymbol{x}+\Sigma_{22})).$$

Case 2: $s$ fixed. Instead of decomposing (III) into two components as in previous case, we apply CLT to term (III) as a whole. Note that the summands in (III) are i.i.d. with mean zero. Moreover,

$$\begin{aligned}
&\mathbb{E}\Big[\varepsilon^2\big(\boldsymbol{x}^T L_U + s^{-1/2}d(\lambda)^{-1/2}N_U(z_0)\big)^2\Big] \\
&= \sigma^2\mathbb{E}\big[(\boldsymbol{x}^T L_U)^2\big] + s^{-1}\sigma^2\mathbb{E}\big[(d(\lambda)^{-1/2}N_U(z_0))^2\big] \\
&\qquad + 2s^{-1/2}\sigma^2\mathbb{E}\big[d(\lambda)^{-1/2}\boldsymbol{x}^T L_U N_U(z_0)\big]
\end{aligned}$$

The first two terms are considered in Case 1. For the third term, we have

$$\mathbb{E}\big[\boldsymbol{x}^T L_U d(\lambda)^{-1/2}N_U(z_0)\big] = \mathbb{E}\big[d(\lambda)^{-1/2}\boldsymbol{x}^T L_U\big(\widetilde{K}_Z(z_0) - L_U^T\boldsymbol{A}(z_0)\big)\big]$$

From Case 1, we have $d(\lambda)^{-1/2}\mathbb{E}\big[\widetilde{K}_Z(z_0)L_U\big] \to \boldsymbol{\alpha}_{z_0}$, $\mathbb{E}[L_U L_U^T] \to \boldsymbol{\Omega}^{-1}$, and

$\boldsymbol{A}(z_0)/d(\lambda)^{1/2} \to -\boldsymbol{\gamma}_{z_0}$. It follows that

$$\mathbb{E}\big[d(\lambda)^{-1/2}\boldsymbol{x}^T L_U N_U(z_0)\big] \to \boldsymbol{x}^T\boldsymbol{\Omega}^{-1}(\boldsymbol{\alpha}_{z_0} + \boldsymbol{\gamma}_{z_0}). \tag{E.14}$$

Therefore, we have

$$\mathbb{E}\Big[\varepsilon^2\big(\boldsymbol{x}^T L_U + s^{-1/2}d(\lambda)^{-1/2}N_U(z_0)\big)^2\Big] \to \sigma^2(\boldsymbol{x}^T\boldsymbol{\Omega}^{-1}\boldsymbol{x} + s^{-1}\Sigma_{22} + 2s^{-1/2}\boldsymbol{x}^T\boldsymbol{\Sigma}_{12}).$$

Hence by central limit theorem, we have

$$(III) \rightsquigarrow N(0, \sigma^2(\boldsymbol{x}^T\boldsymbol{\Omega}^{-1}\boldsymbol{x} + s^{-1}\Sigma_{22} + 2s^{-1/2}\boldsymbol{x}^T\boldsymbol{\Sigma}_{12})).$$

Similarly, we have

$$(IV) \rightsquigarrow N\big(0, (1-s^{-1})\sigma^2\Sigma_{22}\big).$$

As (III) and (IV) are independent, so are their limits. Therefore in the case that $s$ is fixed, we have

$$(I) \rightsquigarrow N(0, \sigma^2(\boldsymbol{x}^T\boldsymbol{\Omega}^{-1}\boldsymbol{x} + \Sigma_{22} + 2s^{-1/2}\boldsymbol{x}^T\boldsymbol{\Sigma}_{12})).$$

This finishes the proof. □

### E.3. Proof of Lemma 7.4.

Proof. Recall that $\Delta m^{(j)} = \widehat{m}^{(j)} - m_0^{(j)}$. As $\Delta m^{(j)}$ minimizes the objective function (3.2), we have

$$\frac{1}{n}\sum_{i\in S_j}(\widehat{m}^{(j)}(U_i) - Y_i)^2 + \lambda\|\widehat{f}\|_{\mathcal{H}}^2 \le \frac{1}{n}\sum_{i\in S_j}(m_0(U_i) - Y_i)^2 + \lambda\|f_0\|_{\mathcal{H}}^2,$$

On the $j$-th sub-population, we have $Y_i = m_0^{(j)}(U_i) + \varepsilon_i$, hence it follows that

$$\frac{1}{n}\sum_{i\in S_j}(\widehat{m}^{(j)}(U_i) - m_0^{(j)}(U_i))^2 + \frac{2}{n}\sum_{i\in S_j}\varepsilon_i(\widehat{m}^{(j)}(U_i) - m_0^{(j)}(U_i)) + \lambda\|\widehat{f}^{(j)}\|_{\mathcal{H}}^2 \le \lambda\|f_0\|_{\mathcal{H}}^2.$$

Adding and subtracting $\mathbb{E}_U[\Delta m^{(j)}(U)^2]$, we transform the above inequality to

$$\frac{1}{n}\sum_{i\in S_j}\Delta m^{(j)}(U_i)^2 - \mathbb{E}_U[\Delta m^{(j)}(U)^2] + \mathbb{E}_U[\Delta m^{(j)}(U)^2] + \lambda\|\Delta f^{(j)}\|_{\mathcal{H}}^2$$
$$+ \frac{2}{n}\sum_{i\in S_j}\varepsilon_i\Delta m^{(j)}(U_i) - 2\lambda\|f_0\|_{\mathcal{H}}^2 + 2\lambda\langle\widehat{f}^{(j)}, f_0\rangle_{\mathcal{H}} \le 0.$$

As we have $\mathbb{E}_U[\Delta m^{(j)}(U)^2] + \lambda\|\Delta f^{(j)}\|_{\mathcal{H}}^2 = \|\Delta m^{(j)}\|_{L_2(\mathbb{P}_U)}^2 + \lambda\|\Delta f^{(j)}\|_{\mathcal{H}}^2 =$

$\|\Delta m^{(j)}\|_{\mathcal{A}}^2$. It follows that

$$\begin{aligned}
\|\Delta m^{(j)}\|_{\mathcal{A}}^2 &\leq -2\Big(\frac{1}{n}\sum_{i\in S_j}\varepsilon_i \Delta m^{(j)}(U_i) - \lambda\langle \Delta f^{(j)}, f_0\rangle_{\mathcal{H}}\Big) \\
&\qquad -\frac{1}{n}\sum_{i\in S_j}\langle \Delta m^{(j)}(U_i)R_{U_i} - \mathbb{E}_U[\Delta m^{(j)}(U)R_U], \Delta m^{(j)}\rangle_{\mathcal{A}} \\
&= -2\langle\frac{1}{n}\sum_{i\in S_j}\varepsilon_i R_{U_i} - P_\lambda m_0^{(j)}, \Delta m^{(j)}\rangle_{\mathcal{A}} - \langle Rem^{(j)}, \Delta m^{(j)}\rangle_{\mathcal{A}}. \qquad \text{(E.15)}
\end{aligned}$$

Define the following two events:

$$\mathcal{B}_1 := \Big\{\Big\|\frac{1}{n}\sum_{i\in S_j}\varepsilon_i R_{U_i}\Big\|_{\mathcal{A}} \leq C\log^2 n(d(\lambda)/n)^{1/2}\Big\},$$

$$\mathcal{B}_2 := \big\{\|Rem^{(j)}\|_{\mathcal{A}} \leq 2c_r d(\lambda)n^{-1/2}\big(CJ(\mathcal{F},1)+\log n\big)\|\Delta m^{(j)}\|_{\mathcal{A}}\big\}.$$

We bound the two terms in (E.15) respectively. First, note that

$$\begin{aligned}
\|P_\lambda m_0^{(j)}\| &= \sup_{\|m\|_{\mathcal{A}}=1}|\langle P_\lambda m_0^{(j)}, m\rangle_{\mathcal{A}}| = \sup_{\|m\|_{\mathcal{A}}=1}\lambda|\langle f_0, f\rangle_{\mathcal{H}}| \\
&\leq \sup_{\|m\|=1}\sqrt{\lambda\|f_0\|_{\mathcal{H}}^2}\sqrt{\lambda\|f\|_{\mathcal{H}}^2} \leq \lambda^{1/2}\|f_0\|_{\mathcal{H}},
\end{aligned}$$

where the last inequality follows from the fact that $\lambda\|f\|_{\mathcal{H}}^2 \leq \|m\|_{\mathcal{A}}^2 = 1$. Therefore on event $\mathcal{B}_1$, the first term in (E.15) can be bounded by

$$\begin{aligned}
\Big|\langle\frac{1}{n}\sum_{i\in S_j}\varepsilon_i R_{U_i} - P_\lambda f_0, \Delta m^{(j)}\rangle_{\mathcal{A}}\Big| &\leq \Big\|\frac{1}{n}\sum_{i\in S_j}\varepsilon_i R_{U_i} - P_\lambda f_0\Big\|_{\mathcal{A}}\|\Delta m^{(j)}\|_{\mathcal{A}} \\
&\leq C\big(\log^2 n(d(\lambda)/n)^{1/2} + \lambda^{1/2}\big)\|\Delta m^{(j)}\|_{\mathcal{A}}. \qquad \text{(E.16)}
\end{aligned}$$

Furthermore, on the event $\mathcal{B}_2$, the second term in (E.15) can be bounded by

$$\begin{aligned}
\Big|\frac{1}{n}\sum_{i\in S_j}\langle Rem^{(j)}, \Delta m^{(j)}\rangle_{\mathcal{A}}\Big| &\leq \|\Delta m^{(j)}\|_{\mathcal{A}}\|Rem^{(j)}\|_{\mathcal{A}} \\
&\leq 2c_r d(\lambda)n^{-1/2}\big(CJ(\mathcal{F},1)+\log n\big)\|\Delta m^{(j)}\|_{\mathcal{A}}^2. \qquad \text{(E.17)}
\end{aligned}$$

Therefore, by (E.15), (E.16) and (E.17), it yields that on the event $\mathcal{B}_1 \cap \mathcal{B}_2$ there exist constants $C, C'$,

$$\begin{aligned}
\|\Delta m^{(j)}\|_{\mathcal{A}}^2 &\leq C'\big(\sqrt{d(\lambda)/n}\log^2 n + \lambda^{1/2}\big)\|\Delta m^{(j)}\|_{\mathcal{A}} \\
&\qquad + 2c_r d(\lambda)n^{-1/2}\big(CJ(\mathcal{F},1)+\log n\big)\|\Delta m^{(j)}\|_{\mathcal{A}}^2.
\end{aligned}$$

By the conditions in the statement of the lemma, we have $d(\lambda)n^{-1/2}\big(CJ(\mathcal{F},1)+$

$\log n\big) = o(1)$. Therefore it follows that

$$\|\Delta m^{(j)}\|_{\mathcal{A}} \le C_1\big(\sqrt{d(\lambda)/n}\log^2 n + \lambda^{1/2}\big) + o(1)\|\Delta m^{(j)}\|_{\mathcal{A}}$$

which implies that for sufficiently large $n$,

$$\|\Delta m^{(j)}\|_{\mathcal{A}} \le C(\log^2 n\sqrt{d(\lambda)/n} + \lambda^{1/2}).$$

Now we are left to bound the probability of $\mathcal{B}_1^c \cup \mathcal{B}_2^c$. For $\mathcal{B}_1$, define $\mathcal{Q}_i = \{|\varepsilon_i| \le \log n\}$. Since $\|R_U\| \le c_r d(\lambda)^{1/2}$, we have that on the event of $\cap_{i\in S_j}\mathcal{Q}_i$, $\big\{\varepsilon_i R_{U_i}\big\}_{i\in S_j}$ is a sequence of random variables in Hilbert space $\mathcal{A}$ that are independent with mean zero and bounded by $cd(\lambda)^{1/2}\log n$. Therefore we have

$$\begin{aligned}
\mathbb{P}(\mathcal{B}_1^c) &= \mathbb{P}\Big(\|\frac{1}{n}\sum_{i\in S_j}\varepsilon_i R_{U_i}\|_{\mathcal{A}} > C\log^2 n(d(\lambda)/n)^{1/2}\Big)\\
&\le \mathbb{P}\Big(\cap_{i\in S_j}\mathcal{Q}_i, \|\frac{1}{n}\sum_{i\in S_j}\varepsilon_i R_{U_i}\|_{\mathcal{A}} > C\log^2 n(d(\lambda)/n)^{1/2}\Big) + \mathbb{P}\big((\cap_i\mathcal{Q}_i)^c\big)\\
&\le 2\exp(-\log^2 n) + 2n\exp(-\log^2 n), \qquad \text{(E.18)}
\end{aligned}$$

where the first term in the last inequality is by Lemma G.1, and the second term is by union bound and the fact that $\varepsilon_i$ are i.i.d. sub-Gaussian.

Now we turn to $\mathcal{B}_2$. Define $\widetilde{m} := (2c_r)^{-1}d(\lambda)^{-1/2}\frac{\Delta m^{(j)}}{\|\Delta m^{(j)}\|_{\mathcal{A}}}$. Then it follows that

$$\|\widetilde{m}\|_{\sup} \le c_r d(\lambda)^{1/2}\|\widetilde{m}\|_{\mathcal{A}} \le 1/2.$$

By the same argument as in Section 7.4, it follows that $\|\Delta f^{(j)}\|_{\sup} \le 1/2$ and $|\boldsymbol{x}^T\Delta\boldsymbol{\beta}| \le 1$ for all $\boldsymbol{x}$. Moreover, we have

$$\|\widetilde{f}\|_{\mathcal{H}} \le \lambda^{-1/2}\|\widetilde{m}\|_{\mathcal{A}} \le (2c_r)^{-1}d(\lambda)^{-1/2}\lambda^{-1/2}.$$

Hence we proved that $\widetilde{m} \in \mathcal{F}$. By Lemma F.1, it follows that

$$\mathbb{P}\Big(\|Z_n(\widetilde{m})\|_{\mathcal{A}} \ge CJ(\mathcal{F},1) + \log n\Big) \le C\exp\Big(-\log^2 n/C\Big).$$

By definition of $Z_n(m)$ and $\widetilde{m}$, we have $Z_n(\widetilde{m}) = (2c_r)^{-1}d(\lambda)^{-1}n^{1/2}\|\Delta m^{(j)}\|_{\mathcal{A}}^{-1}Rem^{(j)}$. Hence it follows that

$$\begin{aligned}
\mathbb{P}(\mathcal{B}_2^c) &= \mathbb{P}\Big(\|Rem^{(j)}\|_{\mathcal{A}} \ge 2c_r d(\lambda)n^{-1/2}\big(CJ(\mathcal{F},1) + \log n\big)\|\Delta m^{(j)}\|_{\mathcal{A}}\Big)\\
&\le C\exp\Big(-\log^2 n/C\Big). \qquad \text{(E.19)}
\end{aligned}$$

Combining (E.18) and (E.19), we have that for some universal constants $c, C$

and sufficiently large $n$,

$$\mathbb{P}\Big(\|\Delta m^{(j)}\|_{\mathcal{A}} \geq C\big(\sqrt{d(\lambda)/n}\log^2 n + \lambda^{1/2}\big)\Big) \leq \mathbb{P}(\mathcal{B}_1^c) + \mathbb{P}(\mathcal{B}_2^c)$$

$$\lesssim n\exp(-c\log^2 n). \tag{E.20}$$

This completes the proof. □

## APPENDIX F: AN EMPIRICAL PROCESS LEMMA

**Lemma F.1** (Local chaining inequality)**.** If $Z_n(m) \in \mathcal{A}$ is a separable process on the metric space $(\mathcal{F}, \|\cdot\|_{\sup})$ and satisfies (7.15):

$$\mathbb{P}\Big(\|Z_n(m_1) - Z_n(m_2)\|_{\mathcal{A}} \geq t\Big) \leq 2\exp\Big(-\frac{t^2}{8\|m_1 - m_2\|_{\sup}^2}\Big)$$

Then for all $m_0 \in \mathcal{F}$ and $x \geq 0$, we have

$$\mathbb{P}\Big(\sup_{m\in\mathcal{F}}\|Z_n(m) - Z_n(m_0)\|_{\mathcal{A}} \geq CJ(\mathcal{F}, \mathrm{diam}(\mathcal{F})) + x\Big) \leq C\exp\big(-x^2/C\mathrm{diam}(\mathcal{F})^2\big), \tag{F.1}$$

where $C$ is a generic constant.

Proof. The proof follows by modifying the proof of Theorem 5.28 in van Handel (2014). Let $k_0$ be the largest integer such that $2^{-k_0} \geq \mathrm{diam}(\mathcal{F})$. Then $N(\mathcal{F}, d, 2^{-k}) = 1$ for all $k \leq k_0$. We employ a chaining argument, and start at the scale $2^{-k_0}$. For every $k > k_0$, let $N_k$ be a $2^{-k}$ net such that $|N_k| = N(T, d, 2^{-k})$. We define the singleton $N_{k_0} = \{m_0\}$. We claim that

$$\begin{aligned} Z_n(m) &= \lim_{k\to\infty} Z_n(\pi_k(m)) - Z_n(m_0) \\ &= \sum_{k>k_0}\big\{Z_n(\pi_k(m)) - Z_n(\pi_{k-1}(m))\big\} \quad \text{a.s.} \end{aligned}$$

where $\pi_k(m)$ is the closes point in $N_k$ to $m$. To prove this identity, note that the sub-Gaussian property of $\{Z_n(m)\}_{m\in\mathcal{F}}$ implies that $Z_n(m) - Z_n(\pi_k(m))$ is $d(m, \pi_k(m))$-sub-Gaussian. Thus

$$\sum_{k=k_0}^{\infty}\mathbb{E}\big[\|Z_n(m) - Z_n(\pi_k(m))\|_{\mathcal{A}}^2\big] \leq \sum_{k=k_0}^{\infty} d(m, \pi_k(m))^2 \leq \sum_{k=k_0}^{\infty} 2^{-k} < \infty.$$

It follows that $\|Z_n(m) - Z_n(\pi_k(m))\|_{\mathcal{A}} \to 0$ a.s. as $k \to \infty$, and the chaining identity follows readily using the telescoping property of the sum. By the chaining identity and separability of $\mathcal{F}$, we obtain

$$\sup_{m\in\mathcal{F}}\|Z_n(m) - Z_n(m_0)\|_{\mathcal{A}} \leq \sum_{k>k_0}\sup_{m\in\mathcal{F}}\|Z_n(\pi_k(m)) - Z_n(\pi_{k-1}(m))\|_{\mathcal{A}}.$$

By union bound and sub-Gaussian property, it follows that

$$\mathbb{P}\Big(\sup_{m\in\mathcal{F}}\|Z_n(\pi_k(m))-Z_n(\pi_{k-1}(m))\|_{\mathcal{A}}>t\Big)\le 2|N_k|\exp\Big(-\frac{t^2}{8\cdot 2^{-2k}}\Big)$$
$$=2\exp\Big(\log|N_k|-\frac{t^2}{8\cdot 2^{-2k}}\Big).$$

For large enough $t$, let $u=\frac{t^2}{8\cdot 2^{-2k}}-\log|N_k|$, and we have

$$\mathbb{P}\Big(\sup_{m\in\mathcal{F}}\|Z_n(\pi_k(m))-Z_n(\pi_{k-1}(m))\|_{\mathcal{A}}>2\sqrt{2}\cdot 2^{-k}(\sqrt{\log|N_k|}+u)\Big)\le 2\exp^{u^2/2},$$

which implies the link $\|Z_n(\pi_k(m))-Z_n(\pi_{k-1}(m))\|_{\mathcal{A}}$ at scale $k$ is small. To show that the links at all scale of $k$ are small simultaneously, we again use the union bound. Define the event $\mathcal{D}:=\big\{\exists k>k_0 \text{ s.t.} \sup_{m\in\mathcal{F}}\|Z_n(\pi_k(m))-Z_n(\pi_{k-1}(m))\|_{\mathcal{A}}>2\sqrt{2}\cdot 2^{-k}(\sqrt{\log|N_k|}+u_k)\big\}$, where $u_k=x+\sqrt{k-k_0}$. Then

$$\begin{aligned}&\mathbb{P}(\mathcal{D})\\&\le\sum_{k>k_0}\mathbb{P}\Big(\sup_{m\in\mathcal{F}}\|Z_n(\pi_k(m))-Z_n(\pi_{k-1}(m))\|_{\mathcal{A}}>2\sqrt{2}\cdot 2^{-k}(\sqrt{\log|N_k|}+u_k)\Big)\\&\le\sum_{k>k_0}\exp(-u_k^2/2)\le\exp(-x^2/2)\sum_{k>0}\exp(-k/2)\le C\exp(-x^2/2).\end{aligned}$$

Moreover, by the fact that $2^{-k_0}\le 2\mathrm{diam}(\mathcal{F})$ and

$$2^{-k_0}\le C2^{-k_0-1}\sqrt{\log N(\mathcal{F},d,2^{-k_0-1})}\le C\sum_{k>k_0}\sqrt{\log|N_k|},$$

we have on the event $\mathcal{D}^c$,

$$\begin{aligned}&\sup_{m\in\mathcal{F}}\|Z_n(m)-Z_n(m_0)\|_{\mathcal{A}}\\&\le\sum_{k>k_0}\sup_{m\in\mathcal{F}}\|Z_n(\pi_k(m))-Z_n(\pi_{k-1}(m))\|_{\mathcal{A}}\\&\le 2\sqrt{2}\sum_{k>k_0}2^{-k}(\sqrt{\log|N_k|}+u_k)\\&\le 2\sqrt{2}\sum_{k>k_0}2^{-k}\sqrt{\log|N_k|}+2\sqrt{2}\cdot 2^{-k_0}\sum_{k>k_0}2^{-k}\sqrt{k}+2\sqrt{2}\sum_{k>k_0}2^{-k}x\\&\le C\int_0^{\mathrm{diam}(\mathcal{F})}\sqrt{\log N(\mathcal{F},d,\epsilon)}d\epsilon+C\mathrm{diam}(\mathcal{F})x\\&=CJ(\mathcal{F},\mathrm{diam}(\mathcal{F}))+C\mathrm{diam}(\mathcal{F})x.\end{aligned}$$

Therefore

$$\mathbb{P}\Big(\sup_{m\in\mathcal{F}} \|Z_n(m) - Z_n(m_0)\|_{\mathcal{A}} \ge CJ(\mathcal{F}, \mathrm{diam}(\mathcal{F})) + C\mathrm{diam}(\mathcal{F})x\Big) \le \mathbb{P}(\mathcal{D}) \le C\exp(-x^2/2).$$

Replacing $C\mathrm{diam}(\mathcal{F})x$ with a new variable $x$, we reach the conclusion of the lemma. □

## APPENDIX G: AUXILIARY LEMMAS

**Lemma G.1.** (Pinelis, 1994) If $\Xi_1, \ldots, \Xi_s$ are zero mean independent random variables in a separable Hilbert space and $\|\Xi_i\| \le M$ for $i = 1, \ldots, n$, then

$$\mathbb{P}\Big(\Big\|\frac{1}{n}\sum_{i=1}^n \Xi_i\Big\| > t\Big) < 2\exp\big(-\frac{nt^2}{2M^2}\big).$$

**Lemma G.2.** We have for all $j = 1, \ldots, s$,

$$\mathbb{E}\big[\|\Delta f^{(j)}\|_{\mathcal{H}}^2 \,|\, \mathbb{U}\big] \le 2\sigma^2/\lambda + 4\|f_0\|_{\mathcal{H}}^2. \tag{G.1}$$

Proof. By the zero order optimality condition, we have

$$\begin{aligned}\lambda\|\widehat{f}^{(j)}\|_{\mathcal{H}}^2 &\le \frac{1}{n}\sum_{i\in S_j}(Y_i - \widehat{m}^{(j)}(U_i))^2 + \lambda\|\widehat{f}^{(j)}\|_{\mathcal{H}}^2\\ &\le \frac{1}{n}\sum_{i\in S_j}(Y_i - m_0^{(j)}(U_i))^2 + \lambda\|f_0\|_{\mathcal{H}}^2\\ &= \frac{1}{n}\sum_{i\in S_j}\varepsilon_i^2 + \lambda\|f_0\|_{\mathcal{H}}^2.\end{aligned}$$

Hence taking expectation conditioned on $\mathbb{U}$, we get

$$\lambda\mathbb{E}[\|\widehat{f}^{(j)}\|_{\mathcal{H}}^2 \,|\, \mathbb{U}] \le \sigma^2 + \lambda\|f_0\|_{\mathcal{H}}^2.$$

Then, applying triangular inequality along with the inequality $(a+b)^2 \le 2a^2 + 2b^2$, we have

$$\begin{aligned}\mathbb{E}[\|\Delta f^{(j)}\|_{\mathcal{H}}^2 \,|\, \mathbb{U}] &\le 2\|f_0\|_{\mathcal{H}}^2 + 2\mathbb{E}[\|\widehat{f}^{(j)}\|_{\mathcal{H}}^2 \,|\, \mathbb{U}]\\ &\le \frac{2\sigma^2}{\lambda} + 4\|f_0\|_{\mathcal{H}}^2,\end{aligned}$$

as desired. □

**Lemma G.3** (Matrix Heoffding in Tropp (2012))**.** Consider a finite sequence $\{\mathbf{A}_i\}_{i=1}^n$ of independent, random, symmetric matrices with dimension $p$. Assume that each random matrix satisfies

$$\mathbb{E}[\mathbf{A}_i] = \mathbf{0} \quad \text{and} \quad \|\mathbf{A}_i^2\| \le M \text{ almost surely.}$$

Then, for all $t > 0$,

$$\mathbb{P}\Big(\Big\|\frac{1}{n}\sum_{i=1}^{n}\mathbf{A}_i\Big\| \geq t\Big) \leq p\exp\Big(-\frac{nt^2}{8M}\Big).$$

**Lemma G.4.** Let $\mathbf{A}, \mathbf{E} \in \mathbb{R}^{k\times k}$ be given. If $\mathbf{A}$ is invertible, and $\|\mathbf{A}^{-1}\mathbf{E}\| < 1$, then $\widetilde{\mathbf{A}} := \mathbf{A} + \mathbf{E}$ is invertible, and

$$\|\widetilde{\mathbf{A}}^{-1} - \mathbf{A}^{-1}\| \leq \frac{\|\mathbf{E}\|\|\mathbf{A}^{-1}\|^2}{1 - \|\mathbf{A}^{-1}\mathbf{E}\|}$$

PROOF. See Theorem 2.5, p. 118 in Stewart and Sun (1990). □

Department of operations research
and financial engineering
Princeton University
Princeton, New Jersey 08544
USA
E-mail: tianqi@princeton.edu
hanliu@princeton.edu

Department of Statistics
Purdue University
West Lafayette, IN 47906
USA
E-mail: chengg@purdue.edu